\title{Towers of recollement and bases for diagram algebras:
planar diagrams and a little beyond}
\documentclass[12pt]{article}
\usepackage{latexsym}
\usepackage{amssymb}
\usepackage[all]{xy}
\UseComputerModernTips
\usepackage{graphicx}
\usepackage{epsf}
\usepackage[dvips]{rotating}
\providecommand{\noglossaryignore}[1]{}
\newcommand{\globalglossaryentry}[3]{\makebox[1.5in][l]{\tt $\backslash${#1}} 
\makebox[1.1in][l]{{$#2$}} \makebox[2.5in][l]{{#3}}\newline} 
\newcommand{\newcommandabbreviation}[3]{\newcommand{#1}{#2}%
\noglossaryignore{\globalglossaryentry{#3}{#2}{}}}
\newcommand{\renewcommandabbreviation}[3]{\renewcommand{#1}{#2}%
\noglossaryignore{\globalglossaryentry{#3}{#2}{}}}
\newcommand{\newcommandmacro}[4]{\newcommand{#1}{#2}%
\noglossaryignore{\globalglossaryentry{#3}{#2}{#4}}}
\newcommand{\gge}[3]{\noglossaryignore{\globalglossaryentry{#1}{#2}{#3}}}
\newcommand{\myaddress}%
{\parbox{3in}{\footnotesize \begin{center} 
Mathematics Department, City University, \\  
Northampton Square, London EC1V 0HB, UK.\end{center}}}
    
\newcounter{minidef}[section]

\newcounter{minicapt}

\newtheorem{theo}{Theorem}[subsection]       
\newtheorem{de}[theo]{Definition}     
\newtheorem{pr}[theo]{Proposition} 
\newtheorem{co}[theo]{Corollary}  
\newtheorem{rem}[theo]{Remark} 
\newtheorem{lem}[theo]{Lemma} 
\newtheorem{claim}{Claim}
\noglossaryignore{GREEK ETC.\newline}
\newcommandabbreviation{\e}{\epsilon}{e}        
\newcommandabbreviation{\lam}{\lambda}{lam}  
\newcommandabbreviation{\la}{\langle}{la}        
\newcommandabbreviation{\ran}{\rangle}{ran}
\newcommandabbreviation{\ha}{\#}{ha}             
\newcommandabbreviation{\rmap}{\rightarrow}{rmap}
\newcommandabbreviation{\aaa}{\alpha}{aaa}        
\newcommandabbreviation{\ab}{\alpha,\beta}{ab}
\newcommandabbreviation{\aab}{a(\ab )}{aab}       
\noglossaryignore{\newline RINGS\newline}
\newcommandabbreviation{\HH}{H \!\!\! I}{HH}              
\newcommandabbreviation{\C}{\mathbb C}{C}
\newcommandabbreviation{\N}{\mathbb N}{N}  
\newcommandabbreviation{\Z}{\mathbb Z}{Z}     
\renewcommandabbreviation{\Re}{\mathbb R}{Re}
\newcommandabbreviation{\R}{{\mathbb R}}{R}
\newcommandabbreviation{\Q}{\mathbb Q }{Q}
\renewcommandabbreviation{\H}{\mathbb H }{H}
\noglossaryignore{\newline SYMMETRIC GROUP\newline}
\def\Sym(#1){\Sigma(#1)}                  
\gge{Sym(-)}{\Sym(-)}{symmetric group on - objects}
\def\Sy(#1){\Sigma_{#1}}                  
\gge{Sy(-)}{\Sy(-)}{symmetric group irreducible -}
\def\sym(#1){\mbox{\LARGE s}(#1)}       
\gge{sym(-)}{\sym(-)}{symmetric group on - objects (variant)}
\def\sy(#1){\mbox{\LARGE s}({#1})}       
\gge{sy(-)}{\sy(-)}{symmetric group irreducible - (variant)}
\newcommandmacro{\cs}{\C \, \sy(n)}{cs}{symmetric group algebra over $\C$}
\noglossaryignore{\newline PARTITIONS/SETS\newline}
\newcommand{\Nset}[1]{\underline{#1}}
\gge{Nset\{-\}}{\Nset{-}}{set of natural numbers to -}
\def\nset(#1){ \{ #1 \}_{ \underline{n} }}
\gge{nset(-)}{\nset(-)}{a set $-\times\Nset$}
\def\ul(#1){_{\underline{#1}}}            
\gge{ul(-)}{{}\ul(-)}{subscript underline -}
\def\Ee(#1){{\bf E}_{#1}}                 
\gge{Ee(-)}{\Ee(-)}{set of equivalence relations on set -}
\def\Eee(#1){{\bf E}_{\{ #1 \}_{\underline{n}}}}  
\gge{Eee(-)}{\Eee(-)}{ditto for nset}
\def\Een(#1,#2){{\bf E}_{\{ #1 \}_{\underline{#2}}}}  
\def\Ssn(#1,#2){{\bf S}_{\{ #1 \}_{\underline{#2}}}}  
\def\Ss(#1){{\bf S}_{#1}}                 
\def\Sss(#1){{\bf S}_{\{ #1 \}_{\underline{n}}}}  
\def\bbc(#1){((\beta_1)(\beta_2)...(\beta_{#1}))}     
\newcommandmacro{\Ln}{{\Gamma}^{n}}{Ln}{large index set}
\newcommandmacro{\LnQ}{{\Gamma}^{n}_Q}{LnQ}{index set}
\newcommandmacro{\Zz}{\zeta}{Zz}{`shape' function}
        
\noglossaryignore{\newline PARTITION ALGEBRA\newline}
\def\ka(#1){\kappa_{#1}}                  
\def\Sm(#1){\Sigma_{#1}}                  
\newcommandmacro{\com}{\bullet}{com}{bullet composition}
\newcommandmacro{\enm}{\; e^n(\! m\! ) \;}{enm}{product of idempotents}
\def\Ai(#1){ A^{ #1 \cdot } }             
\def\Aij(#1,#2){ A^{ #1  #2 } }           
\newcommandmacro{\One}{\mbox{\bf $1 \!\!\! 1$}}{One}{algebra unit 1}

\newcommandmacro{\Bp}{B_p}{Bp}{partition basis}
\def\Bb(#1){B_p[#1]}                      
\def\Pp(#1){P_n[#1]}                      
\def\Ps(#1){P_n[#1] \! /}                 
\newcommandmacro{\Ph}{\hat{P}}{Ph}{P hat  algebra}
\def\Is(#1){\sim^{#1}}                    
\noglossaryignore{\newline MODULES\newline}
\def\Wm(#1){{\cal S}_{#1}}                
\gge{Wm(-)}{\Wm(-)}{Weyl module with index -}
\def\wm(#1,#2){{}_{#1}{\cal S}_{#2}}      
\gge{wm(-1,-)}{\wm(-1,-)}{Weyl module with index -}
\def\Ind(#1,#2,#3){\mbox{Ind}_{#1}^{#2}#3}
\gge{Ind(-1,-2,-)}{\Ind(-1,-2,-)}{induction}
\def\Res(#1,#2,#3){\mbox{Res}_{#1}^{#2}#3}
\gge{Res(-1,-2,-)}{\Res(-1,-2,-)}{restriction}
\newcommandabbreviation{\weyl}{standard}{weyl}
\newcommandabbreviation{\mod}{\mbox{mod}}{mod}
\newcommandabbreviation{\head}{\mbox{head }}{head}
\newcommandabbreviation{\Weyl}{Weyl}{Weyl}
\def\SS(#1){{\cal S}_{#1}}                
\gge{SS(-)}{\SS(-)}{Specht/Weyl module index -}
\def\LL(#1){{\cal L}_{#1}}                
\gge{LL(-)}{\LL(-)}{simple module index -}
\noglossaryignore{\newline FUNCTORS/MAPS\newline}
\newcommandmacro{\Gg}{{\cal G}}{Gg}{G Functor}
\newcommandmacro{\Fg}{{\cal F}}{Fg}{F Functor}
\newcommandmacro{\ra}{\rightarrow}{ra}{}
\def\ses(#1,#2,#3){0\ra #1 \ra #2 \ra #3 \ra 0}  
\gge{ses(1,2,-)}{\ses(1,2,-)}{\hspace{.5in} short exact sequence}
\def\starr(#1){ \stackrel{ #1 }{\longrightarrow} }
\gge{starr(-)}{\starr(-)}{}
\newcommandmacro{\doublerightarrow}{\; -\!\!\! -\!\!\!\!\!\! \gg \;}
{doublerightarrow}{}
\noglossaryignore{\newline PARTITION ALGEBRA MAPS\newline}
\newcommandmacro{\smap}{s}{smap}{`inclusion' map}
\newcommandmacro{\tmap}{t}{tmap}{$ P_n -> S_n$}
\newcommandmacro{\pmap}{\psi}{pmap}{$ S_n -> P_n $}
\noglossaryignore{\newline MISC.\newline}
\def\Amap(#1){{\cal A}_{#1}}              
\gge{Amap(-)}{\Amap(-)}{}
\def\Rr(#1){R_{#1}}                       
\gge{Rr(-)}{\Rr(-)}{restriction of E}
\def\Cr(#1){C_{#1}}                       
\gge{Cr(-)}{\Cr(-)}{restriction of E to N}
\newcommandmacro{\Tm}{{\cal T}}{Tm}{Transfer Matrix}
\def\On(#1){{\cal I}_{#1}}
\gge{On(-)}{\On(-)}{}
\newcommandmacro{\UU}{\underline{\sqcup}}{UU}{}  
\newcommandmacro{\UUU}{\sqcup}{UUU}{}  
\newcommandmacro{\Vq}{V_Q^{\otimes n}}{Vq}{Potts config. space}
\def\bs(#1,#2){\mbox{{\Large $\ast$}}^{#1}_{#2}} 
\gge{bs(-,-)}{\bs(-,-)}{general plumbing multiplier}
\newcommand{\ignore}[1]{}
\def\choo(#1,#2){ \left( \begin{array}{c} #1 \\ #2 \end{array} \right) }
\gge{choo(-1,-)}{\choo(-1,-)}{choose}
\newcommand{\Qed}{$\Box$}
\gge{Qed}{\mbox{\Qed}}{QED}
\def\staq(#1){\stackrel{#1}{=}}           
\gge{staq(-)}{\staq(-)}{}
\def\stam(#1){\stackrel{#1}{\rightarrow}} 
\gge{stam(-)}{\stam(-)}{}
\def\mat{ \left( \begin{array} }    
\def\tam{ \end{array}  \right) }
\gge{mat/tam}{...}{matrix delimiters}
\newcommand{\beq}{\begin{equation} }
\def\eql(#1){ \begin{equation} \label{#1} 
}
\newcommand{\eq}{\end{equation} }
\def\eqal(#1){\begin{eqnarray} \label{#1} }
\def\eqa{\end{eqnarray} }
\def\lab(#1){\label{#1}
}
\def\prl(#1){ \begin{pr} \label{#1} 
}
\def\del(#1){ \begin{de} \label{#1} 
}
   
\gge{smeq\{-\}}{...}{small equation}
  
\gge{fneq\{-\}}{...}{very small equation}
\noglossaryignore{\newline HECKE/BLOB\newline}
\newcommandmacro{\Hnq}{H_n(q)}{Hnq}{ * freestanding symbol}
\newcommandmacro{\Hn}{H_n}{Hn}{      *-mod etc.}
\newcommandmacro{\A}{{\cal A}}{A}{}
\newcommandmacro{\Cwts}{C}{Cwts}{}
\newcommandmacro{\CA}{{\cal A}}{CA}{}

\newcommandmacro{\calA}{{\cal A}}{calA}{}
\newcommandmacro{\modi}{\mbox{Mod} }{modi}{was mod not modi!}
\newcommandmacro{\Wgen}{{\Bbb S}}{Wgen}{}
\def\ol(#1){\overline{#1}}
\newcommandmacro{\st}{\mbox{St}}{st}{}
  
\def\CMult(#1,#2){(#1:#2)}
\def\CM(#1,#2){( #1 : #2 )}
\def\FMult#1,#2{(#1:#2)}
\def\CF#1,#2{(#1:#2)}

\newcommandmacro{\Top}{\mbox{Top}}{Top}{}
\newcommandmacro{\Soc}{\mbox{Soc}}{Soc}{}
\newcommandmacro{\Head}{\mbox{Head}}{Head}{}
\newcommandmacro{\Filt}{{\cal F}}{Filt}{}
\newcommandmacro{\Mod}{\mbox{mod}}{Mod}{}
\newcommandmacro{\Resi}{\mbox{Res }}{Resi}{was without i!}
\newcommandmacro{\Indi}{\mbox{Ind }}{Indi}{was without i!}
  
\def\RR(#1,#2){R^{#1}_{#2}}  
\def\TT(#1,#2){T^{#1}_{#2}}

\newcommandmacro{\Ann}{\mbox{Ann}}{Ann}{}
\newcommandmacro{\Cen}{\mbox{Cen}}{Cen}{}
\newcommandmacro{\End}{\mbox{End}}{End}{}
\newcommandabbreviation{\semisimple}{semisimple}{semisimple}
\newcommandabbreviation{\Bratteli}{Bratteli}{Bratteli}
\newcommandabbreviation{\JBC}{Jones Basic Construction}{JBC}
\newcommandabbreviation{\pa}{partition algebra}{pa}
\newcommandabbreviation{\TM}{transfer matrix}{TM}
\newcommandabbreviation{\PM}{Potts model}{PM}
\newcommandabbreviation{\QSC}{quantum spin chain}{QSC}
\newcommandabbreviation{\Hamiltonian}{Hamiltonian}{Hamiltonian}
\newcommandabbreviation{\YS}{Young symmetrizer}{YS}

\newcommand{\be}{\begin{eqnarray}}
\newcommand{\eeq}{\end{eqnarray}}
\newcommand{\non}{\nonumber}

\newcommand{\TL}{Temperley--Lieb}

\newcommand{\abe}{{\bf e}}

\newcommand{\Id}{{\mathbb I}}

\newcommand{\pdef}[1]{\medskip \noindent {\small {\bf #1 }} \newline }

\begin{document}
\newlength{\fred}
\setlength{\fred}{10pt}
\newcommand{\aar}{\ar@{-}}         
\newcommand{\headroom}{80}       
\newcommand{\raised}{-\headroom} 
\newcommand{\pp}[2]{\begin{picture}(#1,\headroom)(0,\raised)
    \put(-20,0){#2} \end{picture}}
\newcommand{\ppp}[4]{\begin{picture}(#1,#2)(0,-#3) \put(-20,0){#4}
    \end{picture}}       

\newcommand{\putch}{\put(-.3,-26)}
\newcommand{\poutch}{\put(-.3,-20)}
\newcommand{\hash}{\#}
\newcommand{\resp}{respectively}
\newcommand{\parker}[1]{\[ \framebox[7in]{ \mbox{\Large #1}} \]}
\newcommand{\aipic}{isotopic}
\newcommand{\aipy}{isotopy}
\newcommand{\pseudo}{ur}
\newcommand{\RX}{C}
\newcommand{\DD}{{D}}
\newcommand{\JJ}{{J}}
\newcommand{\Jn}[1]{{\JJ_{#1}}}
\newcommand{\DV}{{\DD(V)}}
\newcommand{\DVe}[2]{{\DD(V^{#1}_{#2})}}
\newcommand{\DVnm}{{\DD(V^n_m)}}
\newcommand{\DSVnm}{{\DD_S(V^n_m)}}
\newcommand{\Do}{{\DD^o}}
\newcommand{\DoV}{{\DD^o(V^{}_{})}}
\newcommand{\DoVe}[2]{{\DD^o(V^{#1}_{#2})}}
\newcommand{\DoVnm}{{\DD^o(V^n_m)}}

\newcommand{\Dz}[1]{{\DD^{z}_{#1}}}
\newcommand{\Dnz}{{\Dz{n}}}
\newcommand{\Dnzl}[2]{{\DD_{#1}^{z,#2}}}
\newcommand{\Dnp}{{\DD^p_{n}}}%
\newcommand{\Dnpc}{{\DD^{pc}_{n}}}%
\newcommand{\Dppc}[1]{{\DD^{pc'}_{#1}}}%
\newcommand{\Dph}{\DD^{\phi}}  
\newcommand{\Dzo}[1]{{\JJ^{z}_{#1 }}}
\newcommand{\Jp}[1]{{\JJ_{(#1)}}}
\newcommand{\Jpp}[1]{{\JJ'_{(#1)}}}
\newcommand{\DB}[1]{{\JJ^{B}_{#1}}}
\newcommand{\DBe}[1]{{\JJ^{Be}_{#1}}}
\newcommand{\Beta}{{\mathfrak B}}
\newcommand{\pseud}{{\mathcal H}}
\newcommand{\achira}{symplectic}
\newcommand{\achiral}{achiral}
\newcommand{\Achira}{Symplectic}
\newcommand{\Achiral}{Achiral}
\newcommand{\achiralb}{\achira\ blob}
\newcommand{\Achiralb}{\Achira\ blob}
\newcommand{\ASTL}{Affine symmetric TL}
\newcommand{\aSTL}{affine symmetric TL}
\newcommand{\CC}{CC}
\newcommand{\Deltag}{\Gamma} 
\newcommand{\sS}{\mathcal{S}}
\newcommand{\hati}{\tilde}
\newcommand{\V}{V} 

 \author{Paul Martin%
\footnote{Centre for Mathematical Science, City University, Northampton
  Square, London EC1V 0HB, UK},
R. M. Green
\footnote{Department of Mathematics, University of Colorado,
Campus Box 395, Boulder, Colorado 80309-0395, USA}, 
and Alison Parker%
\footnote{Department of Mathematics, University of Leicester,
  Leicester, LE1 7RH, UK}.
}
\date{}
\maketitle
\vspace{5cm} \begin{abstract}
The recollement approach to the representation
theory of sequences of algebras is extended to pass basis information 
directly through the globalisation functor. 
The method is hence adapted to treat sequences that are not
necessarily towers by inclusion, such as \achiralb\ algebras
(diagram algebra quotients of the type-$\hati{C}$ Hecke algebras).

By carefully reviewing the diagram algebra construction,
we find a new set of functors interrelating 
module categories of ordinary  blob algebras
(diagram algebra quotients of the type-${B}$ Hecke algebras)
at {\em different} values of the algebra parameters.
We show that these functors generalise
to determine the structure of 
\achiralb\ algebras, and hence of certain
two-boundary Temperley-Lieb algebras
arising in Statistical Mechanics. 

We identify the diagram basis with a cellular basis for each
symplectic blob algebra, and prove that these algebras are
quasihereditary over a field for almost all parameter choices, 
and generically semisimple. 
(That is, we give bases for all cell and standard modules.)
\end{abstract} \newpage
{\small \tableofcontents}
\section{Introduction}
\newcommand{\Ring}{K} 

The idea of recollement \cite{BeilinsonBernsteinDeligne81} 
is applied to categories of modules in
\cite{ClineParshallScott88}. 
Iterated `towers' of recollement are used in 
algebraic representation theory in \cite{Martin91} and
formalised, for example, in \cite{CoxMartinParkerXi03}. 
(The {\em tower} here refers to the algebraic structure needed for
statistical mechanics \cite{Martin91}, although an elementary connection can be made
in the semisimple case to Jones basic construction \cite{GoodmanDelaharpeJones89}.)

If $A$ is an algebra, and 
$e \in A$ an idempotent, then the category $eAe$-mod 
of left $eAe$-modules fully embeds in $A$-mod. 
At its most basic the idea is that if $eAe$-mod
may be relatively simply analysed, the embedding then gives partial
knowledge of $A$-mod. 
The standard tower picture has $A$ as
part of a tower of algebras by inclusion, such that $eAe$ may be identified
isomorphically with one of the subalgebras. 
The interplay of induction/restriction and
globalisation/localisation functors facilitates representation
theory in such a tower. 
This is the device discussed in \cite{CoxMartinParkerXi03}.


All this would be of academic interest only, were it not for the 
ubiquity of such towers `in nature'. 
Transfer matrix algebras are algebras whose
representation matrices build statistical mechanical transfer
matrices \cite{Martin91}. The stability of the thermodynamic limit corresponds 
(to cut a long story short) to the
existence of a tower of module-categorical embeddings. 
However \cite{CoxMartinParkerXi03} addresses only one way 
among many in which 
such a tower could occur. 

A further limitation 
in the formalisation of  \cite{CoxMartinParkerXi03} is
that it concentrates on the abstract module category tower, and does
not incorporate the tower of special module bases found in concrete
examples (such as in \cite{MartinSaleur94a,Martin96},
and cf. \cite{GrahamLehrer96}). 
This paper 
describes two ways in which the framework formalised in 
\cite{CoxMartinParkerXi03}
may be extended, so as to treat 
the representation theory of a wider class of algebras. 
Firstly we integrate the module category tower framework 
with special algebra bases
(as in diagram algebras, for example).
This allows us to enumerate explicit bases for modules and
   algebras, rather than simply to generate structure theorems.
Secondly we show that the framework is useful even when the tower is
not (necessarily) a tower of inclusion. 
Indeed the control of basis compensates for the lack of induction and
restriction functors, so the framework will work for algebra towers
without induction and restriction.
This latter is important for treating 
our motivating examples: 
families of algebras arising
recently in the Physics of systems with special boundaries \cite{DeGier02}, 
which do not include
(generalising the ordinary Temperley-Lieb
algebras \cite{TemperleyLieb71}
and their immediate family, which do include \cite{Martin91}). 
We demonstrate the method by determining the structure of these algebras.
(In the process we also introduce and make use of functors relating
module categories for algebras differing by the choice of 
specialisation of a deformation parameter --- certain choices of 
specialisation then being treatable together in `meta-categories'.)

\bigskip

\noindent
One may
deform the ordinary \TL\ diagram algebra by two-colouring the diagrams
(seen as maps --- see later) and assigning different parameters to  
shaded and unshaded loops. 
The result is isomorphic to the original. 
However deforming the B-type (left-right symmetric) subalgebra similarly 
constructs an algebra with a true extra parameter --- the blob algebra. 
Varying the extra parameter, the blob algebra may be used to build 
the representations of the periodic \TL\ algebra 
(\cite{MartinSaleur93,MartinSaleur94a,FanGreen99,Green98,GrahamLehrer03}).
Thus all these algebras can be analysed using the included-tower
technology \cite{CoxMartinParkerXi03}.
The next natural generalisation is the \achiralb\  
algebra
(left-right symmetric {\em and} periodic). 
This sequence of algebras cannot be made to include in an appropriate way, 
and so presents a suitable challenge for our method.


In this paper we implement a {\em towers of recollement} programme,
and variations, to determine the structure 
(bases and representation theory) of three interesting algebras. 
Many workers have considered wreath-like extensions for Brauer,
Temperley--Lieb diagram 
\cite{JonesPlanar,MartinSaleur94a,Green98b,FanGreen99} and even
partition algebras \cite{Bloss04}. 
These extensions 
are of interest
as testing grounds for techniques intended to be applied in the representation
theory of more classical objects, such as the symmetric group.
For example, one approach to systems of algebras with
Jucys--Murphy-like elements \cite{Jucys74,Murphy81,Cherednik91} 
is to consider extensions of the
algebras by new commuting generators which obey relations emulating
identities obeyed by the Murphy elements.
With careful preparation these extensions behave like wreaths 
(confer \cite{Nazarov96,HaringOldenburg99,CoxMartinParkerXi03}). 
Here, considering the most general case of nonabelian wreath algebras,
we also explore a new and intriguing set of interconnections,
taking us to the \achiralb\ algebra.  
This leads in particular to applications in boundary integrable
statistical mechanics \cite{DeGier02}. 


Planar diagram algebras (such as ordinary \TL) have been much studied
(with useful consequences in both representation theory and physics),
as have `fully non-planar' algebras such as the partition algebra.
These represent simplifying extremes in a range of generally harder
problems.  
As mentioned above, 
certain annular algebras can be brought into the
planar framework, using the blob algebra \cite{MartinSaleur94a}. 
The algebras we focus on
here are (in a suitable sense --- see later) mildly non-planar diagram
algebras. 
These are harder to treat, but not intractable, as we shall
demonstrate. 


The motivating aims of the paper are: 
\\
1. To provide an organisational framework for unifying the 
representation theory of various
forms of periodic/annular/type-B/boundary TL algebras 
studied in the literature 
\cite{MartinSaleur93,MartinSaleur94a,tomDieck94,DeGier02}.
(The ordinary \TL\ algebra is a nexus for many branches of mathematics,
with isomorphic algebras constructed in areas such as:
factors \cite{JonesPlanar}, 
representation theory, combinatorics, 
statistical mechanics \cite{TemperleyLieb71,Baxter82}.
Variants appear in these contexts, but are no longer all isomorphic,
and connections between them are not yet fully understood.)
\\
2. To provide the representation theoretic formalism for analysing these
algebras.


The layout of the paper is as follows.
In section~\ref{cat} we introduce the general theory necessary to
pass specific basis elements, for special types of module, 
between layers of a globalisation tower
(irrespective of inclusion). 
We also discuss how certain other features 
of modules which we shall use later
behave under globalisation. 
\\
In section~\ref{s diag} we collect the definitions of 
one of the families of diagram algebras which we shall need. 
All of these are based on Weyl's diagrams
for the Brauer algebra \cite{Weyl46}.  
To help prepare the ground for later more elaborate constructions
we also point out some paradigms for diagram algebra construction. 
For example:
combinatorial sets with diagrammatic realisations 
(which form bases for algebras via diagram concatenation), 
containing topologically characterised subsets respected by concatenation;
with the resultant subalgebras amenable to deformations not tolerated
by the original combinatorial algebra. 
\\
In section~\ref{blob'} we focus on deformations of Temperley--Lieb
algebras --- again looking at subalgebras and deformations.
We use these to establish homomorphisms between various families of
algebras. 
\\
In section~\ref{recol b} we use the alternate realisations established above
to construct new functors in the tower of blob algebras, and hence to
relate categories of modules for blob algebras with different
values of the defining parameters (and hence to
analyse their representation theory). 
\\
In sections~\ref{ASTLA} and~\ref{ASTLA'} 
we define affine symmetric  Temperley--Lieb
algebras and \achiralb\ algebras, and relate their categories of
modules. 
\\
In section~\ref{rep ASTLA} we investigate the representation theory of affine
symmetric  Temperley--Lieb algebras, using results from all the
previous sections. The main results here are perhaps the simple indexing
Theorem~\ref{index wangy}, and the generic structure results of
section~\ref{wa}. Beyond the semisimple cases, we show that the algebra
is almost always quasihereditary, and give bases for the standard modules. 


Several of the incarnations of the \TL\ algebra have 
an interesting `periodic' generalisation, as noted above, 
but these are much harder to treat. 
The blob algebra is a device that largely solves this problem
--- casting the representation theory of (infinite) periodic \TL\ 
algebras into the
setting of a (finite) generalisation of TL with properties much closer
to the original. 
The blob algebra suggests various generalisations of its own, but
these
are once again rather harder to treat, and have (until now) lacked the
motivations of the blob algebra (i.e. 
its simple but beautiful representation theory;
its application to periodic \TL\ and hence affine Hecke algebras). 
Recently, though, both the blob algebra and its two-boundary \TL\
algebra generalisation have arisen in the treatment of boundary
integrable systems in Statistical Mechanics 
\cite{DeGier02,NicholsRittenbergdeGier05},
suggesting that category embedding methods should work here.

\subsection{Preliminary definitions}\label{INot}

A Coxeter graph is any finite undirected
graph without loops (that is, without edges that begin and end on the
same vertex).
For example:
\[
A_n: \;\; \raisebox{-.1in}{\includegraphics{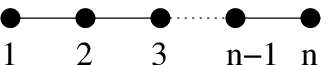}}
\qquad\qquad
B_n: \;\; \raisebox{-.1in}{\includegraphics{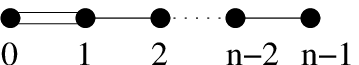}}
\]
The Coxeter Artin system of Coxeter graph $G$ is a pair $(B,S)$
consisting of a group $B$ and a set of pairs 
$g_v, g_v^{-1}$ of generators of $B$ labelled by
the vertices $v$ of $G$, with relations
\[
g_v g_{v'} ... \; = \; g_{v'} g_v ...
\] 
where the number of factors on each side is two more than the number
of edges between $v$ and $v'$.

The Coxeter system of $G$ is a pair $(W,S)$ where $W=W(G)$ is the quotient
of $B$ by the relation $g_v = g_v^{-1}$. 

For $K$ a ring, let $q_v$ be an invertible element in $K$, for each 
$v \in G$, such that $q_v = q_{v'}$ if $g_v$ conjugate to $g_{v'}$ in $B$. 
Let $Q_v = (g_v - q_v)(g_v + q_v^{-1})$.
The generic Hecke algebra of $G$ is 
\[
H(G) = K B / \langle Q_v \; | \; v \in G \rangle
\]
Examples: If $G=A_n$ then all generators are conjugate and we have the
(one-parameter) Hecke algebras of type-$A$.
If $G=B_n$ there is one generator not conjugate to the rest  and we have the
(two-parameter) Hecke algebras of type-$B$.
If $G=\hati{C_n}$ 
there are two generators not conjugate to the rest,
or to each other,  and we have the
(three-parameter) Hecke algebras of type-$\hati{C}$.


Each of these algebras has an algebra homomorphism onto $K$ defined by 
\[
\rho_+ (g_v) = q_v 
\]
Note that for $q_v=1$ the relation $Q_v=0$ is $g_v^2=1$ and hence a
group relation, so that $H(G)$ is the  group algebra of $W(G)$ in this case. 
It will be convenient to write $\sigma_v$ for $g_v$ in the group
case. We have \cite[Ch.7]{Humphreys90} that if 
$w= \sigma_{i_1} ...\sigma_{i_l}$ 
is a reduced expression in $W(G)$ then 
$\{ T_w = g_{i_1} ...g_{i_l} \; | \; w \in W(G) \}$ is a basis for $H(G)$ in general. 
Define the {\em symmetrizer} 
\[
E_G = \sum_{w \in W(G)} \rho_+(w) T_w
\]
in $H(G)$. 
For example if the vertices of $A_2$ are labelled from $\{1,2\}$ we
have
\[
E_{A_2} = 1+ q (g_1 + g_2) + q^2 (g_1g_2 + g_2 g_1) +q^3 g_1 g_2 g_1 . 
\]
Defining $u_v = (g_v + q_v^{-1})$, so that $g_v u_v = q_v u_v$, we have 
\[
E_{A_2}
= q^3 ( u_1 u_2 u_1 - u_1 ) . 
\]
If the vertices of $B_2$ are labelled from $\{ 0,1 \}$ we have 
\[
E_{B_2} = 1 + q_1 g_1 +q_0 g_0 + q_0 q_1 (g_0 g_1 + g_1 g_0 )
+q_0 q_1 (q_0 g_0 g_1 g_0 + q_1 g_1 g_0 g_1 ) +(q_0q_1)^2 g_0 g_1 g_0 g_1 
\]
\eql(uuuu)
=  u_0 u_1 u_0 u_1 - \frac{q_0^2 + q_1^2}{q_0 q_1} u_0 u_1 . 
\eq

Note that the relations of $H(G)$ are invariant under the parameter transformation 
\[
s_v : q_v \mapsto -q_v^{-1}  
\qquad \mbox{ and } 
s_v : q_w \mapsto q_w \qquad  \mbox{ for $g_v,g_w$ not in the same class}
\]
(that is, there is one such transformation for each parameter).
For each parameter transformation $s_w$ there are, in addition to $\rho_+$, further
homomorphisms of $H(G)$ onto $K$:  
\[
\rho_w (g_v) = s_w (q_v) . 
\]


Any subset of $S$ generates a parabolic subalgebra of $B$ or of $H(G)$. 
If  $v,v' \in G$ have at least one edge between them define
$E_{vv'}$ as the symmetrizer on their parabolic subalgebra of $H(G)$;
else $E_{vv'} =0$. 
Then define
\[
T(G) = H(G)/  \langle E_{v,v'} \; | \; v,v' \in G \rangle . 
\]
Example: $T(A_n)$ is the ordinary Temperley--Lieb algebra 
\cite{TemperleyLieb71}. 

Not much is known about $H(G)$ or $T(G)$ for general $G$ (although see 
\cite[Ch.9]{Martin91}), but the cases in which $G$ is positive
definite or positive semidefinite are relatively tractable (although
still interesting) 
\cite{MartinSaleur94a,GrahamLehrer98,FanGreen99,ErdmannGreen99,Green03}.

\medskip

\newcommand{\vardelta}{\delta_e}
\newcommand{\vargamma}{\gamma}%
\newcommand{\uU}{U}
\newcommand{\ue}{e}

For $K$ a ring, $x$ an invertible element in $K$, 
$q=x^2$, and $\vargamma,\vardelta \in K$, 
define $TLb_n^K$ to be the $K$--algebra with 
generators $\{ 1,e,\uU_1,\ldots,\uU_{n-1} \}$ and relations
\begin{eqnarray} 
 \uU_i \uU_i &=& (q+q^{-1}) \uU_i     \label{TL001} \\
 \uU_i \uU_{i\pm 1} \uU_i &=& \uU_i    \label{TL002} \\
 \uU_i \uU_j &=& \uU_j \uU_i \hspace{1in} \mbox{($|i-j|\neq 1$)}  \label{TL003}
\end{eqnarray}
\begin{eqnarray} 
 \uU_1 e \uU_1 &=&  \vargamma \uU_1   \label{TL004} \\
 e e &=& \vardelta e \\
 \uU_i e &=& e \uU_i   \hspace{1in} \mbox{($i \neq 1$)} . \label{TL006}
\end{eqnarray}


Note that $e$ can be rescaled to change $\vargamma$ and
$\vardelta$ by the same factor. 
Thus, if we require that $\vardelta$ be invertible, then we might as well
replace it by 1. 
This brings us to the original two--parameter presentation of the
algebra,
sometimes known as 
the blob algebra by presentation, or 
the two--parameter \TL\ algebra of type B.
The subalgebra generated by  $\{ 1,\uU_1,\ldots,\uU_{n-1} \}$
is $T(A_n)$, the \TL\ algebra of type A,
sometimes denoted $TL_n^{K}$.

For $k$ a field that is a $\Ring$--algebra define 
$TLb_n = k \otimes_{\Ring} TLb_n^{\Ring}$. 


The type-A algebra is isomorphic to the well known 
ordinary \TL\ diagram algebra \cite{Martin91,Graham95}; 
and the algebra $TLb_n$ is isomorphic to the blob diagram algebra $b_n$
\cite{MartinSaleur94a,Graham95,Green98,FanGreen99}
(see section~\ref{dir homs}).
Because of these isomorphisms
it is common to use generators and diagrams interchangeably.



There is also 
a {\em periodic} \TL\ diagram algebra (TLDA). 
That is, a TLDA 
defined using certain periodic TL diagrams.
Continuing the above `duality', 
the {\em periodic TLA},
on the other hand, is defined by abstract generators and relations. 
See 
\cite{PasquierSaleur90,Levy91,Martin91,MartinSaleur93,MartinSaleur94a,%
FanGreen99,DeGier02} 
and references therein for details of both. 
The relationship between the two versions is not so
straightforward as in the ordinary TLA case. 
See also \cite{Jones94b,GrahamLehrer98}.

To summarize the
naming conventions: \TL\ algebras are defined by
  generators and relations. Blob, contour, partition algebras (and
  others with diagram suggestive names) are defined via bases of
  diagrams and diagrammatic composition rules. 

\prl(hecke)
The map $u_i \mapsto U_i$ ($i>0$), $u_0 \mapsto e$,
extends to an  algebra homomorphism $\phi$ from  $T(B_n)$
to $TLb_n^K$ in the case 
$q=q_1$, 
$\delta_e= q_0+q_0^{-1}$
and $\gamma = \frac{q_0^2 + q^2}{q_0 q}$. 
\end{pr}
{\em Proof:} The relation checking is largely routine. 
Note from (\ref{uuuu}) that 
$\phi(u_1u_0u_1 -  \frac{q_0^2 + q_1^2}{q_0 q_1} u_1)$ vanishing is
sufficient to ensure 
$\phi(E_{B_2})=
\phi(u_0(u_1u_0u_1 -  \frac{q_0^2 + q_1^2}{q_0 q_1} u_1))=0$. 
\Qed

\medskip

Several authors have used this so-called `blob' homomorphism to investigate
Hecke algebra representation theory in the type-$B$ and type-$\hati{A}$
cases
\cite{GrahamLehrer03,CoxGrahamMartin03,MartinWoodcock03}. 
One final way to view the present paper is as
 a similar `blob' approach to type-$\hati{C}$. 
(The three parameter affine-$C$ Hecke algebra itself is of interest for a
 variety of reasons --- 
see for example \cite{Lusztig89C,Sahi99} and references therein.)

\subsection{Summary of notations}\label{SINot}

For the reader's reference we summarize here the 
notations for algebras used in  this paper 
(and indicate the  section in which each is defined):
\\
$b_n$ blob algebra, section~\ref{dir homs} \\
$b'_n$ \achiral\ algebra
  section~\ref{blob'} \\ 
$b_n^x$  \achiralb\ algebra,  section~\ref{ASTLA} \\
$b_n^{x'}$ big \achiralb\ algebra,  section~\ref{ASTLA} \\ 
$b_{2m}^{\phi}$ affine symmetric \TL\ diagram algebra,  section~\ref{ASTLA'} \\ 
$\Beta_n$ Brauer algebra, section~\ref{s diag} \\
$C_{n,m}(l)$ contour algebra, section~\ref{s contour} \\ 
$C_{n}^{\sim}(l)$ generalised contour algebra, section~\ref{s contour} \\ 
$H(G)$ generic Hecke algebra of graph $G$, section~\ref{INot} \\
$TL_n$  \TL\ algebra of type-A, section~\ref{INot} \\  
$TLb_n$ \TL\ algebra of type-B, section~\ref{INot} \\
$T(G)$ a quotient of $H(G)$, section~\ref{INot} \\

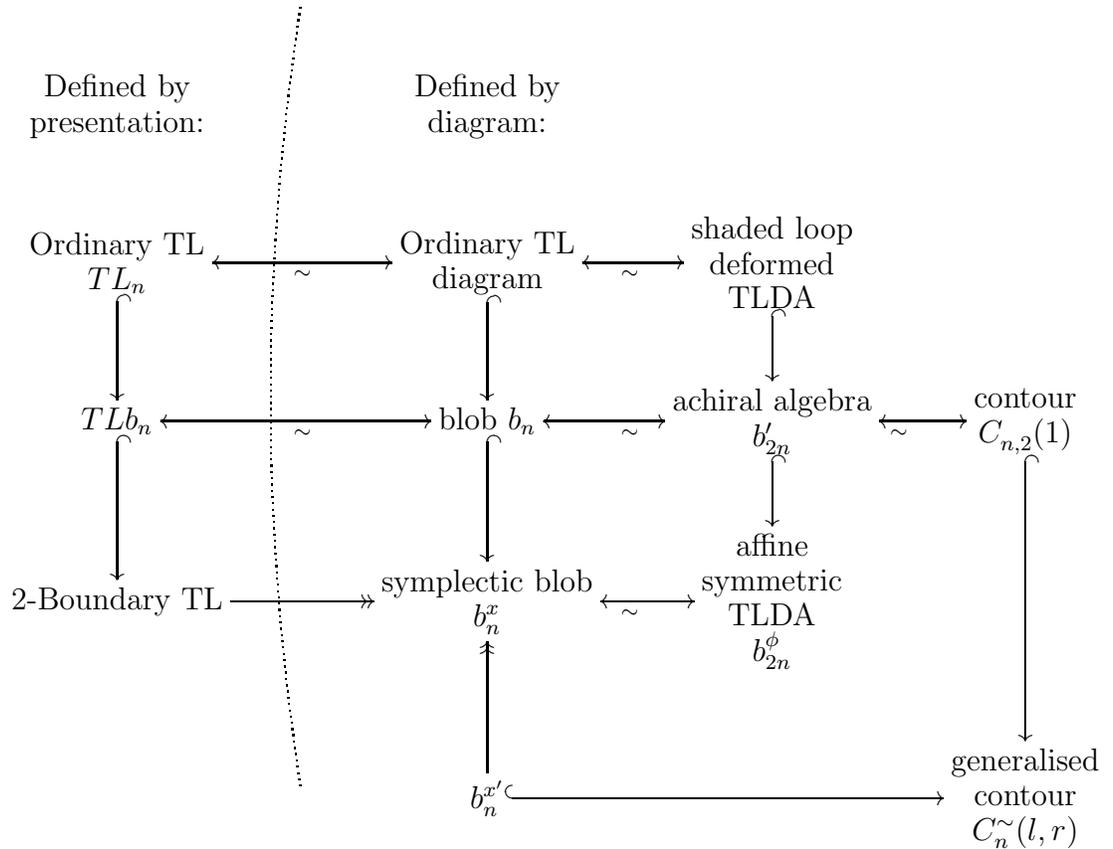
\begin{figure}
\[
\xymatrix{
 && \ar@/_1pc/@{.}[ddddd] 
\\
 & \txt{Defined by \\presentation:} && \txt{Defined by \\diagram:}
\\
   & \txt{Ordinary TL \\ $TL_n$} \ar@{<->}_{\sim}[rr]  \ar@{^{(}->}[d] 
               && \txt{Ordinary TL \\ diagram}  \ar@{<->}_{\sim}[r]  \ar@{^{(}->}[d] 
                          & \txt{shaded loop \\ deformed \\ TLDA}  \ar@{^{(}->}[d] 
\\
   & \txt{$TLb_n$}  \ar@{<->}_{\sim}[rr]  \ar@{^{(}->}[d] 
           && \txt{blob $b_n$}   \ar@{<->}_{\sim}[r]  \ar@{^{(}->}[d] 
                & \txt{\achiral\ algebra \\ $b'_{2n}$} 
                     \ar@{^{(}->}[d]  \ar@{<->}_{\sim}[r] 
		        & \txt{contour \\ $\RX^{}_{n,2}(1)$} \ar@{^{(}->}[dd] 
\\
   & \txt{2-Boundary TL} \ar@{->>}_{}[rr]
             && \txt{\achiralb\ \\ $b^x_n$}  \ar@{<->}_{\sim}[r]
                          & \txt{affine \\ symmetric \\ TLDA \\ $b^{\phi}_{2n}$}
\\
 &&& \txt{$b^{x'}_n$}  \ar@{->>}_{}[u] \ar@{^{(}->}[rr]
         && \txt{generalised \\ contour \\ $\RX^{\sim}_{n}(l,r)$}  }
\]
\caption{\label{schemat} Algebra relationships}
\end{figure}

The relationship between these algebras 
is indicated by the
schematic in figure~\ref{schemat}.

\newcommand{\gloss}[3]{$#1 \;$ #2, section~\ref{#3}\\ } 

\noindent
A glossary of other notations used for sets may also be useful:
\\
\gloss{B_n}{set of blob diagrams}{dir homs}
\gloss{B_n^x}{set of left-right blob diagrams}{pdcat}
\gloss{B_n^{\phi}}{set of left-right symmetric reduced periodic pseudodiagrams}{reduc}
\gloss{\DV}{set of beaded diagrams (given vertex set $V=V_n$ and bead set $S$)} {brauer}
\gloss{\DoV}{beaded pseudodiagrams}{brauer} 
\gloss{\Dnz}{planar diagrams}{s contour}
\gloss{{\Dnzl{n}{l}}}{planar $l$-exposed diagrams}{s contour}
\gloss{D_n=D_{n,m}}{planar diagrams with $\Z_m$-beads}{s contour}
\gloss{\Dnp}{cylinder embeddable diagrams}{deform}
\gloss{\Dnpc}{isotopy classes of concrete cylinder diagrams}{deform}
\gloss{\Dppc}{
   classes of concrete cylinder diagrams
   including non-contratible loops}{deform}
\gloss{\Dph_n}{classes of left-right symmetric concrete cylinder diagrams}{Dph}
\gloss{{\Jn{n}}}{unbeaded diagrams}{brauer}
\gloss{{\Dzo{n}}}{planar unbeaded diagrams}{dir homs}
\gloss{\Jp{n}}{periodic unbeaded diagrams}{deform0}
\gloss{\DB{n}}{symmetric planar unbeaded diagrams}{deform2}
\gloss{\DBe{n}}{symmetric planar unbeaded $ede$ diagrams}{recol b}

\section{Category theory preliminaries} \label{cat}

The starting point at the category theory level is as follows.  
Given an algebra
$A$ and an idempotent $e \in A$ then $eAe$ is also an algebra, and 
$Ae$ is a left $A$ module and a right $eAe$ module. 
Thus we may define functors between the category $A$-mod, 
of left $A$-modules, and $eAe$-mod:
\newcommand{\glob}{G}
\newcommand{\loc}{F}
\begin{eqnarray}
\label{glob}
\glob: \mbox{$eAe$-mod} & \longrightarrow & \mbox{$A$-mod}
\\ \non
M & \mapsto & Ae \otimes_{eAe} M
\end{eqnarray}
\begin{eqnarray}
\label{loc}
\loc:  \mbox{$A$-mod} & \longrightarrow & \mbox{$eAe$-mod}
\\  \non 
N & \mapsto & eN
\end{eqnarray}
with various powerful properties 
(summarized in \cite{Green80,MartinRyom02}). 
In particular if $A$ is an algebra over a field then 
$\loc$ is exact and $\glob$ is right exact. 
Further, the image of a simple
module under $\glob$ has simple head. 
\prl(grebo)\mbox{\rm \cite{Green80}}
Let $\{ S_{\lambda} \; | \; \lambda \in \Lambda \}$ be a complete set
of inequivalent simple (left) modules of $eAe$ over some field. 
Then $\{ \head(G(S_{\lambda})) \; | \; \lambda \in \Lambda \}$ 
is a set of inequivalent simples of $A$, and any simple $A$-module $S$
in an equivalence class not represented in this set obeys $eS=0$.  
\end{pr}

\newcommand{\prestandard}{prestandard} 
\newcommand{\Prestandard}{Prestandard} 
\subsection{\Prestandard\ modules}

The functor $\loc$ is called localisation, and $\glob$ is
globalisation, with respect to $e$.
(We may write $\loc_{e},\glob_{e}$ where convenient.)
Suppose that we are given an idempotent $e$ in an algebra $A$, and
that $S$ is a simple module of $eAe$. 
Then $\glob(S)$ is called the {\em \prestandard} module of $A$ associated
to $S$ by $e$.

Indeed, suppose that $e=e_1e_2=e_2e_1$ ($e_1,e_2$ also
  idempotent). Then the sequence of idempotents $1,e_1,e$ defines a
  sequence of algebras $A,e_1Ae_1,eAe$. 
A \prestandard\ module $M = \glob_{e_2}(S)$ of $e_1Ae_1$
  will globalise to a \prestandard\ module of $A$.
On the other hand if $M \neq \head\!({M}) $ then these two 
will not necessarily globalise to the same module
in $A$ 
(although both are of  \prestandard\ type by construction).
This makes the \prestandard\ notion less canonical, 
although more general, than the {\em standard} modules of
quasihereditary algebras, for example.  

\prl(orestes)
If $M$ is a \prestandard\ $A$-module then 
\\
(i) it has simple head
$L_M$ (say), and if $M_0$ is the maximal proper submodule then $M_0$
does not contain $L_M$ as a simple composition factor;
\\
(ii) if $A$ has an involutive antiautomorphism defined on it 
fixing $e$
then $M$ has at most one  
contravariant form defined on it,
up to scalars,
and the rank of any such form is the dimension of $L_M$.
\end{pr}
{\em Proof:} (i): Only the last claim remains to be proven. Suppose that
$M=G(S)$ for simple $eAe$-module $S$, then $F(M)=eM \cong S$.
In particular $e$ acts as zero on all but one simple factor in $M$.
Now suppose there exists a proper submodule $M'$ of $M$. 
If $eM' \neq 0$ then $eM'=S=eG(S)$, so $AeM'=AeG(S)=G(S)$,
which would imply $M' \supseteq G(S) \supset M'$ --- a contradiction.
Thus $F$ kills every proper submodule of $M$, so $e L_M \neq 0$.
\\
(ii): There is a one-to-one correspondence between such forms and 
homomorphisms from $M$ to its contravariant dual $M^{\times}$, 
but by (i) 
there is at most one such homomorphism (up to scalars),
whose image is $L_M$. 
To see this note that by the proof of (i)
$eL_N=0$ for every simple factor $L_N$ of $M$ not in the head.
Let $L_M^{\times}$ denote the contravariant dual of $L_M$, a simple module.
Since $eL_M \neq 0$ we have $eL_M^{\times}\neq 0$ so neither $L_M$ nor
$L_M^{\times}$ appears below the head in $M$. 
Thus $L_M$ does not appear above the socle in $M^{\times}$ and a homomorphism 
$M \rightarrow M^{\times}$ is only possible if it maps the head $L_M$ of $M$
to the socle $L_M^{\times}$ of $M^{\times}$, with $L_M \cong L_M^{\times}$.
\Qed

Propositions~\ref{grebo} and~\ref{orestes}
and the exactness properties make \prestandard\ 
modules potentially useful modules to study,
in representation theory.
In this paper we will encounter various 
modules (for algebras) with useful natural bases. 
It will be convenient, where possible, 
 to be able to identify these as \prestandard. 


\subsection{Globalisation and balanced maps}

The exactness properties of $\glob$ and $\loc$ are standard
results (see \cite{MartinRyom02,Green80}). 
However it is worth unpacking a little before we go on, 
since some of the mechanics will be used later. 
The first thing to recall is the notion of balanced map
\cite{Chevalley56,CurtisReiner62}. For $M$ a right module and $N$  a
left module of a ring $R$ with 1, a balanced map $f$ of $M \times N$ into
an additive abelian group $P$  is a map such that
$f (m+m',n)=f(m,n)+f(m',n)$,
$f(m,n+n')=f(m,n)+f(m,n')$,
$f(m,rn)=f(mr,n)$. 

The map that takes $(m,n)$ to $m \otimes n  \in M \otimes_R N$ is a
balanced map. 
If $f:M \times N \rightarrow P$ is balanced then there is a
homomorphism $f^* : M \otimes_R N \rightarrow P$ such that 
$f^*(m \otimes n) = f(m,n)$ 
(in fact $f^*$ is uniquely determined by $f$) \cite{CurtisReiner62}. 

Now consider 
\[
\loc( \glob (N)) = eAe \otimes_{eAe} N \stackrel{\mu}{\rightarrow} N
\]
where $\mu$ is derived from the (NB, balanced) map
\[
(a,n) \mapsto an . 
\]
We may define a homomorphism $\nu: N \rightarrow eAe \otimes_{eAe} N$
by 
\eql(inve)
\nu (n) = e \otimes n
\eq
Obviously $\mu \nu$ is the identity map on $N$; and
\[
\nu ( \mu ( a \otimes n) = \nu( an) = e \otimes an = a \otimes n
\]
since $a \in eAe$, so $\mu$ is an isomorphism. 


Similarly we have that 
\eql(iso x)
G(eAe) = Ae \otimes_{eAe} eAe \cong Ae . 
\eq
More generally, suppose that $S$ is a left sub-$eAe$-module of $eAe$
(i.e. a left ideal),
then there is a multiplication map
\begin{eqnarray*}
\mu :  Ae \otimes_{eAe} S  & \rightarrow &  AeS   \\
ae \otimes s  & \mapsto  &  aes 
\end{eqnarray*}
(in the rest of this section, $\mu$ applied to a tensor product 
of this form will
always be the appropriate multiplication map)
however, this surjection need not be an injection in general. 
(There is a grotesque example in \cite{CurtisReiner62}.)
The issue is the construction of the `inverse' as in (\ref{inve}). 
The `identity element' $e \in eAe$ will not generally lie in $S$. 
On the other hand, suppose that there are $f,g \in eAe$ such that
$S=eAef$ and $fgf=f$. 
(Such an $f$ is said to satisfy the {\em return condition}.
We call a left $eAe$-ideal $S$ of form $eAef$ with $f$ satisfying
the return condition a {\em return ideal}.)
Then there is a map 
$\nu: AeS \rightarrow Ae \otimes_{eAe} S$ given by 
\[
\nu(x) = x \otimes gf
\]
so that $\mu(\nu(x)) = xgf=x$ and 
$\nu(\mu(a \otimes s))= \nu(as) = as \otimes gf= a \otimes s$. 
Therefore
\prl(lemmin)
If $S=eAef$ is a left ideal of $eAe$ generated by $f \in eAe$ such
that $fgf=f$ for some $g \in eAe$, then the multiplication map $\mu$
is an isomorphism
\[
G(S) \cong AeS = AS = Af
\]
(NB, $\mu$ and its inverse are given explicitly). 
In particular the set inclusion of $S$ in $AS$ passes to an injection
$\nu$ of $S$ into $G(S)$. This is not an algebra-module map,
but if $D$ is a linearly independent set in $S$ then it is linearly
independent in $AS$ and $\nu(D)$ is in $G(S)$. 
\end{pr}
Note that $fg$ is idempotent, so 
\[
S = eAef \doublerightarrow eAefg
\]
is a surjective map to a projective $eAe$-module. 


\subsection{Module bases under globalisation}

The functors $\loc,\glob$ are tools for analysing categories of
modules,
rather than specific bases or representations. 
However, 
following the discussion above, there are realistic cases in which 
one can use a basis for $S$ to construct a basis for $\glob(S)$. 
(This is particularly so for diagram algebras, which come with a
diagram basis for elements $f$ of which the return condition $fgf=f$ is
always true for some algebra element $g$. Indeed $g$ can usually be
chosen a basis element, or else a scalar multiple thereof.)


\prl(basis thang)
Let $S$ be a submodule of the (left) regular module of $eAe$, 
and suppose that this submodule has basis $D$. 
Then the concrete set of elements $e \otimes D$ 
is independent in and 
will also generate $\glob(S)$. 
Further, 
$$
\glob(S) \stackrel{\mu}{\doublerightarrow} AS \hookrightarrow \glob(eAe) 
\stackrel{(\ref{iso x})}{\cong} Ae \hookrightarrow A
$$ 
as left $A$-modules. 
\end{pr}
{\em Proof:} 
For $d \in D \subset S$ then $e \otimes d \in  Ae \otimes_{eAe} S = \glob(S) $. 
The image  $\mu(e \otimes d) = d$, so $\mu(e \otimes D) = D$. 
On the other hand, the direct set map $D \hookrightarrow AD$ is an inclusion, 
so $D$ is linearly independent in $AS$.
The multiplication map 
$\mu : \glob(S) \rightarrow AS$
is surjective, not necessarily bijective, 
but it is still a module homomorphism.
Thus if  $e \otimes D$ {\em were} to be linearly {\em dependent} in $\glob(S)$,
the image $D$ would be linearly dependent in $AS$ --- a contradiction. 
So the set $e \otimes D$ is linearly independent. 

On the other hand $D$ spans $S$, so $Ae\otimes D$ spans $\glob(S)$. 
Thus $e\otimes D$ extends to a basis of $\glob(S)$ by the exchange
theorem. 
If the multiplication map $\mu$ is bijective then $\glob(S)$ has a
natural isomorphic image in $A$, with a basis which contains $D$ as a
subset.

Finally
$eD=D \subset eAe$, so $AeD \subset Ae$. 
\Qed


\prl(flat)
Suppose that $ S_1 \stackrel{\psi}{\hookrightarrow} S_2 $ 
is an inclusion of return ideals. 
Then $G(S_1) \stackrel{G(\psi)}{\hookrightarrow} G(S_2)$,
i.e. $G$ behaves as if left exact.
\end{pr}
{\em Proof:} Unpacking the assumptions
then $eAef_1 \hookrightarrow eAe f_2$, so $f_1 \in eAe f_2$.
Let us say (WLOG) $f_1 =  f f_2$. 
We have 
\[
\xymatrix{
eAe f f_2 & \stackrel{\psi}{\hookrightarrow} & eAe f_2  \\
\downarrow G &             & \downarrow G  \\
Ae \otimes_{eAe} f f_2  & 
            \stackrel{G(\psi)}{\rightarrow} & Ae \otimes_{eAe} f_2  \\
\downarrow \mu_1 &             & \downarrow \mu_2  \\
A f f_2  & \hookrightarrow &  A f_2
}
\]
Since both $\mu$-maps are isomorphisms, 
there are two ways of constructing a homomorphism in the middle:
via the functor $\glob$; or via the bottom row inclusion
$\nu_2(\mu_1(a\otimes ff_2)) = \nu_2(aff_2) 
= aff_2 \otimes g_2f_2 = a \otimes ff_2 = a f \otimes f_2$.
Again since the multiplication maps are isomorphisms, the latter
construction 
is an injection, i.e. $\glob(S_1) \hookrightarrow \glob(S_2)$. 
On the other hand 
$G(\psi)(a \otimes ff_2) = a \otimes \psi(ff_2)
  =a\otimes ff_2 \in G(S_2)$,
so $G(\psi)$ and $\nu_2 \circ \mu_1$ are the same map.
(NB, for $M \hookrightarrow N$ the element
$a \otimes m \in G(M)$ is not the same thing as 
$a \otimes m \in G(N)$ in general --- a set of such objects can
be independent in $G(M)$ and not in $G(N)$, but here
we can also build such objects on the $G(S_2)$ side 
by the kernel-free $\nu_2 \circ \mu_1$ route.)
\Qed

Indeed, suppose that 
\eql(inc seq1)
S_1 \hookrightarrow S_2 \hookrightarrow \ldots \hookrightarrow eAe
\eq
is a nested sequence of return ideals,  
and $D_i$ a basis for each such
that $D_i \subset D_{i+1}$. Consider the sections defined by this
sequence: 
\[
0 \rightarrow S_{i} \stackrel{\psi}{\rightarrow} S_{i+1} 
  \rightarrow S_{i+1}/S_{i} \rightarrow 0
\] 
then $\glob(S_{i}) \hookrightarrow \glob(S_{i+1})$ 
and $\glob(S_{i+1}/S_{i}) = \glob(S_{i+1}) / \glob(S_{i})$.
Thus  in particular if $S_{i+1}/S_{i}$ is simple then
$\glob(S_{i+1})/\glob(S_{i})$ is \prestandard. 

\medskip

It may be that (\ref{inc seq1}) is valid over a ground ring 
that 
specialises to a field in a number of different ways. 
Then the sections of the sequence make sense (both before and after
globalisation) over the ring, and \prestandard\ modules have the flavour of
Specht modules \cite{JamesKerber81}. 

Diagram algebras come with bases with special properties. 
We will use this concrete construction to decompose the regular module
explicitly throughout entire towers of algebras. 

\section{Diagram algebras: initial examples} \label{s diag}


In order to define diagram algebra quotients of Hecke algebras later,
we start by defining related algebras which have both diagram and
`linear' realisations. 

\subsection{Brauer algebra wreaths} \label{brauer}

Fix $n,m \in \N$ with $n+m$ even, and let 
$\V^n_m=\{1,2,\ldots,n,1',2',\ldots,m'\}$, 
called the set of vertices.
Write $\V=\V_n$ for $\V^n_n$. 
Write $J^n_m$ for the set of pair partitions of $\V^n_m$
(so $J_n := J^n_n$ is the usual basis of the Brauer algebra $\Beta_n$ 
\cite{Brauer37,Weyl46}). 
For $S$ a set, 
an $S$-decorated pair partition is 
an element of $J^n_m$ 
together with a map from the set of pairs to the set of words in $S$. 
Fixing $S$, 
write $\DVnm = \DSVnm$ for the set of such objects. 
Thus $\DD_{\emptyset}(\V^n_n) \cong J_n$. 

Let
\[
p = \{\{ i, j \},\{ k, l \}, \ldots \}
\]
be a pair partition of $\V^n_m$. Then we may write $d \in \DVnm$ as 
\[
d = \{\{ i, j \}_{w_1}, \{ k, l \}_{w_2}, \ldots \}
\]
where $w_l$ is the word in $S$ associated to the $l$-th pair. 
(We adopt the convention of omitting the subscript when the word is
the empty word.)
We call this the {\em serial} realisation of  $d \in \DVnm$.
We next describe two further useful realisations. 


It will be helpful to think of the following mild extension to Weyl's
\cite{Weyl46} diagram realisation of the pair-partition basis in the 
Brauer algebra. Consider: 
\newline (i)
the vertices as arranged on a rectangular frame, 
$1$ to $n$ across the top edge, $1'$ to $m'$ across the bottom; 
\newline (ii)
the pairings as pieces of string (called {\em lines}) 
appropriately connecting vertices;
\newline (iii)
the accompanying elements of $S$ as threaded beads, 
threaded in the order indicated by the word 
(reading from vertex $i$ to $j'$, or from $i$ to $j$ if $i<j$, 
or from $i'$ to $j'$ if $i<j$). 
\newline
(Note that any tangling of strings, 
perhaps arising from some perceived embedding in an underlying space, 
is irrelevant here --- it is only
the pairings they define that matter.)


A third realization is achieved by arbitrarily embedding 
\cite{Armstrong79}
each string $\{i,j \}$ as a line from $i$ to $j$ 
in the plane region bounded by the frame rectangle. For example
\[
\{\{1,3'\}_{ab}, \{2,3 \}_{c}, \{1',2' \}_{} \}
=
\raisebox{-.21in}{
\includegraphics{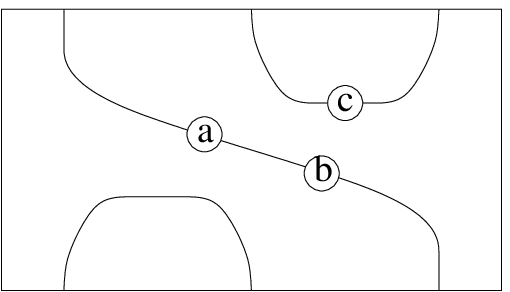}
}
\]
NB, it is not possible in general to do this without distinct lines
crossing (see later, and cf. \cite{JonesPlanar}).

In summary: we will call objects in the first realisation 
decorated pair partitions; 
objects in the string/bead realisation diagrams;
and objects in the third realisation concrete diagrams.

\medskip

\noindent {\bf Pseudodiagrams}
\\
Consider the idea of closed loops of string in the string/bead picture
--- that is, strings that do not end on any of the vertices. 
Corresponding to this, 
it will be convenient to extend the notion of 
decorated pair partitions to include 
(possibly multiple) copies of the empty set in the partition.
(Hence we make mild abuse of this terminology.)
Following \cite{MartinWoodcock2000} we call such diagrams augmented by
zero or more loops {\em (Brauer) pseudodiagrams},
and the extended partitions {\em pseudopartitions}. 
For example (with $\{\}^l$ denoting $l$ copies)
\[
d = \{\{1,3'\}_{ab}, \{2,3 \}_{c}, \{1',2' \}_{}, 
\{\}_{def}, \{\}^2    
\}
= \{\{1,3'\}_{ab}, \{2,3 \}_{c}, \{1',2' \}_{}, 
\{\}_{edf}, \{\}^2  
\}
\]
As before, any perceived embedding of these loops in an underlying
space is irrelevant --- thus in particular a loop does not have an
orientation or starting point for the reading off of bead sequences.
Thus $\{ \}_{def} = \{ \}_{edf}$.

\del(c_d)
We write $\DoVe{n}{m}$ for the set of $(n,m)$-pseudodiagrams.
If $d \in \DoVe{n}{m}$ is a pseudodiagram 
let $c_d \in \DVe{n}{m}$ denote the underlying diagram,
that is, the diagram obtained by omitting any loops;
and $ c^o_d \in \DoVe{0}{0}$ the complement, obtained by keeping only loops.
\end{de}


For $w$ a word in $S$ let $w^o$ denote the opposite word 
(the 
word with the same letters, but written in the reverse order).
In the serial (pair-partition) realisation, 
if we write a vertex pair in a definite
(not necessarily canonical) order: $(i_1, i_2)$, then the string/bead
datum for this string obeys
\eql(string op)
(i_1, i_2)_{w} = (i_2, i_1)_{w^o}
\eq
On the other hand if a pseudodiagram has a closed loop 
then this may again have beads on it, but reading the
sequence of beads depends on an arbitrary choice of starting point and
direction round the loop. 
We say two words are {\em loop equivalent} 
if one can be changed to the other or its opposite 
by any cyclic permutation. (Note that loop equivalence is an
equivalence relation on the set of words.)
Thus the bead sequence on a loop is only defined up to loop equivalence.  
 
\medskip

We {\em concatenate} pseudodiagrams 
$d_1 \in \DoVe{n}{m}$,  
$d_2 \in \DoVe{m}{l}$ 
to form a pseudodiagram 
$d_1d_2 \in \DoVe{n}{l}$ 
as follows.  
Pass to the string/bead realization and there juxtapose the vertices
$1',2',\ldots,m'$ in $d_1$ with the corresponding unprimed vertices in
$d_2$. 
Some of 
the chains of string resulting from this concatentation
will connect pairs among the
unprimed vertices in $d_1$ and the primed vertices in $d_2$
(defining a new diagram $c_{d_1d_2}$ on these vertices), 
and some will form closed loops. 
Because of the string/bead realisation we call this the 
{\em abacus product} on pseudodiagrams. 

For example
\eql(frubeddd)
 d_1 . d_2 . d_3 = \; 
\raisebox{-.521in}{
\includegraphics{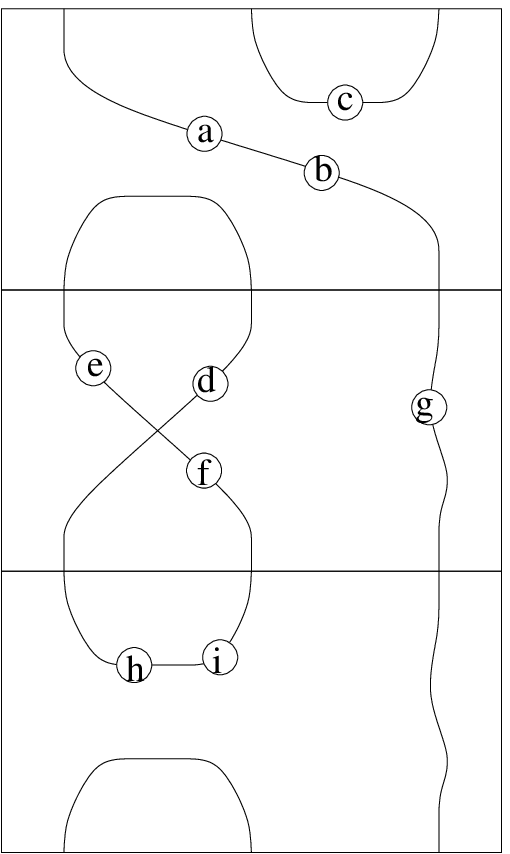}
}
\eq
is 
$$
\{\{1,3'\}_{ab}, \{2,3 \}_{c}, \{1',2' \}_{} \}
. \{\{1,2'\}_{ef}, \{3,3' \}_{g}, \{2,1' \}_{d} \}
. \{\{1,2\}_{hi}, \{3,3' \}_{}, \{1',2' \}_{} \}
$$ $$
= \{\{1,3'\}_{abg}, \{2,3 \}_{c}, \{1',2' \}_{}, \{ \}_{efihd} \}
 .$$ 

To confirm that this product is well defined note that 
in the ordered pair form of the serial realisation (\ref{string op})
composition is given by a sequence of $m$ replacements of form: 
\[
\underbrace{\{ \ldots (i_1, i_2')_{w} \ldots \}}_{d_1}
\underbrace{\{ \ldots (i_2, i_3)_{w'} \ldots \}}_{d_2} \;  \leadsto 
\{ \ldots (i_1, i_3)_{ww'} \ldots \}  
\]
and 
\[
(i_1, i_1)_{w} \leadsto ()_{w} = \{ \}_{w}
\]
the order of application of which, when not forced, produces no ambiguity.


\prl(abacus) The abacus product is associative. 
The closed case with $n=m=l$ is unital. 
Hence $\Do$ is a category with object set $\N$ and $(n,m)$ morphism set
$\DoVe{n}{m}$.  
\end{pr}
{\em Proof:} Associativity follows from the construction. 
For an example, consider $(d_1 d_2) d_3$ and 
$d_1 (d_2 d_3)$ in equation~(\ref{frubeddd})
--- one draws the same picture in each case. 

The identity element is the pair partition $\Id$: 
\eql(Id)
\Id = \{\{1,1'\},\{2,2'\},\ldots,\{n,n'\}\}
\eq
(all words empty). 
\Qed


Fix $\Ring$ a ring and $\delta_w \in \Ring$ for each $w$ a 
loop class representative word in $S$. 
\del(k_d)
For each $d \in \DoVe{n}{m}$ define a scalar $k_d \in \Ring$ by 
$$
k_{d} = \prod_{l} \delta_{w_l}
$$ 
with a factor
$\delta_{w_l} $ for each closed loop $l$ in $d$ with bead sequence $w_{l}$. 
\end{de}
An example follows shortly. 


Define a map
\begin{eqnarray} \label{compose}
\DVe{n}{m} \times \DVe{m}{l} & 
    \rightarrow & \Ring \times \DVe{n}{l} \rightarrow \Ring\DVe{n}{l} 
\\ \non
(d_1,d_2) &  \mapsto  &   (k_{d_1 d_2}, c_{d_1 d_2}) 
                \mapsto  k_{d_1 d_2} c_{d_1 d_2}
\end{eqnarray}
(recall that the strings of diagram $c_{d_1 d_2}$ are the open chains 
from the abacus product, 
each carrying the accumulated beads of this chain in the natural order). 

Since the underlying abacus product is associative, 
this (\ref{compose}) extends to  
an associative unital product on $\Ring\DVe{n}{n}$. 
For  example, 
in equation~(\ref{frubeddd}) $(d_1 d_2) d_3 = \delta_{w}   
\{\{1,3'\}_{abg}, \{2,3 \}_{c}, \{1',2' \}_{} \} $ 
where $w$ represents the class containing $efihd$.

Remark:
This product is amenable to massive generalisation, 
which we will largely ignore,
but see for example \cite{Martin94,JonesPlanar,AlvarezMartin05}. 
A milder generalisation is the {\em cyclotomic} variant, 
in which 
\eql(t def)
(i_1,i_2)_{s_1 s_2 \ldots} =  (i_2,i_1)_{\ldots s_2^t s_1^t} ,
\eq
where $t$ is an involutive map on $\langle S \rangle$ that does not
necessarily fix elements of $S$. 


\subsection{Periodic pair-partitions} \label{deform0}

Consider the `infinite' rectangular frame in which the primed and
unprimed vertices are labelled by corresponding 
copies of the set of integers. 
The set of arbitrary pair-partitions of this vertex set 
(call it $V_{\infty}$) is rather
unmanageable, but there are a number of more manageable subsets which
are closed under composition (ignoring loops). 
A pair-partition is said to be $n$-periodic if for every pair
$\{i,j\}$ there are pairs $\{ i \pm n , j \pm n \}$
(with $m' \pm n \; := \; (m\pm n)'$). 
It follows that there are only $n$ distinct orbits of pairs.
An $n$-periodic pair-partition can be specified by listing a
fundamental subset of $n$ pairs.

The first element in each such pair (at least) can be chosen to lie in the
fundamental set $\V^n_n$. 
Then if $(i_1,i_2)=(i_1, \bar{i_2}+mn)$ 
where $\bar{i_2}$ also lies in  $\V^n_n$ 
we might write  $(i_1,\bar{i_2})_m$ for  $(i_1,i_2)$. 
Using this notation we can demonstrate
a map from decorated pair partitions on $\V^n_n$ 
with bead set $S=\{ L_+, L_- \}$, 
and $L_+^t = L_-$ the involution as in equation~(\ref{t def}), 
to $n$-periodic pair partitions. We take a string with $m$ beads 
$(i_1,i_2)_{L_{\pm}^m}$ to $(i_1,i_2)_{\pm m}$ ($m \geq 0$).
In other words each bead $L_+$ corresponds to winding once clockwise
round the period, and $L_-$ is anticlockwise (thus we 
 take the quotient with  $L_+ L_- =L_-L_+=1$).

We write $\Jp{n}$ for the set of  $n$-periodic pair-partitions. 
There are infinitely many of these. 
For example with $n=1$ 
we have   
$\{\{\{1,m'\}\} \; | \; m \in \Z \}$
(writing only a fundamental subset for each partition).

For later convenience it will be useful sometimes to index vertices by odd
integers rather than all integers.
When we do this we will write the pair partition $\{\ldots \}_{o}$. 
For $n=2$ examples include
\[ \{\{ 1,7 \},\{ 1',7' \}\}_{o} . 
\]
Note that this pair-partition can be realised by vertex-connecting lines
embedded in the infinite rectangular interval, as in the finite case.
But note that in this example these lines necessarily cross
(the orbit of $\{ 1,7 \}$ includes $\{-3,3\}$ and $\{5,11\}$ for example).
For $n=4$ examples include
\[ \{\{ 1,-1 \},\{ 3,5 \},\{1',-1' \},\{3',5' \}\}_{o}
\qquad  \{\{ 1,3 \},\{ -1,-3 \},\{1',3' \},\{-1',-3' \}\}_{o} .
\]
Note that these particular examples can be realised by non-crossing lines.


Under composition, ignoring loops for a moment,
it is a simple exercise to show that $n$-periodicity is preserved. 
Two types of loops can appear: an `orbit' of loops individual members
of which are periodic images of one another; 
and an individual `non-contractible' loop which is 
mapped into {\em itself} by periodicity.
Finitely many instances of each type may be created in composition. 
We will see 
in section~\ref{deform}  
a natural way to keep track of these 
(and any possible decorations). 
Thus, given two different types of loop, we may introduce a set of
periodic pseudodiagrams, which is then closed under composition.
It will also be convenient to introduce the set 
$\Jpp{n}$ 
of periodic diagrams
augmented just by the non-contractible type of loops.
A suitable collection of relations removing the orbits of loops then 
makes $\Jpp{n}$ a basis for an algebra generalising the Brauer algebra
(cf. \cite{OrellanaRam01}). 

There is an injective homomorphism from $\Jn{n}$ into $\Jp{n}$ which
simply uses $p \in \Jn{n}$ as the fundamental subset. This can be
extended to an algebra map. 

We write $\Jp{n}^{\!\! S}$ for the beaded version of $\Jp{n}$. 


\subsection{Planar embeddings and contour algebras} \label{s contour}

\begin{de}\label{dplanar}
(1) A diagram in $\DVe{n}{m}$ is called 
{\em planar} if it is possible to embed the
strings in the plane interior to the rectangle 
(touching the boundary only at the vertices)
in such a way that they do
not (self-intersect or) touch one another. 
\\
(2) Any specific such embedding is called a {\em concrete} planar diagram. 
\\
(3) If one concrete planar diagram may be continously deformed into another,
with all the intermediate stages concrete planar diagrams, then the diagrams
are said to be {\em \aipic} \cite{Moise77}. 
\end{de}

\prl(post3)
Two concrete planar diagrams are isotopic if and only if they are
realisations of the same underlying diagram. 
\Qed
\end{pr}

\begin{rem} \label{cave} {\em
If we consider planar pseudodiagrams in the same way, then a given
Brauer pseudodiagram may have more than one isotopy class of planar
embeddings. Note however that both $c_d$ and $k_d$ can be considered
as applying to planar pseudodiagrams via their underlying Brauer
pseudodiagrams. 
}
\end{rem}

Write $\Dnz \subset \DVe{n}{n}$ 
for the set of planar diagrams. It will be evident that
the restriction of diagram composition (\ref{compose})
 to $\Dnz$ closes on $\Ring\Dnz$. 


\del(covered)
A string in a planar diagram is called {\em exposed} (or 0-covered) if it may
be deformed \aipic ally to touch the western frame edge,
and $l$-covered if it may be deformed to touch an $(l-1)$-covered line
and no lower. 
\end{de}

Let $\Dnzl{n}{l} \subset \Dnz$ 
denote the subset of planar diagrams in which only the
$l'$-covered lines with $l' \leq l$ may be decorated 
(i.e. beaded, i.e. map to other than the empty word). 
\begin{pr}
The restriction of diagram composition (\ref{compose})
 to $\Dnzl{n}{l}$ closes on
$\Ring\Dnzl{n}{l}$. 
\end{pr}
{\em Proof:} Composition may expose new line segments, but it cannot
cover any that were previously exposed (since the relevant part of the
western frame is still  in place). Thus any decorated (hence no more
than $l$-covered) line remains no more covered (hence decorable) in
composition. \Qed


\medskip
\noindent
We now note some specialisations with interesting finite dimensional quotients.

Fix  bead set $S$ of order one 
($S= \{ L \}$, say).
It follows that words in $S$ are all of form $L^i$, 
and that each is in a separate loop class.
Write $\delta_i \; := \; \delta_{L^i}$. 
Fix $m \in \N$ and
let $D_{n}=D_{n,m}$ denote the subset of $\Dnz$ in which no
string carries more than $m-1$ beads. 
Consider the $\Ring$-algebra with basis $\Dnz$ 
and $\delta_{i+m}=\delta_i$. 
With this specialisation of the parameters we may 
impose the quotient relation 
that $m$ beads together may be cancelled ($L^m =1$). 
This produces an algebra $\RX_{n,m}$ with basis $D_{n,m}$ 
($D_{n,m}$ is clearly spanning; to see that it is independent
note that the relation cannot be used to change the {\em shape} of a
diagram).  
\del(contour)
Define $\RX_{n,m}(l)$ as the subalgebra 
of  $\RX_{n,m}$ 
spanned by $\Dnzl{n}{l}$
(and hence with basis $\Dnzl{n,m}{l} \; := \; \Dnzl{n}{l} \cap D_{n,m}$). 
These are called contour algebras (see \cite{CoxMartinParkerXi03}). 
\end{de}


More generally, the `algebra' taking place on a single string is the
free monoid on the generators $S$. 
Every quotient by some set of relations $\sim$ to a finite monoid
(or even $\Ring$-algebra with basis a finite subset of the free
monoid),
together with a consistent specialisation of the parameters,
induces a finite generalised contour algebra $\RX^{\sim}_{n}(l)$.
(Here $\Ring$-linear combinations of words on a string pass linearly to
corresponding $\Ring$-linear combinations of diagrams.) 

As usual the identity in these contour algebras is the pair partition
$\Id$. 
Suppose that $L \in S$. 
We define $L_i \in \DVe{n}{n}$ as the diagram that is $\Id$ 
as a pair partition, but has the single letter word $L$ 
on the $i$-th string, with all other words empty.

Define 
$U_i$ as the diagram differing from $\Id$ in having the pairs 
$\{i,i+1 \}, \{ i', (i+1)' \}$. 

\prl(claim1)
The algebra $\RX_{n,m}(l)$ with $m>1$, $l<n$, $S=\{ L \}$, is generated by the set
$$
\{ \mathbb{I} \} \cup   \{ L_i \}_{i=1}^{l+1} 
 \cup   \{ U_i \}_{i=1}^{n-1} .
$$
\end{pr}
{\em Proof:} See Appendix. \Qed
\newcommand{\toappendix}[2]{\newcommand{#1}{{ #2 }}}
\toappendix{\stuffb}{{
\section{Proof of Proposition~\ref{claim1}}

To prove: that 
the set $D_{n,m}^{z,l}$, can be generated by the 
set 
$$B := \{ \mathbb{I} \} \cup   \{ L_i \}_{i=1}^{l+1} 
 \cup   \{ U_i \}_{i=1}^{n-1} .$$

Clearly the set of diagrams generated by $B$ is contained in
$D_{n,m}^{z,l}$, since concatenation can never increase the level of
coveredness of a particular string.
So we need only prove that $D_{n,m}^{z,l}$ is contained in the set of
diagrams generated by $B$.
We sketch a proof of this by induction on $l$.

If $l=-1$ and there are no decorated lines, then this is just the
result for the diagram version of the Temperley--Lieb algebra
\cite{Martin91}.

Now suppose $l=0$ and we have a diagram with a decorated $0$-covered
line with $j$ beads.
Since the decorated line is $0$-covered there are two
possibilities. Either the line is the string starting at the first
position or ending at the first position
- in which case we can decompose in the diagram into a
product of $j$ $L_0$'s together with a smaller diagram with one less
$0$-decorated line, 
or vice versa, 
or the line is a starting at a greater position
than
the first. In this case we get a diagram that looks like figure 
\ref{fig:gens0} where we have only decorated
the $0$-covered line with one bead rather than $j$ beads for simplicity.
We have drawn a propagating $0$-covered line; the dashed line
represents a non-propagating $0$-covered line. 
\begin{figure}[ht]
\begin{center}
\epsfbox{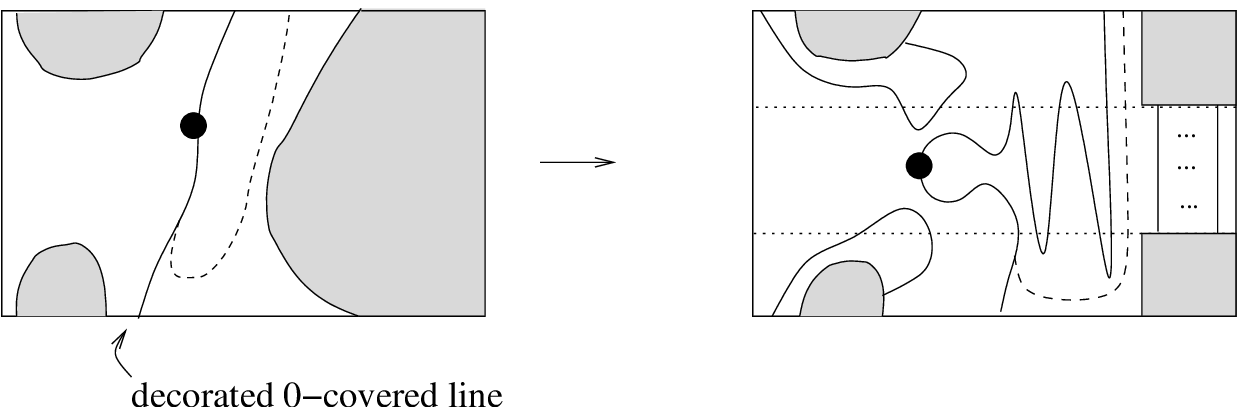}
\end{center}
\caption{\label{fig:gens0}
}
\end{figure}
Now note that the grey regions to the left of the decorated line
cannot contain any propagating lines - since the decorated line is
$0$-covered.
But both these regions must contain at least one string that is
$0$-covered and so we can deform the diagram to look like that on the 
right hand side of figure \ref{fig:gens0}.

We can arrange it so the two grey regions on the right hand 
side are only joined by
propagating lines and the number of these lines $x$, say will be the
same (\resp\ different) parity as $n$ if the decorated $0$-covered line is
non-propagating (\resp\ propagating). Thus the difference $n-x$
is even if the $0$-covered line is non-propagating and odd if the
$0$-covered line is propagating.
We now ``wiggle'' the $0$-covered line enough times 
so that we get the right number
of lines so that the
middle section of the diagram enclosed in dotted lines is now the
diagram product $U_1L_1U_2U_1$ (which has $n-3$ propagating lines). 

So we can decompose the diagram
into a product of three smaller diagrams, the outside diagrams having
a smaller number of $0$-covered lines.
 
The case with $l \ge 1$ is similar and is illustrated 
in figure~\ref{fig:gens1} below.
We have drawn a propagating $l$-covered line; the dashed line
represents a non-propagating $l$-covered line. 
We can assume that the $l-1$-covered line (which may be decorated or
not) is propagating, for otherwise the $l$-covered line would not be
$l$-covered.

\begin{figure}[ht]
\begin{center}
\epsfbox{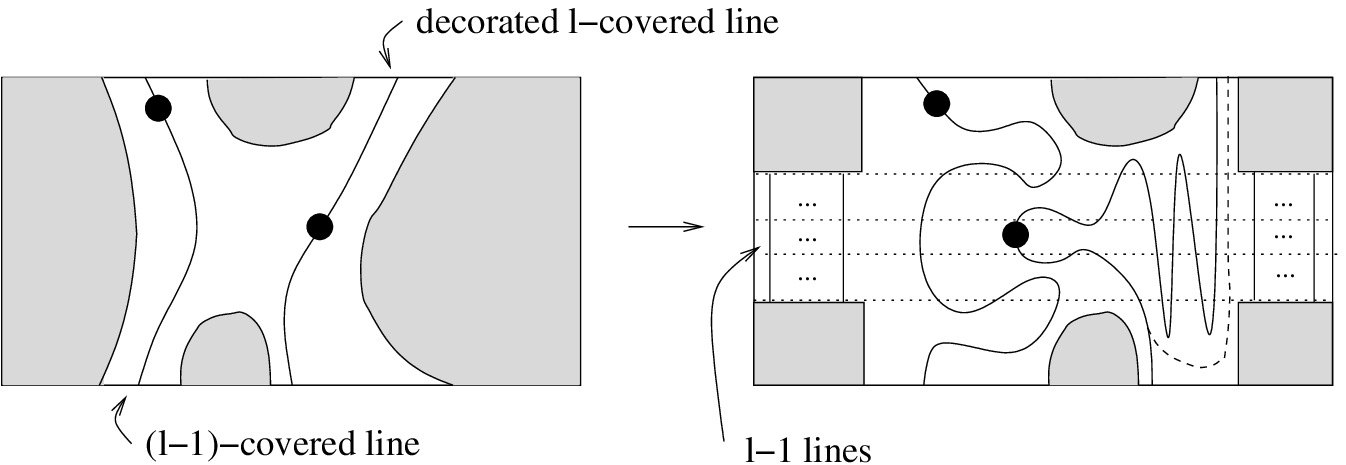}
\end{center}
\caption{\label{fig:gens1}
}
\end{figure}
Note that now when we ``pull apart'' the grey regions we can always
stretch then so that we can get only propagating lines joining the two
smaller
grey regions.
Also note that the number of propagating lines in the grey region on the
left hand side is at exactly $l-1$, for otherwise the $l-1$-covered
line would not be $l-1$-covered.
We again get the right parity, so that we can ``wiggle'' the
$l$-covered line so that 
the middle section of the diagram enclosed in dotted lines is now the
diagram product $U_{l+1}L_{l+1}U_{l+2}U_{l+1}$, 
and so we can decompose the diagram
into a product of three smaller diagrams, the outside diagrams having
a smaller number of $l$-covered lines.


}}


Note that in case $m=2$ we may replace 
the relation set $\sim = \{ LL = 1 \}$
(giving $L_i L_i = \Id$)
with $\sim = \{ LL=L \}$ 
(giving $L_i L_i = L_i$)
and obtain an isomorphic algebra. 
Indeed we may deform the monoid (algebra) via $LL=\kappa L$ similarly
(giving $L_i L_i = \kappa L_i$). 
Only the case $\kappa=0$ departs from the rest.


It will be evident that similar definitions 
to $\RX^{\sim}_{n}(l)$ may be contructed for
{\em eastern} exposure, and for composites
$\RX^{\sim}_{n}(l,r)$. 
Also of interest, as it turns out, are subalgebras of the case
$S=\{ L, R \}$ 
generated by 
$$
\{ \Id, \; L_1, \; R_n \} \cup \{ U_i \}_{ i=1}^{ n-1}  .
$$ 
where $\sim$ defines a certain noncommutative monoid. 
(Note that this choice of generators 
prescribes the way in which $L$ and $R$ can meet,
which leads to some interesting topological effects --- see later.)


\subsection{Direct homomorphisms with known algebras} \label{dir homs}

\del(TL)
A Temperley--Lieb (TL) diagram is an isotopy class of concrete planar
diagrams, or any representative thereof. Ordinary TL diagrams are beadless.  
\end{de}

By Proposition~\ref{post3}   
the subset $\Dzo{n} \hookrightarrow \Dnz \hookrightarrow \DVe{n}{n} $ 
with no beads 
is in bijection with the set of ordinary 
Temperley--Lieb
diagrams on two rows of $n$ vertices. 


\del(blob)
{\rm \cite{MartinSaleur94a}}
The set $B_n$ of {\em blob diagrams}
is the set of decorated  
TL diagrams on two rows of $n$ vertices
in which western exposed lines, only, may be
decorated, with at most a single bead (`blob') on each. 
\end{de}
For example, for any given $n$ the diagram $e \in B_n$ has the shape
of the identity diagram ($n$ vertical lines), but with the leftmost
line decorated with blob. 

The blob algebra
$b_n = b_n(\delta=q+q^{-1},\delta_e,\gamma)$
(as in \cite{MartinSaleur94a}, but as parameterised in 
\cite{MartinRyom02}), 
is generated by TL diagrams and $e$, 
with two blobs on a line appearing in composition replaced via: 
\eql(ee=e)
e e = \delta_e e
\eq
and a loop decorated by a blob replaced by a factor of $\gamma$.
Thus $b_n(\delta,\delta_e,\gamma)$ has basis $B_n$.


It is easy to show that 
\begin{pr} For all $n$:
\newline
(i) The subset $\Dzo{n}$ of $\Dnz$ 
generates a finite dimensional algebra 
isomorphic to
 the ordinary Temperley--Lieb algebra $TL_n(\delta_0=q+q^{-1})$ 
\cite{TemperleyLieb71}. 
\newline
(ii) The set $\Dnzl{n,2}{0}$ is 
essentially identical to the set of blob diagrams $B_n$,
with $L_1=e$. 
The algebra $\RX_{n,2}(0)$ is isomorphic to the blob
algebra $b_n$ \cite{MartinSaleur94a}
and to $TLb_n$ 
(more generally, $\RX_{n,m}(0)$ is the coloured blob algebra mentioned in
\cite{MartinWoodcockLevy00}). 
\newline
(iii) The case $\RX_{n,m}(n)=\RX_{n,m}(\infty)$ is isomorphic to the cyclotomic
Temperley--Lieb algebra \cite{CoxMartinParkerXi03}.
\Qed 
\end{pr}


Note (from section~2 of \cite{CoxMartinParkerXi03})
that the tower of recollement framework applies to all of the above algebras.

Following the largely physically motivated 
investigation of the ordinary Temperley--Lieb algebra in the 1980s,
the blob algebra was introduced in order to allow use of 
the western edge of the frame (actually either one) as a cohomology
seam, 
and hence to address the {\em periodic} Temperley--Lieb algebra
(again, originally, with physical motivation). 
The very first level of exposure is sufficient for this purpose (see
the literature, for example \cite{MartinSaleur94a,GrahamLehrer03}). 
Interest in the `homogeneous' case (i.e. not filtered by exposure) has been
slower to arise, but now see \cite{RuiXi04} and \cite{JonesPlanar}
(whose primary interest is in subfactors). 
The general 
intermediates have yet to find a physical application. 



\subsection{Diagram embeddings, subalgebras and deformations} \label{deform}
This section describes  `topological' realisations of 
certain subsets of diagrams, 
generalising the planar embedding of definition~\ref{dplanar}, 
and the deformations of the algebra product possible in these cases. 


\pdef{Isotopy:}
A brief remark is in order on the general notion of isotopy 
following from the benign paradigm in definition~\ref{dplanar}. 
A realisation of a (pseudo)diagram is called a picture if it is an arbitrary
choice among a continuum of such realisations of the same diagram. 
(The existence of such realisations --- 
faithful but non-canonical {\em drawings} of the diagram --- 
is at the heart of the use of the word diagram to describe these objects.)
Suppose we have a subset of a set of diagrams characterised by 
the existence of pictures satisfying certain properties
(such as the concrete planar embedding in definition~\ref{dplanar}). 
Then given a picture of a diagram, 
another diagram realisation is said
to be {\em isotopic} to it if they belong to a continuum of pictures
all satisfying the characterising property. 

We note as a paradigm for later reference that when restricted to
$\Dnz$ the product in (\ref{compose}) is amenable to deformation.
This is firstly because the orientation of loops becomes an invariant of
isotopy (a loop cannot be flipped without some intermediate
crossing). Thus $\delta_w=\delta_{w'}$ is only necessary if $w,w'$
related by a cyclic permutation. 
\[
\includegraphics{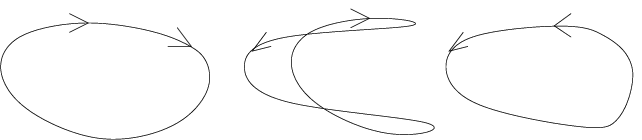}
\]
Secondly, the nonisotopic placement of loops noted in
remark~\ref{cave}
gives scope for further deformation (see section~\ref{blob'}). 


\begin{de}\label{dperiod}
(1) Consider the manifold constructed from the plane interior to the
boundary rectangle of a diagram by identifying the eastern and western
edges (a cylinder). 
A diagram in $\DVe{n}{m}$ is called {\em periodic} 
if it is possible to embed the
strings in this manifold 
(touching the northern and southern boundaries only at the vertices)
in such a way that they do not touch one another. 
\\
(2) Any specific such embedding is called a {\em concrete} periodic diagram. 
\\
(3) If one concrete periodic diagram may be continously deformed into another,
with all the intermediate stages concrete periodic diagrams, then the diagrams
are said to be {\em \aipic} \cite{Moise77}. 
\end{de}
For example, every planar diagram is periodic, while 
\[
\tau \; = \;  \{\{1,2'\},\{2,3'\},\ldots, \{n,1'\}\}
\]
is periodic but not planar. 

Write $\Dnp$ for the set of such periodic diagrams,
so that  
$\Dnz \subset \Dnp \subset \DVe{n}{n}$. 


\toappendix{\stuffbb}{{ 
\section{On constrained isotopy} \label{stuffbb} 
It may be helpful to elaborate on the meaning of isotopy in
Definition~\ref{dperiod}. For planar diagrams, once the vertices are
labelled on the frame there is no isotopy that can obfuscate this
order. The precise location of individual vertices on the frame is not
of any concern in converting between (non-unique) concrete diagrams
and their (unique) underlying diagrams.
For periodic diagrams, if any movement on the frame is allowed then
it might seem that 
isotopy can untwist the full twist (on the identity element for
example).
However, the intermediate objects in the associated
continuum would be formally ill-defined as concrete diagrams 
(in that it would not be possible to regard them
as concrete realisations of the original diagram, or indeed of any
particular diagram). 
It follows that the `natural' embedding of the identity element is not
isotopic to a twisted one in our definition of isotopy.
A more general algebra arises, therefore, if we consider 
the fundamental objects to be (frame-fixing)
isotopy classes of concrete periodic diagrams, than if we regard the
corresponding diagrams as fundamental.
}}

Proposition~\ref{post3} notes that 
two concrete planar diagrams are isotopic if and only if they are
realisations of the same diagram. 
This is not true in general for periodic diagrams. 
In particular there are non-isotopic embeddings of the identity
diagram. 
The picture of $\tau^n$ obtained by composing pictures is an embedding of
$\One$ not isotopic\ to the obvious embedding, for example. 
See appendix~\ref{stuffbb} for a fuller discussion.

Let $\Dnpc$ denote the set of periodic isotopy classes of concrete diagrams
associated to $\Dnp$ (hence a set of periodic TL diagrams). 


Because of the non-isotopic embedding possibility mentioned above the
set $\Dnpc$ is larger than the underlying set of diagrams, 
as defined by their serial realisations. 
It is possible (and useful) to have a serial realisation of isotopy classes,
however. 
Note that embedded periodic diagrams may be drawn as period-$n$ periodically
repeating {\em planar} diagrams in the infinite frame
(a string coming out of vertex 1 produces strings coming out of all
vertices congruent to 1 modulo $n$, and so on),
and hence as a subset of $\Jp{n}^{\!\! S}$:
\[
\Dnpc \hookrightarrow \Jp{n}^{\!\! S} . 
\] 
If we use {\em this} labelling in writing down the serial realisation of
diagrams (and treat vertices on different `sheets' as distinct, even
if they are congruent), then we recover the situation that concrete
periodic diagrams are isotopic if and only if they are realisations of
the same diagram.
(There are still infinitely many such diagrams however, even without beads. 
Our set of examples
$\{\{\{1,m'\}\} \; | \; m \in \Z \}$ 
are all in $\Dnpc$ with $n={1}$.)

The restriction of diagram composition (\ref{compose})
 to $\Dnp$ closes on $\Ring\Dnp$.
This product is amenable to deformation through the non-isotopic
 embeddings mentioned above 
(see \cite{MartinSaleur93,FanGreen99} for example). 
In particular let $\Dppc{n}$ denote the augmentation of $\Dnpc$ 
in case $S=\emptyset$ by classes of concrete diagrams 
including (non-crossing)  
non-contractible loops in the manner of $\Jpp{n}$. 
(NB, There is a drawing error in Fig.10 of \cite{MartinSaleur93}.
This diagram should have 4 non-contractible loops, not 3.)
 
\section{The blob algebra $b_n$ and the \achiral\ algebra $b_n'$} \label{blob'}


\subsection{Two-coloured diagrams}
Concrete TL diagrams may be thought of as partitions of
the plane interior to the frame rectangle, with the lines being the
boundaries of parts. 
The non-crossing rule means that these diagrams may be two-coloured
(in the four colour theorem sense).
For the sake of definiteness let us say that the part whose closure
includes the interval on the frame between the vertices 1 and 2 is
coloured black. 
It follows that the part between 1' and 2' is also black, 
and that concatenation preserves this canonical colouring. 
It also follows that closed loops formed in concatenation may have
either a black or a white interior (or to be more precise, immediate
interior, since they may be nested). 

Note that although colouring requires embedding, which is
noncanonical, the number of loops of each colour formed in
concatenation is invariant under plane \aipy\ 
(note that this is not true for the larger classes of Brauer \aipy). 
It follows that we can generalise the composition rule in
(\ref{compose}) by generalising $k_{a,b}$. 

By rescaling the generators $U_i$ one can see that 
this generalisation is isomorphic to the ordinary case 
(excepting the specialisation in which one of the parameters is not
invertible, which is not appropriate for the context of computation
for Potts and vertex models, where the algebra has its origin 
\cite{Baxter82}). 
However, we shall now show that it leads the way to some 
further rather more useful generalisations 
(in the spirit of the blob algebra 
\cite{MartinSaleur94a} viewed as a generalisation of 
the type-B algebra of \cite{tomDieck94}).


\subsection{Subalgebras and deformations} \label{deform2}
Let $\DB{n} \subset \Dzo{n}$ be 
the subset of TL diagrams that are (\aipic\ to concrete diagrams that
are) invariant under reflection in a central vertical line.
Putting aside for a moment any closed loops that might arise,
this subset is obviously fixed under concatenation. 
That is, the concatenation of two concrete representatives of elements
of $\DB{n}$ is a representative of an element of $\DB{n}$
(ignoring closed loops).

\noindent {\em Remark:}
Indeed we could consider a version of $\DB{n}$ that is not a strict subset of
$\Dzo{n}$ (a set whose elements are isotopy classes, where even
classes containing symmetric elements
include elements which are not concretely symmetric) but to be such
that the elements are `symmetric isotopy' classes - i.e. classes whose
elements maintain exact symmetry. 
So long as we define our algebra composition without reference to
pseudodiagrams this distinction is academic
 (but see section~\ref{reduc}).

If $n$ is odd then $\DB{n}$ is a rather uninteresting subset.
If $n$ is even, redefine canonical colouring to be that in which
the part whose closure includes the central northern interval
is coloured white. 
Note that within this subset the property that the concatenation of
two diagrams forms a loop that {\em crosses} the centre line
(indeed white,
or black, loop that crosses) 
is invariant under \aipy. 
Thus we may deform the algebra spanned by these diagrams by 
generalising the scalar factor $k_{d_1 d_2}$ so that 
\eql(new k)
k_{d_1 d_2} = \delta^{l_0} \delta_e^{l_w} \kappa^{l_b}
\eq
where $\delta, \delta_e, \kappa \in \Ring$, $l_0$ is the number of
(pairs of) noncrossing loops, $l_w$ is the number of white crossing
loops and $l_b$ the number of black crossing loops. 
\del(b')
Let $b'_n(\delta, \delta_e, \kappa)$ denote the algebra which is 
$\Ring \DB{n}$ with product defined by (\ref{new k}). 
\end{de}


As usual the identity in this algebra is the pair partition $\Id$.
For any given $n=2m$ let $\abe'$ denote the diagram corresponding to the
pair partition differing from 1 only in $\{m,m+1\},\{m',m+1'\}$,
that is (in case $m=3$):
\[ 
\abe' = \; 
\raisebox{-.21in}{
\includegraphics{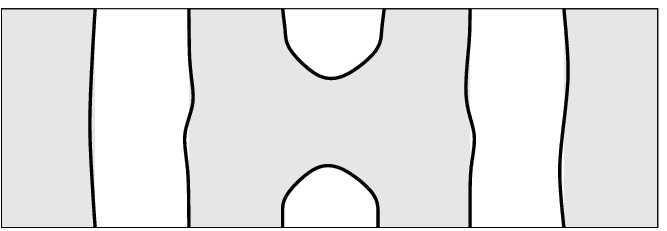}}
\]

\subsection{The unfolding map $\mu$}

\begin{de}
An {\em \pseudo-diagram} is a TL diagram in which strings may end on
the east or west edge of the rectangle as well as north and south.
\end{de}

Recall that $B_n$ denotes the set of blob diagrams. 
Given a blob diagram $d$, we may define a \pseudo-diagram from it by 
deforming every string with a blob until an arc in the immediate
neighbourhood of the blob is just on the outside of the western edge
of the rectangle, and then discarding this arc. 
If we compose this \pseudo-diagram with its mirror
image in the western edge 
(NB, this is a well defined construct on isotopy classes)
we have a left-right symmetric TL diagram.
\newcommand{\Mu}{{\Large \mu}}%
Let us call it $\Mu(d)$. 
For example $\Mu(e)=\abe'$. 
\prl(blob base)
The map $\Mu:B_m \longrightarrow \Dzo{2m}$ defined above is an injection,
and the range is the set $\DB{2m}$ of left-right symmetric  diagrams.
\Qed
\end{pr}
(This generalises the combinatorial map described in
\cite{CoxGrahamMartin03}.) 
\prl(blob iso)
The map $\Mu$ extends to an algebra homomorphism, so that 
the algebra $b'_{2m}(\delta, \delta_e, \kappa)$
is isomorphic to the blob algebra $b_m(\delta, \delta_e, \kappa)$.
\end{pr}
{\em Proof:} We need to check that $\Mu(a) \Mu(b) = \Mu(ab)$,
where $a,b$ are blob diagrams. If we consider $ab$ as a concatenation,
without (for a moment) imposing the blob relations, 
it can have two blobs on the same line, and it can have closed loops,
with and without blobs. The map $\mu$ makes sense on such an $ab$, and
commutes with concatenation. It is thus necessary to check that the
imposition of the blob relations on $ab$ produces the same factor as 
$k_{\Mu(a), \Mu(b)}$.
If two blobs come together we use the blob relation $ee=\delta_e e$
from (\ref{ee=e}).
The image under $\Mu$ is (locally): 
\[
\includegraphics{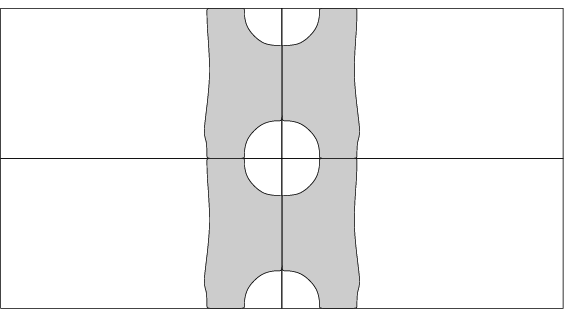}
\]
which has one white loop, so $k_{\Mu(a), \Mu(b)}=\delta_e$ as
required. 
A decorated loop is replaced by a factor $\kappa$ in the blob algebra
$b_m(\delta, \delta_e, \kappa)$ (cf. section~\ref{dir homs}),
and passes to a black loop under $\Mu$. 
From (\ref{new k}) this again gives a factor $\kappa$. 
An undecorated
loop ($\delta$ in the blob algebra) passes to a pair of off-axis loops
under $\Mu$.
\Qed


\section{Recollement and $b_n'$ representation theory} \label{recol b}
\begin{figure}
\[
\includegraphics{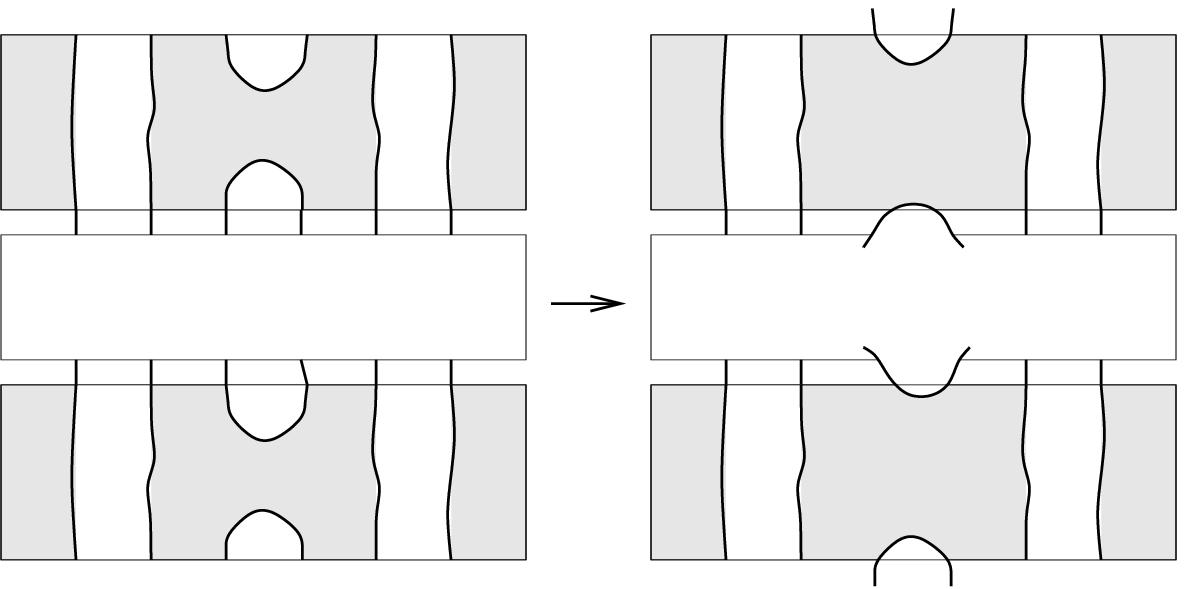}
\]
\caption{\label{foobre} Removing the central cup and cap.}
\end{figure}

Suppose that $\delta_e$ is invertible. 
Then the subset of $\DB{2m}$ consisting of 
symmetric diagrams in which the inner central `cup' and `cap' appear
(i.e. diagrams that contain 
$\{m,m+1\},\{ m',m+1' \}$
as pair partitions) are a basis for a subalgebra of $b'_{2m}$.
Since $\abe' \abe' = \delta_e \abe'$ then $\frac{1}{\delta_e} \abe'$ 
is idempotent,
and is the unit in this subalgebra.
It may be identified with the idempotent subalgebra  
$b''_{2m} = \frac{\abe'}{\delta_e} b'_{2m}(\delta,\delta_e,\kappa)
 \frac{\abe'}{\delta_e}$.
Since a white loop appears automatically in every composition in this
algebra, a better basis is the set 
$\DBe{2m}$ 
of such diagrams each multiplied by 
$\frac{1}{\delta_e}$. 
Then if no {\em other} loops appear in composition the
product of two basis elements is another basis element.%
\\
{\em Remark:}
Actually the algebras $b'_{2m}(\delta,\delta_e,\kappa)$
  and  $b'_{2m}(\delta,\alpha \delta_e,  \alpha \kappa)$ are readily seen to be
  isomorphic for any invertible $\alpha$ 
(consider the isomorphic \cite{Martin95pres}
presentational form (\ref{TL001}-\ref{TL006}), for example).
It is thus possible to replace 
  $b'_{2m}(\delta,\delta_e,\kappa)$ with  
  $b'_{2m}(\delta,1,\delta_e^{-1}\kappa)$ without loss of generality.
\medskip


It will be evident from figure~\ref{foobre} that 
there is 
a set map
\[
\rho^- : \DBe{2m} \rightarrow \DB{2m-2}
\]
defined by simply removing the central upper `cup' and lower `cap' from
the diagram underlying each basis element (and discarding the factor 
$\frac{1}{\delta_e}$). Indeed $\rho^-$ is a bijection. 
We next extend this to an algebra map. 


\prl(recool)
The map $\rho^-$ extends $K$--linearly to an algebra isomorphism 
$$
\rho^-: 
\frac{\abe'}{\delta_e} b'_{2m}(\delta,\delta_e,\kappa) \frac{\abe'}{\delta_e} 
   \longrightarrow b'_{2m-2}(\delta,\kappa,\delta_e) .
$$ 
\end{pr}
{\em Proof:}
Consider $\rho^-$ extended in the obvious way to pseudodiagrams.
Concatenation of diagrams may be thought of as 
the first step in computing composition 
on both sides, and commutes with  $\rho^-$ 
(given the automatic cancellation of one loop with a normalising
factor on the domain side).
The final step is interpretation of the pseudodiagram as a scalar
multiple of the underlying basis element. 
This differs on the range side, 
precisely in that the colour assigned to each loop is reversed 
(by the cup/cap removal).
\footnote{ 
Consider the colouring of a pseudodiagram 
on the domain side of $\rho^-$, which
determines the scalar factors there.
If we think of trying to keep this colouring, then
once the cup and cap are removed the colour of
the upper central interval (as used to determine scalar factors)
is black, not white.
Obviously then, {\em all} colours are inverted,
compared to what they would normally be on the target side.  
Thus in particular central loops which would properly be 
white will be black, and vice versa. 
Thus the map $\rho^-$ reverses all colours
and will not extend naively to an algebra homomorphism. 
Swapping the colours exchanges the roles of $\delta_e$ and $\kappa$,
so we can fix this by defining the extension
as above for basis elements, 
but which maps scalar coefficients by exchanging the 
roles of $\delta_e$ and $\kappa$.
}
This difference is thus itself reversed by exchanging the roles of 
 $\delta_e$ and $\kappa$. 
\\
\Qed


For example, let $a$ be a basis element with underlying diagram
as shown in the upper left of:
\[
\includegraphics{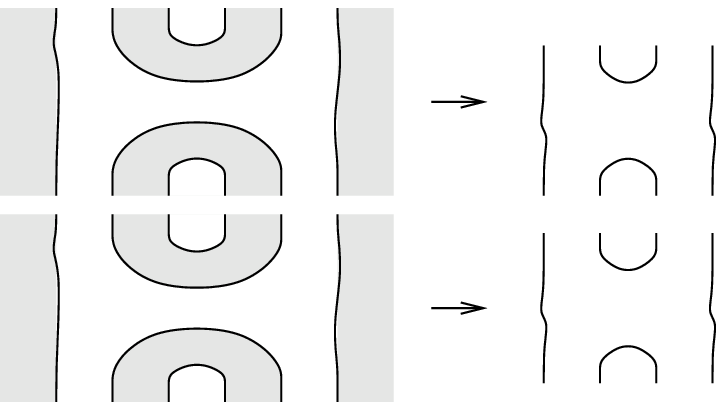}
\]
We see that the product on the left is $a.a = \kappa a$
(taking account of factors of $\frac{1}{\delta_e}$). 
The image $\rho(a)$ is shown on the right, together with 
$\rho(a).\rho(a) = \tilde \delta_e \rho(a)$ 
(writing $\tilde \delta_e$ for the $\delta_e$ parameter in the image). 
Thus in order for $\rho(a.a)=\rho(a).\rho(a)$ we require 
$\tilde \delta_e = \kappa$. 

On the other hand, with the basis elements indicated on the left in:
\[
\includegraphics{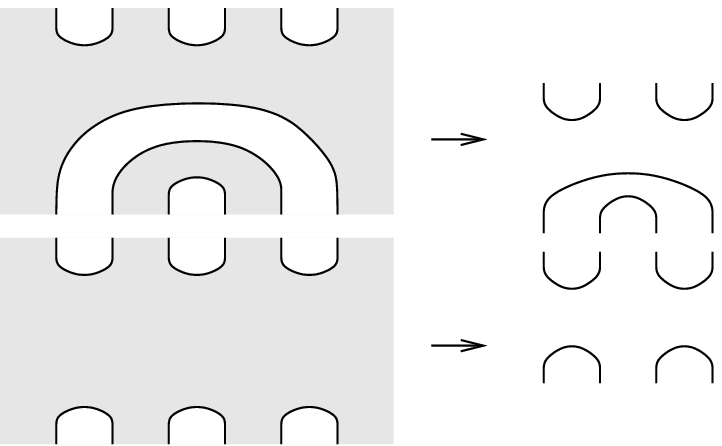}
\]
we have $b.c = \delta_e c$, while on the left 
$\rho(b).\rho(c) = \tilde \kappa \rho(c)$. 
So we require $\tilde \kappa = \delta_e$, as stated. 


For the moment let us use the shorthand $A_m$ for 
$b'_{2m}(\delta,\delta_e,\kappa)$ and $B_{m}$ for
$b'_{2m}(\delta,\kappa,\delta_e)$,
and $\Lambda^A_m$ for an index set for simple modules of $A_m$
(and similarly for $B_m$). 
We may now bring to bear the recollement part of the machinery 
described in \cite{CoxMartinParkerXi03} 
(summarized in section~\ref{cat}). 
In particular note that
Proposition~\ref{recool} tells us the following. 
\begin{theo} \label{b' fullembed}  
The category of left
$b'_{2m-2}(\delta,\delta_e,\kappa)$-modules 
is {\em fully embedded} in the category of left
$b'_{2m}(\delta,\kappa,\delta_e)$-modules. 
\Qed
\end{theo}   
It follows (via proposition~\ref{grebo}) that the simple modules $S$
of the latter {\em not} obeying $\abe'S=0$ may be indexed by 
$\Lambda^B_{m-1}$. 
The simple modules obeying  $\abe'S=0$ are also simple modules of the
quotient
$A_m/A_m \abe' A_m$, but it is easy to see that this has precisely one
simple module
(since $A_m \abe' A_m$ contains every diagram with fewer than $2m$ propagating lines). 

In this way we derive the (well known) index set for simple modules of
$A_m$. We also derive some striking results relating the homomorphisms
between standard modules for $A_m$ and $B_{m-1}$ which, while not
revealing any homomorphisms which were not already known, do reveal a
layer of symmetry in the organisation of these homomorphisms which has
not been noted before. 
This structure is not needed in the analysis of the blob algebra's
representation theory, but it raises the very intriguing possibility
of similar symmetries in (affine) Hecke algebra representation theory.
We will discuss this further elsewhere. 

\subsection{The $b_n$ version} \label{remi}

With the benefit of hindsight we see that the recollement can be
invoked directly in $b_n$. 
Let $B_n'$ be the subset of $B_n$ consisting of elements in which both
the string containing vertex 1 and that containing vertex $1'$ are decorated
(of course this could be the same string). 
Let $B^e_n = \{ \frac{1}{\delta_e} d \; | \; d \in B_n' \} \; \subset b_n   $.

Define a map   
\eql(rho_1)
\rho_1: B^e_n \longrightarrow B_{n-1}
\eq
as follows. 
For $d \in B_n'$,  
consider the region of $d$ with the western edge in its closure:
we have a sequence of one or more decorated lines reading clockwise around
this region, starting from the vertex 1
(ignore the undecorated ones).  
This sequence is of the general form 
$\{ 1 , i_1 \}, \{i_2,i_3\},\ldots, \{ i_l, 1'  \}$ (some possibly
primed vertices $i_1,\ldots,i_l$), 
or simply  $\{ 1, 1' \}$. 
In the latter case, simply erase the (decorated) line $\{ 1, 1' \}$. 
Otherwise, erase the sequence and replace with 
decorated lines $\{i_1,i_2 \}, \ldots, \{ i_{l-1}, i_l \}$. 
After suitable renumbering we have an element of $B_{n-1}$. 
This is $\rho_1( \frac{1}{\delta_e} d)$. 

The map is illustrated by the following example:
\[
\includegraphics{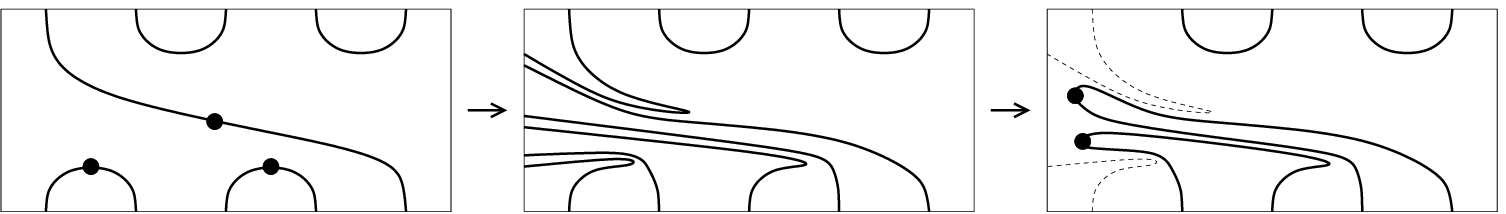}
\]
It will be evident that this map is a bijection. 


\prl(b case)  
The map $\rho_1$ extends $K$--linearly
to an algebra isomorphism from the subalgebra $b^e_n$ 
of $b_n(\delta,\delta_e,\kappa)$ 
spanned by $B^e_n$ (with unit $\frac{1}{\delta_e} e$) to 
$b_{n-1}(\delta,\kappa,\delta_e)$.
That is, $\rho_1$ is an algebra isomorphism with $\delta_e$
and $\kappa$  interchanged. 
\end{pr} 
{\em Proof:} Consider Proposition~\ref{blob iso} and 
Proposition~\ref{recool}. 
A short manipulation of diagrams (in the diagram algebra
sense) confirms that the following diagram is commutative:
\[
\xymatrix{
b^e_{n}   \ar[d]_{\rho_1} \ar[r]^{\mu} & b''_{2n} \ar[d]^{\rho^-} \\
b_{n-1}        \ar[r]_{\mu}          & b'_{2n-2}
}
\]
for the appropriate parameter values. 
\Qed
\\
The parameter change can be seen  
directly 
by considering multiplications in low $n$ cases.
For example:
\eql(pc1)
\xymatrix@R=1pt%
@C=25pt@M=0pt{
\includegraphics{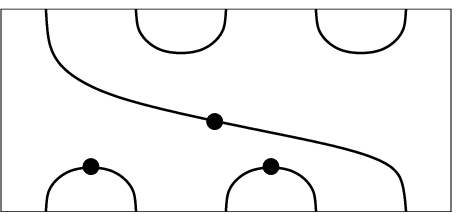} 
\raisebox{.21in}{$\rightarrow$}
\includegraphics{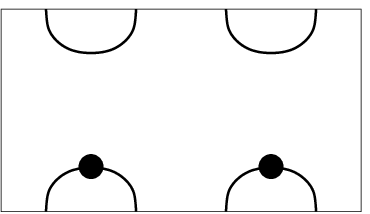}
\\
\includegraphics{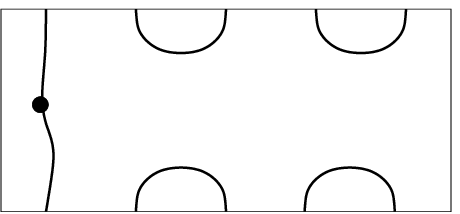} 
\raisebox{.21in}{$\rightarrow$}
\includegraphics{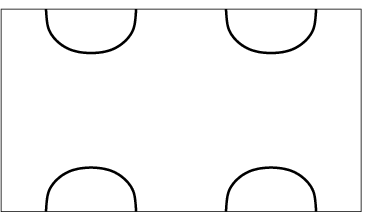}
}
\eq
To implement the proposed normalisation, let us call the top left
diagram $d_1$ and the bottom left $d_2$. Then on the left,
implementing the normalisation, we want to consider 
$\frac{d_1}{\delta_e} \frac{d_2}{\delta_e}
= \frac{\delta_e^3}{\delta_e^2} e U_2 U_4
= \delta_e^2 \frac{e U_2 U_4}{\delta_e} $.
The figure shows that the left-hand side of this identity passes to 
$\kappa^2 U_1 U_3$ under $\rho_1$, while the right-hand side passes
directly to $ \delta_e^2 U_1 U_3$. 
Here we see that, allowing for the normalisation of the first blob as
an idempotent, on the left we pick up a factor $\delta_e^2$, 
and on the right a factor $\kappa^2$. 
Meanwhile
\eql(pc2)
\xymatrix@R=1pt%
@C=25pt@M=0pt{
\includegraphics{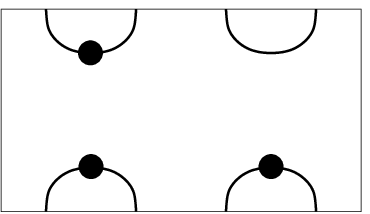} 
\raisebox{.21in}{$\rightarrow$}
\includegraphics{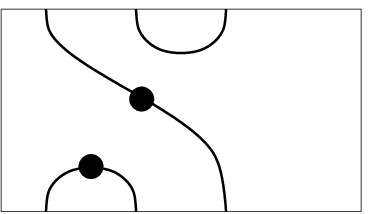}
\\
\includegraphics{xfig/eTLB155.eps} 
\raisebox{.21in}{$\rightarrow$}
\includegraphics{xfig/eTLB156.eps}
}
\eq
has $\kappa^2$ on the left but $\delta_e^2$ on the right.


\section{\Achiralb\ algebras} \label{ASTLA}
In this section we define several new algebras with the same flavour
as $b_n$, and relate them to other algebras studied in the literature
(cf. \cite{OrellanaRam01,NicholsRittenbergdeGier05}).
\\
The first step is to define an appropriate class of pseudodiagrams,
which compose by concatenation. Then we define
a reduction of pseudodiagrams into scalar multiples of basic diagrams,
so that composition reduces to an algebra multiplication. 

\subsection{Pseudodiagram categories}\label{pdcat}

A TL {\em pseudodiagram} is a TL diagram possibly including closed loops. 
NB, loops cannot move isotopically over lines, so the set of $(n,m)$
TL pseudodiagrams is larger than the subset 
of $\DoVe{n}{m}$ with no beads. 

\del(pseud)
{\em \cite[\S2.3]{MartinWoodcock2000}} 
A (left) {\em blob  pseudodiagram} is a TL pseudodiagram in which any
left 0-covered arc may be decorated with a (left-)blob.  
\end{de}
There is a corresponding notion of right blob pseudodiagrams. 
A {\em left-right blob pseudodiagram} is a pseudodiagram which may have 
left and right-blob decorations, so long as every decorated arc may be
deformed to touch its appropriate edge {\em simultaneously}. 
Write $\pseud(n,m)$ for the set of $(n,m)$-pseudodiagrams of this
kind. 


Provided they have the right number of vertices,
these planar 
pseudodiagrams may be composed by extending the usual diagram concatenation
(\ref{frubeddd}).
That is, concatenate concrete representations to give a concrete
representative of the composite.
\begin{lem} \label{weeee}
This composite is well defined
(i.e. independent of the choice of representatives),
associative and unital. 
\end{lem}
{\em Proof:}
The argument of Proposition~\ref{abacus} is not affected by the need to take
account of isotopy classes distinguished by the embedding of closed loops.
\Qed \\
Further
\prl(logic)
The composition of an $(n,l)$-- and an $(l,m)$--pseudodiagram of the same
type (ordinary, blob, left-right blob) is an $(n,m)$--pseudodiagram 
of that type. 
For each type the triple $(\N,\pseud,\circ)$
is a category, where $\pseud(n,m)$ is the set of morphisms for 
$n,m \in\N$.
\end{pr}
Note that there are unboundedly many pseudodiagrams of each type. 
Various features can appear (repeatedly) in pseudodiagrams, such as: 
\newline
($\delta$) undecorated loops;
\newline
($\delta_L$, $\delta_R$) 
consecutive runs of two left or right-blobs on the same arc ($LL$ or $RR$);
\newline
($\kappa_L$, $\kappa_R$) loops decorated with a left or right-blob;
\newline
($\kappa_{LR}$) loops decorated with a left and a right-blob 
(NB, LRLR sequences are not possible on loops, so it is always possible to
arrange blobs on loops into at most two same-type runs);
\newline
($k_L$) in a pseudodiagram with unique propagating line; this line can have 
LRL (resp. RLR) sequences on it. 
\newline
It will be convenient to be able to refer to these features by
the set of shorthand names indicated:
 $P_B \; := \; \{\delta, \delta_L, \delta_R, \kappa_L, \kappa_R,
\kappa_{LR}, k_L \}$. 

\prl(fin of)
For given $(n,m)$, there are only finitely many 
diagrams with {\em none} of the  features in $P_B$. 
\end{pr}
{\em Proof:} Note that there are only finitely many underlying
(undecorated) TL
shapes possible, since all possible decorated loops are excluded.
Then only finitely many decorations of the
lines in these loop free shapes remain. 
\Qed

For example in case $(1,1)$ every diagram has underlying TL diagram
with just one string. The $R,L$-word on this string is one of 
$\{ \emptyset, L, R, LR, RL \}$. 


\del(Bx')
A left-right blob diagram is a pseudodiagram in which none of the
features above occur. 
The set of left-right blob diagrams with 
$m$ vertices on the northern edge and $m$ vertices on the southern edge 
is denoted $B^{x'}_{m}$. 
\end{de}
This set $B^{x'}_{m}$ 
is thus the set of diagrams whose underlying TL diagram has no
loops, and where each western exposed line may be decorated with at
most one left blob, and each eastern exposed line may be decorated
with at most one right blob, subject to the condition that it must be
possible to deform each blob to its appropriate edge simultaneously
without lines crossing.
Each diagram has a representation as a pair partition with
decorations, much as in the $b_n$ case except that decorations $R$ and
even $LR$ and $RL$ may be possible. For example
\[
\{ \{1,2 \}_{LR} ,  \; \{1',2' \}_{LR} \} \in B^{x'}_{2}
\]
More generally, focusing on $L$-decorated pairs, the simultaneous
deformation requirement gives the following. 
\prl(LR)
If $d \in B^{x'}_{m}$ and $i,j \in \{1,\ldots,m\}$ then:
\\ 
$(i,j)_{LR}$ a pair part in $d$ implies $i<j$ are the two largest
unprimed numbers in the list of $L$-decorated pairs;
\\ 
$(i,j')_{LR}$ or $(i,j')_{RL}$
a pair part in $d$ implies $i$ (\resp\ $j'$) is the largest
unprimed (\resp\ primed) number in the list of $L$-decorated pairs;
\\ 
$(i',j')_{LR}$ a pair part in $d$ implies $i<j$ are the two largest
primed numbers in the list of $L$-decorated pairs.
\end{pr}
Thus
\begin{co} \label{ok}
At most two $R$s can appear in the list of $L$-decorated pairs
in $d \in B^{x'}_{m}$,
and then precisely in the situation of the following figure.
\end{co}
\eql(tq RHS)
\includegraphics{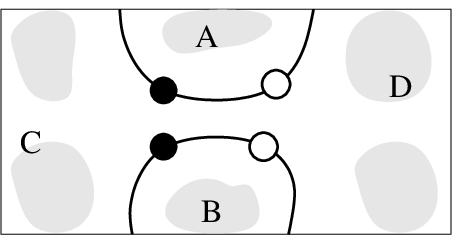}
\eq
Note that in pictures we use a solid blob for $L$ and $\circ$ for
$R$. 

\del(Bxm)
The set $B^{x}_{m}$ is obtained from $B^{x'}_m$ by discarding the
diagrams of the type shown  
in (\ref{tq RHS}).  
\end{de}

\subsection{Initial pseudodiagram reduction}

\newcommand{\mysim}[1]{\stackrel{#1}{\leadsto}}%
Let $d$ be a pseudodiagram. Write 
$$
d \mysim{x} d'
\qquad 
(x \in P_B)  
$$
if $d,d'$ differ by removal of one corresponding loop 
or, respectively, $L,R$-string replacement: 
$LL \leadsto L$, 
$RR \leadsto R$, 
$LRL \leadsto L$,
$RLR  \leadsto R$. 


It will be evident that all maximal chains of relations
$\mysim{\delta}$ starting from $d$ end in a diagram with no
$\delta$-loops, and are of the same length --- call this length
$\hash_{\delta}(d)$; 
and that $\hash_{\delta_L}(d)$, $\hash_{\delta_R}(d)$ may be defined
similarly. 
Indeed 
\prl(reducer)
For any $d$ there is always a chain of relations, with each relation some 
$\mysim{x}$ ($x \in P_B$), 
ending in a diagram with {\em none} of the identified features. 
There are in general multiple such chains from $d$, but every 
one ends in the same `reduced' diagram --- call it $r(d)$.  
Each such chain for $d$ has the same number of links of the form $\mysim{x}$
for {\em given} $x \in P_B$. 
\end{pr}
This number of links defines $\hash_x(d)$ for each $x \in P_B$. 
\\
{\em Proof:} Note that the reductions are of two types: those that
shorten the sequence of decorations on some arc of some line; and
those that remove a loop. 
Both types are localised to individual lines 
(leaving all other structure unchanged), 
so we may talk of individual lines as being `locally' reduced. 
Note also that loop removals only apply to
loops that are locally reduced. 
Thus we may consider the reduction of each individual line.  
On each line we have a simple Bergman diamond 
\cite{Bergman} for the reduction 
of sequences of Ls and Rs.
The only ambiguity is the reduction of LRLR to LR via LRL or RLR
replacement, but both of these has $x=k_L$.

Once line-local reductions are complete, the loop removal process is
immediate. 
\Qed

\subsection{Towards finite and localisable algebras}

Note that if $\hash(a)$ is the number of occurences of any given one
of the features $P_B$ in pseudodiagram $a$, then 
\[
\hash(ab) \geq \hash(a) + \hash(b)
\]
With this in mind, we can set out to define a finite dimensional
algebra of diagrams (with fixed number of northern and southern
vertices) by applying a quotient rule which equivalences a
pseudodiagram with some such feature with a scalar multiple of the
same diagram but with this feature 
excised or replaced as in $\mysim{x}$. 
Thus
\eql(LL)
L L = \delta_L L
\eq
\eql(LRL)
LRL=k_L L  
\eq
and so on, replacing each feature in $P_B$ with 
a correspondingly named scalar:
\eql(Loop)
d=k_d r(d)
\eq
where $k_d = \prod_{x \in P_B} x^{\hash_x(d)}$.
(Note that this is not a unique finitising procedure 
--- for example we could omit
the excision of $RLR$, to leave a slightly larger but still finite
quotient.)

The set $B^{x'}_{m}$ is a basis for the proposed 
finite-dimensional quotient algebra.
To see that this quotient is internally consistent 
note that composition proceeds by concatenation to produce a
pseudodiagram, which is then reduced using the relations.
This reduction is consistent by Proposition~\ref{reducer}. 


Ab initio one might try to assign a different scalar parameter to 
the given LRL and RLR reductions, 
but they are related by commuting diagram:
\[
\xymatrix{
& LRLR \ar@{->}[dl] \ar@{->}[dr]_{RLR=k_R R}
\\
k_L LR  \ar@{=}[rr] & & k_R LR
}
\] 

Note that it is not possible to reduce the occurrence of LRL (or RLR)
in a loop using $k_L$ ($=k_R$) since this requires that the LRL reside
on the only propagating line, so there is no return route to complete
the loop. 
Indeed $k_L$ never occurs in even index 
(even $m$) 
algebras,
while $\kappa_{LR}$ {\em only} occurs in even index algebras 
(no propagating line can pass either to the left or to the right of
the loop). 
It will turn out to be appropriate to set 
\eql(KLR=KL)
\kappa_{LR} = k_L
\eq
as can thus be done without loss of generality. 

\newcommand{\abparams}{\delta,\delta_L,\delta_R,\kappa_L,\kappa_R,\kappa_{LR}}%

At this stage we define  algebra
$b^{x'}_{m}(\abparams)$
to be the quotient of the linear extension of the pseudodiagram
composition by the relations  (\ref{LL}), (\ref{LRL}) and so on
associated to the $\mysim{x}$ relations. 
The following is clear.  

\prl(blob is sub)
There is an algebra monomorphism 
$$
b_{m}(\delta,\delta_L,\kappa_L) \hookrightarrow b^{x'}_{m}(\abparams)
$$
that takes $e \mapsto L_1$. \Qed
\end{pr}

We want to define an algebra (the {\em \achiralb\ algebra}) 
that will be the quotient of $b^{x'}_m$ by some 
small set of further relations, 
chosen so as to make the representation
theory of the algebra tractable by localisability.
(These relations will also turn out  
to make contact with some other interesting algebras; 
see section~\ref{ASTLA'}.) 
That is, following the remarks in section~\ref{remi} we may 
determine the structure of such an algebra at level $m$ by constructing an
idempotent subalgebra isomorphic to some known algebra 
(i.e. some version of the level $m-1$ algebra, known by inductive
hypothesis)
--- the localisation of the algebra under consideration. 

\medskip

Before passing to a localisable quotient, 
note that $b_n^{x'}$ is a quotient of the affine-$C$ Hecke algebra 
defined in section~\ref{INot}: 
\prl(hecke C)
The map $u_i \mapsto U_i$ ($n>i>0$), 
$u_0 \mapsto e$ (left-blob), 
$u_n \mapsto f$ (right-blob),
extends to an  algebra homomorphism $\phi$ from  $T(\hat{C}_n)(q_0,q_1,q_x)$
to $b_n^{x'}$ in case 
$q=q_1$, 
$\delta_L= q_0+q_0^{-1}$,
$\delta_R= q_x+q_x^{-1}$,
$\kappa_L = \frac{q_0^2 + q^2}{q_0 q}$,
and $\kappa_R = \frac{q_x^2 + q^2}{q_x q}$. 
\end{pr}
{\em Proof:} The relation checking is largely routine. 
Note from (\ref{uuuu}) that 
$\phi(u_1u_0u_1 -  \frac{q_0^2 + q_1^2}{q_0 q_1} u_1)$ vanishing is
sufficient to ensure 
$\phi(E_{B_2})=
\phi(u_0(u_1u_0u_1 -  \frac{q_0^2 + q_1^2}{q_0 q_1} u_1))=0$,
and similarly with $u_0,q_0$ replaced by $u_n,q_x$. 
\Qed

\medskip

This is closely analogous to the connection between the blob algebra
and type-$B$ Hecke algebra (Proposition~\ref{hecke}), which has 
allowed the localisability of the {\em blob} to be used to 
investigate the representation theory of that Hecke algebra 
\cite{CoxGrahamMartin03}. 

\subsection{Combinatorial localisation}

By analogy with  section~\ref{remi}, 
we may form the subset $B^{x'e}_{m}$ of diagrams 
(normalised by $1/\delta_e = 1/\delta_L$ as before) 
in which the
string(s) involving vertices 1 and $1'$ are decorated (with $L$).
\\
The set $B^{xe}_{m}$ is the subset of $B^{x'e}_{m}$ 
whose underlying diagrams are elements of  
$B^{x}_{m}$. 
\medskip 

Our relations imply that this subset  $B^{x'e}_{m}$   spans a subalgebra. 
The question is:
{\em How do we construct an isomorphism between this subalgebra 
(or a `large' quotient) and the algebra at level $m-1$?}
(Note that this choice of idempotent subalgebra 
breaks the left-right symmetry of the algebra
--- so there is a corresponding R--based formulation.)

To address this we first construct set maps between these algebras'
basis diagrams. 


\del(def pre tilde rho)
Let $\pseud(m,m)^e$ denote the subset of $\pseud(m,m)$ 
in which the string(s) involving vertices 1 and $1'$ are $L$-decorated,
and $\pseud_o(m,m)$ the subset with no loops and no
multiple $L$-decorations on the same segment, 
and $\pseud_o(m,m)^c = \pseud_o(m,m) \cap \pseud(m,m)^e$. 
\end{de}

Let $\tilde\rho$ denote the direct extension of the $\rho_1$ map 
(\ref{rho_1})  
for $b_n$ to $ B^{x'e}_{m}$.
That is 
\[
\tilde\rho : B^{x'e}_{m} \rightarrow \pseud(m-1,m-1)
\]

Examples:
\eql(top q1)
\includegraphics{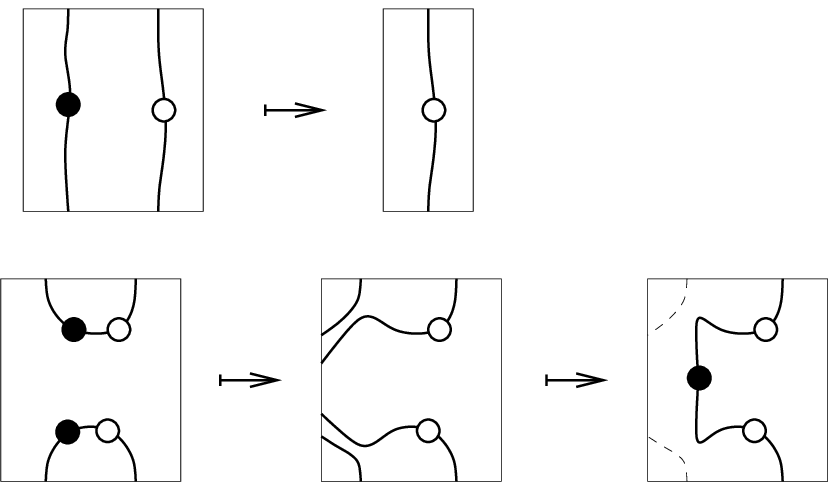}
\eq
\eql(top q1+)
\includegraphics{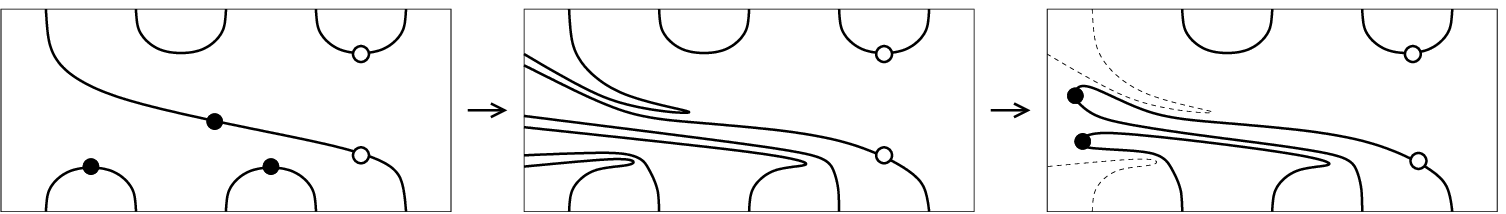}
\eq

Note that the range is not $B^{x'}_{m-1}$, because segments with 
decoration RLR are possible, as (\ref{top q1}) illustrates. 


\del(def rho*)
Let 
\[
\rho^* : B^{x}_{m-1} 
\rightarrow \{ \frac{1}{\delta_L} p \; | \; p \in \pseud(m,m)^e \}
\]
be defined as follows. If $d$ has no $L$ decorations then 
it maps to the same diagram but with a propagating $L$-decorated line
on the left. Otherwise deform all blobs to just outside the western
edge; cut the blobs off --- 
so some positive even number of lines now end at the western edge; 
deform the top and bottom-most of these lines so that they end on the
top and bottom edge respectively; 
$L$-decorate these lines; 
and reclose the remaining western endpoints with $L$-decorated arcs in 
adjacent pairs in the only possible way. 
\end{de}


\prl(p*p)
The map $\tilde\rho$ restricts to a bijection:
\[
\tilde\rho : B^{xe}_{m} \rightarrow B^{x}_{m-1}
\]
with inverse $\rho^*$. 
\end{pr}
{\em Proof:}
First note that if $d$ has $\{1,1'\}$ as $L$-decorated pair the result
$
\rho^* ( \tilde\rho (d)) = d 
$
is clear from the definitions. 

For $d \in B^{xe}_{m}$ let $\{1,i_1 \},\ldots, \{i_l,1'\}$
be the $L$-decorated pairs. 
By Corollary~\ref{ok} at most one of these is $R$-decorated
--- $\{ i_r, i_{r+1} \}$ say. 
Then $\tilde\rho(d)$ has the same
(suitably relabeled) non-$L$-decorated pairs, and the $L$-decorated
pairs $\{i_1,i_2 \},\ldots, \{i_{l-1},i_l\}$;
with the further $R$-decoration (if any) on 
$\{ i_{r-1}, i_{r} \}$ or $\{ i_{r+1}, i_{r+2} \}$
as appropriate:

\eql(move4)
\includegraphics{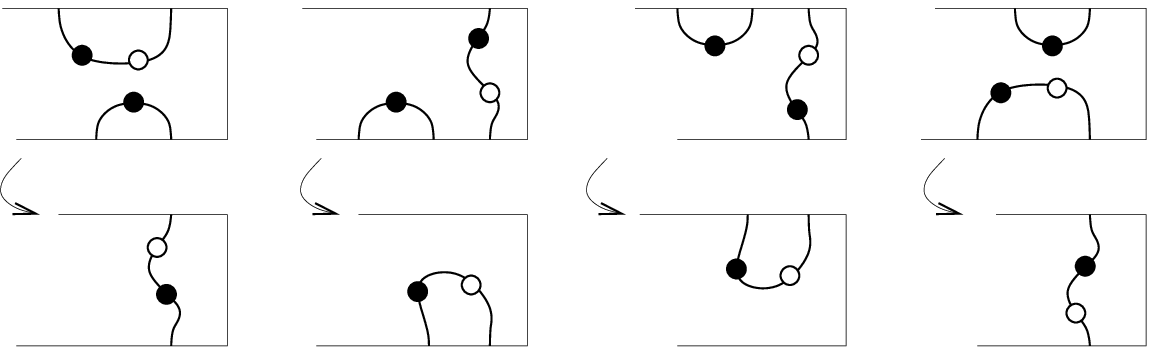}
\eq
Note that this establishes that 
$\tilde\rho(  B^{xe}_{m} ) \subseteq B^{x}_{m-1}$.

For $d \in B^{x}_{m-1}$ let $\{i_1,i_2 \},\ldots, \{i_{l-1},i_l\}$ 
be the list of $L$-decorated pairs. As before there is at most one
$R$-decoration in this list, on $\{ i_{r-1}, i_{r} \}$ say. 
Then $\rho^*(d)$ has the same non-$L$-decorated pairs as $d$, and 
 $L$-decorated pairs $\{1,i_1 \},\ldots, \{i_l,1'\}$.
It has at most one further $R$, on 
$\{ i_{r-2}, i_{r-1} \}$ or $\{ i_{r}, i_{r+1} \}$
as appropriate.
Note that it follows that 
$\rho^*(  B^{x}_{m-1} ) \subseteq B^{xe}_{m}$.

Considering $\rho^* ( \tilde\rho (d))$ then, the lists of $L$-decorated
(and undecorated) pairs are manifestly restored.
The location of the $R$-decoration (if any) in the $L$-decorated list 
may be seen to be restored by considering the four cases
in Proposition~\ref{LR} (illustrated above).

A similar argument shows the right inverse property. 
\Qed

\newcommand{\stuffa}{{
\parker{I THINK WE CAN DROP THIS NOW?:}

A pictorial version of prop~\ref{p*p} is as follows. 
\prl(claim big)
The map $\tilde\rho$ restricts to a bijection:
\[
\tilde\rho : B^{xe}_{m} \rightarrow B^{x}_{m-1}
\]
\end{pr}
\footnote{OLD VERSION:
\begin{claim}
This requirement is satisfied in general by imposing 
the `topological' quotient (\ref{top quot}). 
\end{claim}
}
{\em Proof:}
First we shall show that this is a map --- i.e. that the range may be
restricted as indicated. 
To begin, consider the range of $\tilde\rho$ as originally defined. 
We shall argue that no `reducible component' 
(e.g. closed loop, or LL on the same line segment) appears in any
pseudodiagram in the range, except possibly RLR 
(whose appearance is then excluded by the restriction we have now imposed).

The restriction to $B^e_{n}$ is $\rho_1$, with range $B_{n-1}$ which,
note, has no closed loops.
Every diagram in $B^{x'e}_{n}$ is obtained from one in $B^e_{n}$ by
adding some R decorations, but R decorations do not affect the
underlying shape change produced by $\tilde\rho$, so:
\\
1. the range of $\tilde\rho$ includes no diagram with loops.
\\
In order for a component of form LL to appear as an image, 
the diagram in its
preimage would have to contain a loop, by the definition of $\rho_1$
(which is straightforwardly reversible on appropriate 
components of pseudodiagrams):
\[
{\tilde\rho}^{-1}: \raisebox{-0.542in}{\includegraphics{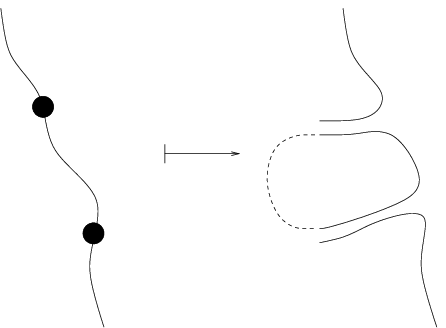}}
\]
and indeed for RR to appear the preimage would also have to contain
RR.
For LRL to appear would similarly require a preimage with a
(decorated) loop,
 so:
\\
2. the range of $\tilde\rho$ includes no diagram with LL or RR or LRL.
\\
Components with RLR can occur, as we have already demonstrated, but an
inverse image analysis similar to the above 
shows that a preimage precisely of the form
eliminated by the topological quotient 
(i.e. not present in the domain in the proposition) is required. 
 
Next we shall show that with the domain in the proposition no diagram
of the form eliminated by the topological quotient can appear as an
image. The component inverse image picture is of the form:
\[
{\tilde\rho}^{-1}: \raisebox{-0.542in}{\includegraphics{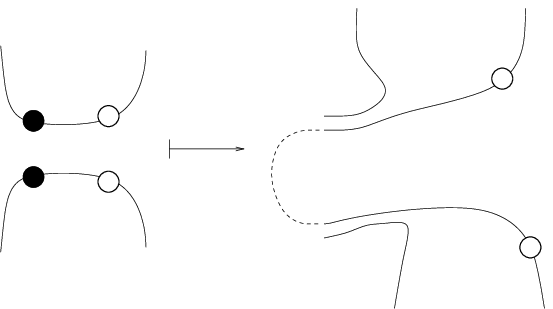}}
\]
The connection shown as a dashed line is forced since their can be no
lines decorated with L to the right of the line(s) decorated with R,
by the definition of LR-pseudodiagram.
The dashed line gives a component which does not occur in any diagram
in the domain, so we are done here.

NOW Think about showing surjectivity. (And 
construct the inverse?)

\parker{PROOF end NEEDED.}
}}

\medskip
\noindent {\em Remark:} 
Map $\tilde\rho$ extends in the obvious way 
to a map from $\pseud(m,m)^e$
to $\pseud(m-1,m-1)$ --- 
it is only necessary further to specify that 
every blob deformation 
(in the sense of the illustration to the definition of $\rho_1$)
that is not on a closed loop, 
except the lowest such, passes 
to the north of every blob deformation that is
--- see the figure below for an illustration. 

Map $\rho^*$ extends in the obvious way to the domain $\pseud(m-1,m-1)$,
provided again that closed loops are unambiguously treated 
in the deformation process --- let us say that they are deemed to 
congregate in the north-west.

The extended map $\tilde\rho$ is illustrated by the following example:
\[
\includegraphics{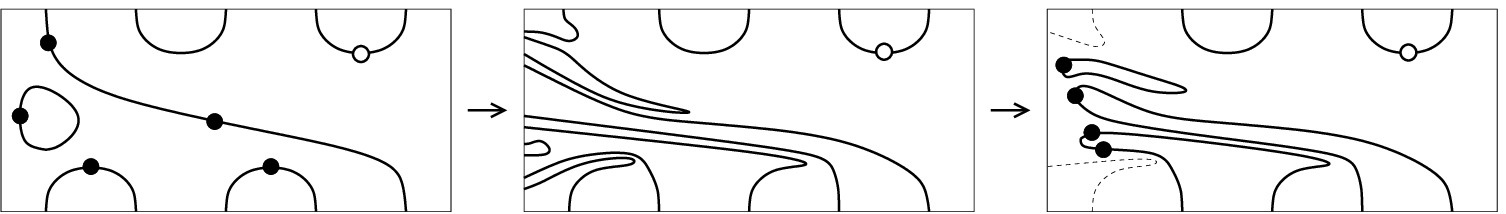}
\]
Note that we have
\[
\pseud(m,m)^e \stackrel{\tilde\rho}{\rightarrow} \pseud(m-1,m-1)
\stackrel{\rho^*}{\rightarrow} \pseud(m,m)^e
\]
but the two maps are not inverse with this domain and range, 
because of the necessarily
arbitrary choice of treatment of closed loops. 
With the choice specified, $\rho^*(\tilde\rho(d))$ is similar to $d$,
except that all `extra' blobs are gathered on the top left line,
regardless of where they where in $d$. 

\subsection{Algebra localisation and the \achiralb\ algebra}

The obvious candidate to construct our algebra isomorphism is
$\tilde\rho$ on $B^{x'e}_{m}$.  
But note that the range is not $B^{x'}_{m-1}$ here, 
and this is not a bijection. 
Consider the examples in (\ref{top q1}) 
(neglecting the $1/\delta_e$ factor for the moment):
Note that the right hand side of the bottom example 
in  (\ref{top q1})  is $k_L$
times the right hand side of the upper example (by the $RLR$ relation). 
Thus if we want to have an isomorphism we must require the same of the
left hand sides, that is, we must impose a further quotient relation. 
This means in particular 
that we eliminate the diagram with two double-decorated
lines from the basis. 


\del(top quot)
The `topological' quotient 
of $b^{x'}_{m}(\abparams)$ 
is defined by 
\eql(topquot)
\kappa_{LR} \;\; \raisebox{-0.2in}{\includegraphics{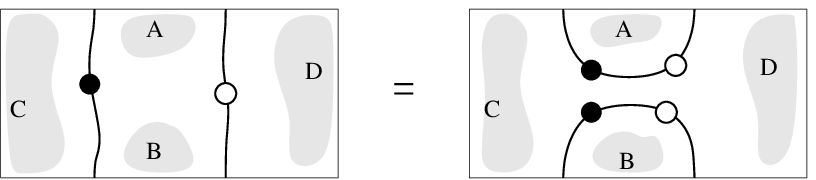}}
\eq
Here each labelled shaded area is shorthand for a certain subdiagram.
Thus each diagram restricts to the same subdiagram in 
the shaded region marked A (resp. B, C, D), 
but (\ref{topquot}) represents an identity for each such 
arrangement.
(Note that in C and D there must be a route for the adjacent blob to
the edge, hence no propagating lines.
Indeed there are no propagating lines at all on the right hand side)
\end{de}


\del(affine blob)
Define the \achiralb\ algebra $b^{x}_{m}(\abparams)$ 
to be the quotient of $b^{x'}_{m}(\abparams)$ by the 
additional set of relations (\ref{topquot}). 
\end{de}

\prl(a b basis)
The set $B^x_m$ is a basis for $b^{x}_{m}(\abparams)$. 
\end{pr}
{\em Proof:} 
We use Bergman's diamond lemma \cite{Bergman}. 

Let us define the {\it height} of a 
(pseudo)diagram to be the sum of the number of
loops and the number of decorations.  We observe that all the nontrivial
relations alter the height of a diagram.  This is easily checked except in
the case of the topological relation, in which case two of the decorations
are removed.  The vertical edges shown on the left hand diagram $D_0$ in
Definition~\ref{top quot} 
are the only vertical edges in $D_0$, and thus these edges can
only participate in a total of one loop in a product $D_1 D_0 D_2$
containing $D_0$.  This means that at most one loop can be deleted when
applying the topological relation to reduce the number of decorations, and
thus that the height decreases by at least $2-1 = 1$.

We may now choose a semigroup order, $<$, on the corresponding free
diagram algebra in such a way that smaller diagrams have strictly smaller
height. We aim to conclude that a minimal diagram in this sense is a basis
element. To do so, we need to worry about inclusion ambiguities (of which
there are none) and overlap ambiguities. (We are assuming familiarity
with Bergman's set-up.)  Imagine that the diagrams shown in 
Definition~\ref{top quot}   are
sandwiched between other diagrams to the top and bottom in a triple
product.

There is nothing to check for $k_L$-type relations, because they cannot
occur in ``even index" algebras, which is where this problem arises.

The $\kappa_L$ and $\kappa_R$ relations cause no problem because they commute
with the topological relation.

The $\delta_L$ and $\delta_R$ relations also commute with the topological
relation, because they never remove the last $L$ (respectively, $R$)
decoration.

The $\delta$ relation is easy to deal with because it cannot interact with
the topological relation, and thus the relations commute.

The only nontrivial case is the $\kappa_{LR}$ relation:

Suppose the top of the right hand side of 
Definition~\ref{top quot}   is part of a
$\kappa_{LR}$-type loop. 
Then we have a choice: we can contract the $\kappa_{LR}$
loop first and then apply the topological relation, or vice versa. The
ambiguity is resolvable here, however, because the region A (plus whatever
is just above it) can only be a disjoint collection of undecorated loops,
so they can be deleted and then recreated anywhere in the diagram where
they will not cause an intersection.

A similar case deals with the region B and a $\kappa_{LR}$ loop at the bottom.

If there is a $\kappa_{LR}$ feature on the left hand side, then there must be
two such features on the right hand side, and a similar argument again
applies.

According to Bergman's diamond lemma, we can conclude that the minimal
diagrams in this ``height" sense are a basis, as desired.
\Qed


\prl(chufster)
Provided that $\delta_L = \delta_e$ is invertible,
the map $\tilde\rho$ extends to an isomorphism
\[
e \; b^{x}_{m}  (\abparams) \; e \; \;  
\cong \;
b^{x}_{m-1}(\delta,\kappa_L,\delta_R,\delta_L,\kappa_R,\kappa_{LR})
\]
\end{pr}
{\em Proof:} We have to check 
$\tilde\rho (d_1 d_2) = \tilde\rho(d_1) \tilde\rho(d_2)$. 
Composition on both sides begins with pseudodiagram concatenation 
--- so it remains to check that pseudodiagram reduction is consistent.
This is a routine `diagram chase' similar to the blob case.
The difference is that there are $R$s present 
--- these largely play no role,
except that the $\kappa_{LR}$ reduction on the left becomes a $k_L$ 
reduction on the right, 
\[
\includegraphics{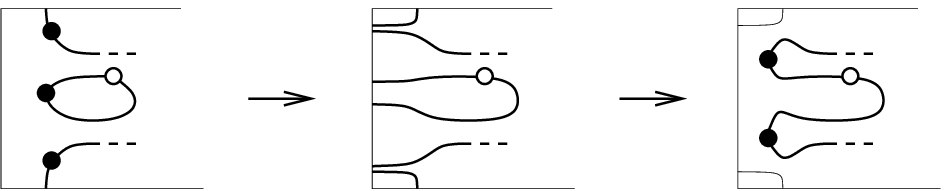}
\]
and vice versa 
(however note that these matters are already
resolved by our identification of these parameters in (\ref{KLR=KL})). 
\Qed

See (\ref{pc1}) and (\ref{pc2}) for diagrams exemplifying the
parameter change.
Also:
\[
\includegraphics{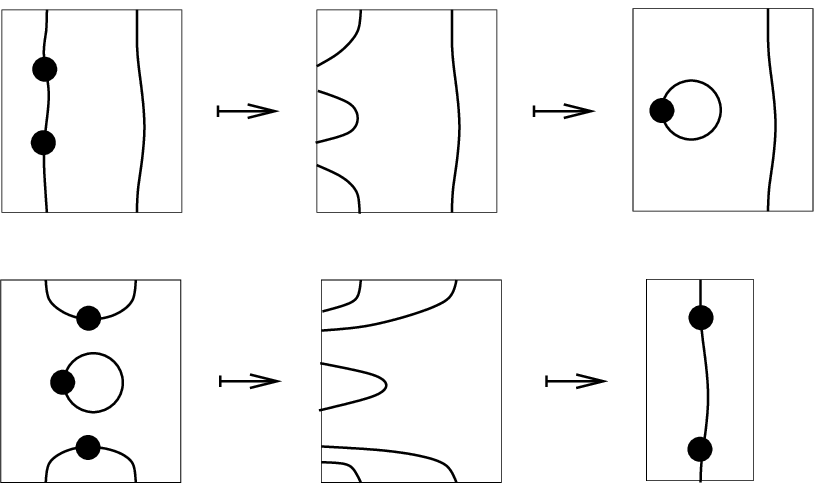}
\]
illustrating that $\delta_L$ and $\kappa_L$ are swapped.

Note that $b^{x}_{m}$ is a radical departure from the original blob algebra,
in that the topological quotient 
mixes between diagrams with different numbers of
`propagating lines'. This appears to deny us a powerful tool in
representation theory (cf. \cite{MartinSaleur94a}). 
However, 
we will determine the structure of this algebra by appealing to a
slightly different realisation.

{\em To reiterate:}
Just as the blob algebra is isomorphic to the 0-cover $m=2$ contour algebra, so
the idea of the east-west composite $m=2$ contour algebra 
(that is the variation in which lines which are 0-covered to the east or
the west may be decorated with a left (respectively right) blob)
also turns out to include a rather interesting case. 
There are a couple of ways in which such an algebra can be defined.
Firstly note that some lines, in some diagrams, 
are 0-covered both to east and to west.
Then left and right-blobs can meet on such a line. In general we will
consider them to be distinct, and even noncommuting on the line.
Consider the case in which we 
disallow multiple decorations with the property that it is not
possible to deform both the east leaning blob to the eastern edge and
the west leaning blob to the western edge simultaneously.
\section{Affine symmetric TL algebra} \label{ASTLA'}
One reason why this algebra $b_m^x$ is interesting is 
the existence of a doubled version of the {\em unfolding} map $\mu$. 
As with $\mu$ on $B_m$, 
the western blobs in an element of $B_m^x$ 
may again be used to map the diagram into left-right symmetric
versions (so far still with eastern blobs,  
now with mirror
images) about a reflection wall corresponding to the western edge. 
If we play the same game with eastern blobs we have another
reflection wall corresponding to the eastern edge, which is thus 
{\em affine} in the affine reflection group sense
\cite{Humphreys90}!
Altogether we have a fundamental domain (as it were) between these two
reflection walls, with a mirror image on each side (and then 
repeated reflections beyond).
Obviously then, the mirror image on the right is a translate of that
on the left and we have {\em periodicity}. 
Here is an example (the embossed letter `R' added to the fundamental
rectangle in this figure 
is only intended to emphasise the mirror images):
\[
\includegraphics{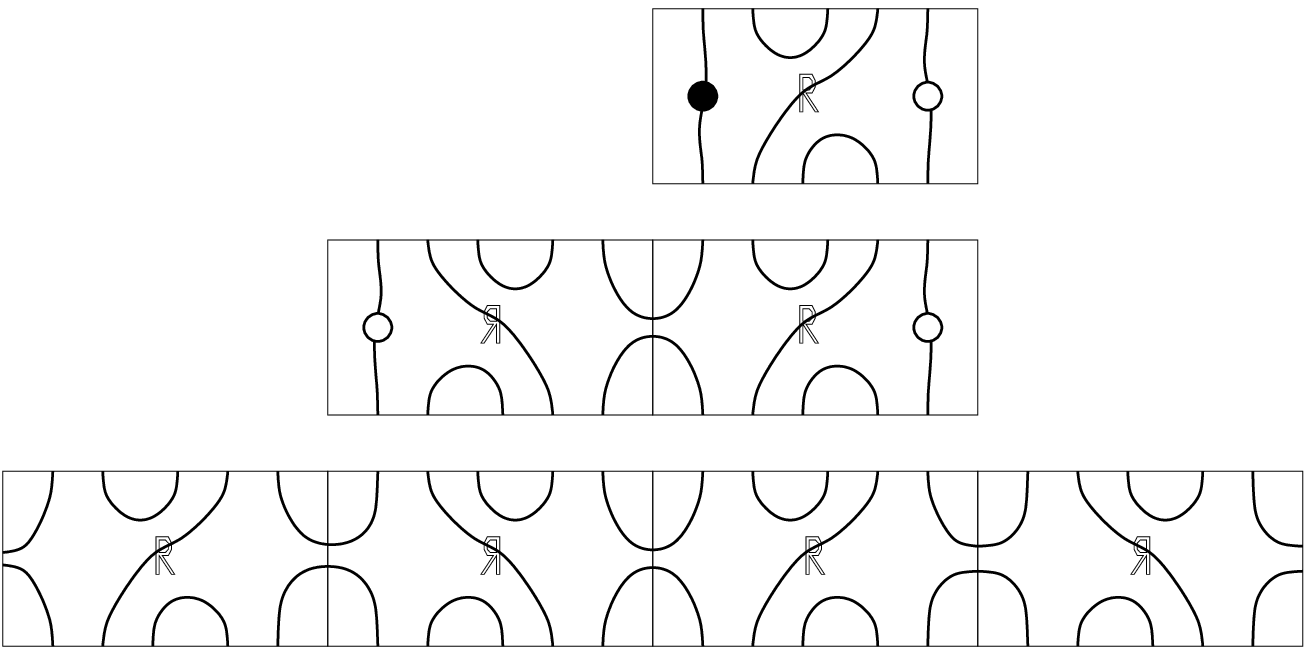}
\]
(Note as before that this is a well defined construct on isotopy classes.)
We observe that the resultant diagram is an element of 
$\Dppc{2m}$, 
the set of periodic TL diagrams with $2m$ vertices along each edge 
in the fundamental period (see \S\ref{deform}).
We denote by $\mu^x$ this unfolding map:
\[
\mu^x : B^x_m \longrightarrow \Dppc{2m}  . 
\] 
\subsection{On properties of the unfolding map $\mu^x$}
We can extend $\mu^x$ $\Ring$-linearly into a map from the algebra $b_m^x$
to $\Ring \Dppc{2m}$.  
It takes the generator $U_i$ in that algebra to a product 
$U_{i}U_{-i}$ (in a suitable labelling) in $\Dppc{2m}$, 
and so on. 
In other words we have a left-right symmetric 
subquotient-algebra of a {\em periodic} TL diagram algebra. 

Here is an example, of $\mu^x$ mapping to the cylinder realisation,
viewing
the cylinder along the axis, so that it appears as an annulus: 
\[
\includegraphics{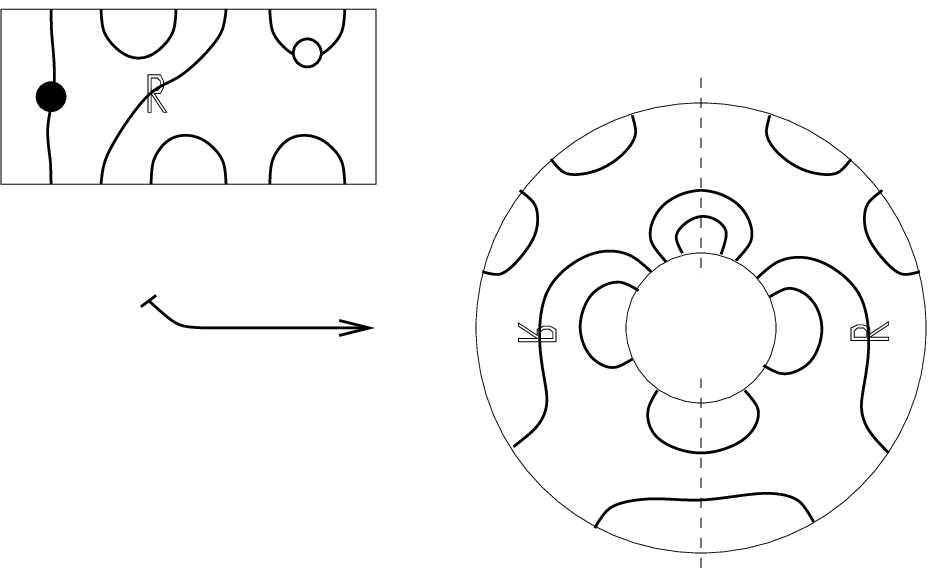}
\]
\label{Dph}
\del(sym period diagrams)
Let   
$\Dph_{2m}$  denote the set of left-right
symmetric periodic diagrams contained in $\Dppc{2m}$. 
(NB, noncontractible loops are still possible in $\Dph_{2m}$.)
\end{de}

There is a subalgebra of the periodic algebra spanned by $D^{\phi}_{2m}$. 
It will be evident from the illustration above
that $\mu^x(B_m^x)$ lies in this set. 

Note further
that these diagrams can be two-coloured like ordinary TL diagrams.
(Periodic diagrams with odd numbers of vertices cannot be
  two-coloured on the cylinder or annulus 
without a cohomology seam, but this need not concern us here.)  
For definiteness we fix that
the region touching the interval of the northern edge 
(which becomes the {\em inner} edge in the annular realisation)
astride the 0-reflection line is coloured {\em white}.
For example
\[
\includegraphics{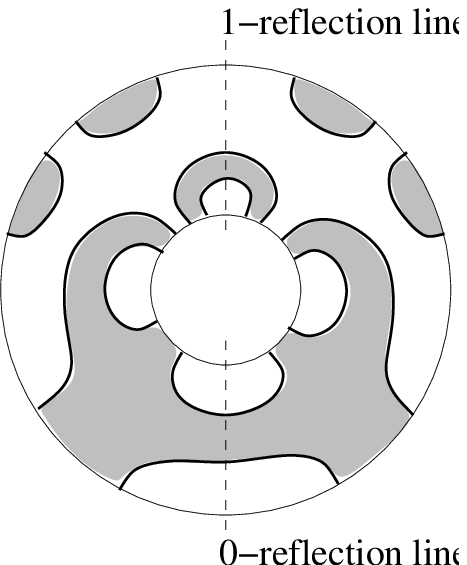}
\]
Note 
that there is a subset of these left-right symmetric periodic diagrams
with the property that the induced
colouring of the intervals of the southern edge coincides with
that on the northern edge. 
(Note that this is a proper subset in general, since it is
  possible to draw a symmetric periodic diagram in which precisely one
  line crosses the reflection line.)
Diagrams with this property are called {\em colouring composable} (CC)
  diagrams because, when they are concatenated in the usual way 
(top edge to bottom edge, or inner edge to outer edge in the annular
  realisation) the colouring we have specified
gives colours that agree across the join. 
It follows that the set 
of colouring composable
diagrams again spans a subalgebra.  
It will be evident that  $\mu^x(B_m^x)$ lies in this set
(since $\mu^x$ drags finite decorated segments of decorated lines 
out of the diagram, creating {\em pairs} of crossings of the 
reflection line). 


To see that the image $\mu^x(B_m^x)$ generates a quotient 
of the corresponding subalgebra of the 
periodic algebra, we should explicitly consider the image 
in the set $D^{\phi}_{2m}$ of 
symmetric periodic diagrams described above. 
\prl(goophy)
The map $\mu^x: B^x_{m} \longrightarrow D^{\phi}_{2m}$ is injective. 
\end{pr}
{\em Proof: }
The map is reversible at the point of 
deforming a blob out of the frame, so the issue is if isotopy on the
target side can equivalence two diagrams. 
However no contractible loops are   
produced from diagrams in 
$B^x_{m}$, so no such isotopy can arise.
\Qed

Note that this map is not surjective on arbitrary
CC symmetric periodic diagrams, 
since the maximum
number of noncontractible loops is 1 in the image. 
The quotients associated to the parameters $k_L$ and $\kappa_{LR}$ on
the blob side both have the effect of replacing a pair of
noncontractible loops with the factor  $\kappa_{LR}$. 

\subsection{Periodic pseudodiagram reduction}\label{reduc}

Recall 
from section~\ref{deform0}
that 
a periodic pseudodiagram is the generalisation of an ordinary
pseudodiagram from a rectangular to a cylindrical geometry 
(or unboundedly wide rectangles with finite periodic repetition). 

Here we will restrict to even period left-right symmetric colouring
composable pseudodiagrams. 
(Periodic left-right symmetric is the same as
affine symmetric, as already noted.)
Write $\CC_{2m}$ for the set of these pseudodiagrams with period $2m$.
We will colour such that the central northern interval is white.

As usual these diagrams are isotopy classes of concrete diagrams.
But as in section~\ref{deform2}, having taken a subset we have the option of
correspondingly strengthening the notion of isotopy. 
Here we consider isotopies that preserve the symmetry. 
(Note, then, that $\CC_{2m}$ is not the same thing
as the subset of general periodic pseudodiagrams with symmetric
representative elements, where a diagram with two contractible loops
on either side of the reflection line is isotopic to one with two
contractible loops both astride the reflection line.)

Composition on $\CC_{2m}$ is defined as before.


\prl(features)
The following list of features 
of concrete pseudodiagrams
are preserved by the 
$\CC_{2m}$
isotopy in this
setting,
and hence can be considered to appear (with well defined multiplicities)
in these pseudodiagrams: 
\newline
($\delta$) symmetric pair of loops (one each side of the symmetry line
--- the 0-reflection line);
\newline
($\delta_L'$) white loop astride the 0-reflection line;
\newline
($\kappa_L'$) black  loop astride the 0-reflection line;
\newline
($\delta_R'$)  white ($m$ even) (\resp\ black ($m$ odd)) 
loop astride the 1-reflection line;
\newline 
($\kappa_R'$)  black ($m$ even) (\resp\ white ($m$ odd)) 
loop astride the 1-reflection line;
\newline
($\kappa_{}'$) pair of noncontractible loops.
(Two such loops are `adjacent' if they may be deformed to touch, and the pair
is called black (\resp\ white) if one is on the black (\resp\ white)
side of the partition formed by the other.)
\Qed
\end{pr}


Let us write $B^{\phi}_{2m}$ for the subset of pseudodiagrams in
$\CC_{2m}$ with 
{\em none} of these features. 
For given period $n$ there are finitely many such pseudodiagrams. 

Define a map $$
\nu :
\CC_{2m} \rightarrow
D^o(V) |_{S=\{ b,w\}}
$$ 
(that is, $D^o(V)$ with two types of decoration) as follows.
Given a diagram in  $CC_{2m}$ consider the ur-diagram which is the
strip between the 0-reflection line on the left and the 1-reflection
line on the right (so the affine reflection group orbit of this strip is the
whole diagram). By the CC condition there are an even number of lines
leaving the strip through the  0-reflection line (and similalrly for
the  1-reflection line). Thus the  lines
leaving the strip through the  0-reflection line may be collected into
pairs,
such that the two lines in a pair are consecutive on the reflection
line.
This means that they can be brought arbitrarily close together
at the reflection line.
Joining each such pair with a blob (and similarly with a white-blob
at the  1-reflection line) we get an element of $D^o(V)$. 

Note that $\nu$ is not injective on $\CC_{2m}$. 
A diagram with two non-contractible loops and a (mirror) pair of
contractible loops above them is mapped to the same element of 
$D^o(V)$ as a diagram  with two non-contractible loops and a (mirror) pair of
contractible loops below them. 


\begin{lem}
The maps $\nu$ and $\mu^x$ induce a bijection between 
$B^x_m$ and $B^{\phi}_{2m}$.
\end{lem}
{\em Proof:} 
Firstly note that $\nu \circ \mu^x$ is the identity map on $B^x_m$
($\nu$ is the reverse of $\mu^x$, which is injective).

Secondly, consider $d \in \CC_{2m} \setminus B^{\phi}_{2m}$. 
It has at least one of the listed features.
It is routine to check
that each of these produces at least one contractible
loop in $\nu(d)$. 
Thus $\nu ( \mu^x(B^x_m) \setminus B^{\phi}_{2m})$ does not
intersect $B^x_m$. 
But by the previous paragraph 
 $\nu ( \mu^x(B^x_m) \setminus B^{\phi}_{2m})$
 is contained in 
 $B^x_m$, so it is empty.
Thus  $\mu^x(B^x_m) \subseteq B^{\phi}_{2m}$.

Next we show that $\nu(B^{\phi}_{2m}) \subset B^x_m$.
First consider $d \in \nu(B_{2m}^\phi)$ that does not have any
non-contractible loops. Any line in $d$ 
starting at the northern edge, say, and 
crossing the $0$ or $1$ line
cannot be propagating, and will have a 
corresponding 
line 
starting at the northern or southern edge 
paired to it by the CC'ness of $d$. Thus 
no string in  $\nu(d)$  has more than one blob on it,
and thus $\nu(d) \in B_m^x$.
If $d \in \nu(B_{2m}^\phi)$ does have a 
(necessarily unique)
non-contractible loop, then
this line under $\nu$ becomes 
part of 
the unique propagating line and is
decorated by exactly one black and one white blob. All other blobs come
from non-propagating lines 
combining in pairs 
as before.

It is easy to see
that $\mu^x \circ \nu |_{B^{\phi}_{2m}}$ is the identity map. 
Thus $\nu$ is injective when restricted to $B^{\phi}_{2m}$.
Thus finally the two sets have the same cardinality.
\Qed


Denote by 
$b_{2m}^{\phi}(\delta,\delta_{L}',\delta_{R}',\kappa_L',\kappa_R',\kappa_{}') $
the quotient of the  $\Ring$-algebra spanned by 
$\CC_{2m}$   
by the relations that each feature 
itemised in Proposition~\ref{features} 
 may be removed at the
cost of introducing a scalar factor as indicated 
 in Proposition~\ref{features} 
 in brackets (each such factor then appearing, note, 
as an argument to $b_{2m}^{\phi}$).

Since
all of the features of pseudodiagrams in Proposition~\ref{features} 
have multiplicity weakly increasing in composition 
we have
\prl(bbbbbasis)
The \aSTL\ algebra 
$b_{2m}^{\phi}(\delta,\delta_{L}',\delta_{R}',\kappa_L',\kappa_R',\kappa_{}')$
has basis $B^{\phi}_{2m}$. 
\end{pr}
{\em Proof:} 
One uses Bergman's diamond lemma much as in Proposition~\ref{a b basis}. 
\Qed


Again this is not the only way to produce a finite rank quotient.
For example $\kappa_{}'$ could be for the excision of black pairs
only
(see also \cite{MartinSaleur93}). 
However
\prl(x-hom)
The map $\mu^x$ extends to an algebra 
isomorphism 
\[
\mu^x: b^x_{m} \longrightarrow b^{\phi}_{2m}  
\]
with the obvious identification of parameters.
\end{pr}
{\em Proof:}
Note from the construction that $\nu$ commutes with diagram
composition, considered as a map from $\CC_{2m}$ to $D^o(V)$. 
It remains to show that the two different kinds of
pseudodiagram reduction yield the same
factors on each side. 
Applying $\nu$ to a diagram with two non-contractible loops will give
a diagram with a loop with both types of decoration on it
--- thus we set $\kappa'  = \kappa_{LR}$. \\
Applying $\nu$ to a diagram with a pair of contractible loops will give
a diagram with an undecorated loop 
--- reduction on either side gives a factor $\delta$. \\
Applying $\nu$ to a diagram with a white loop astride the 0-line will give
a diagram with a line with two left-blobs on it
--- thus we set $\delta_L'  = \delta_{L}$. \\
Applying $\nu$ to a diagram with a black loop astride the 0-line will give
a diagram with a loop with a left-blob on it
--- thus we set $\kappa_L'  = \kappa_{L}$. \\
The loops astride the 1-line  pass across similarly. 
\Qed

{\em Remark:}
Note that $b^{x'}_{m}$ does not map injectively into $b^{\phi}_{2m}$
without the `topological' quotient, since 
\be \non
\includegraphics{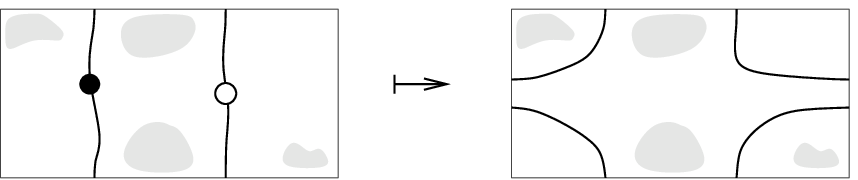}
\\
\includegraphics{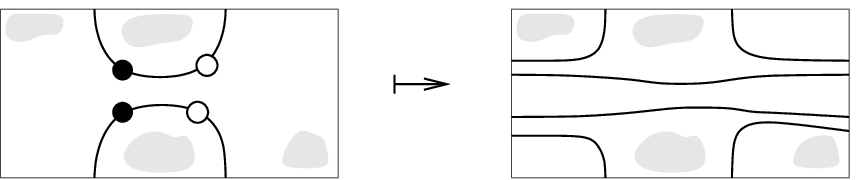}
\eeq

We will see that 
even before the quotient the affine symmetric subalgebra
is not complicated by as many
`infinities' as the ordinary periodic TLA.
We will also see that it is amenable to the same recollement treatment
as the non-affine case above. 

The claim is that the set of diagrams that contain a cup and cap
astride the 0-reflection line is a subset that spans an (idempotent)
subalgebra (a similar statement holds for the 1-reflection line). 
There is a bijective map into the set of all diagrams with one fewer vertex on
each side obtained by simply removing this cup and cap. 
In order to elevate this to the status of an algebra homomorphism we
will again have to take care with the parameters. 

\section{Representation theory of ASTLA} \label{rep ASTLA}
In what follows $f$ corresponds to the right blob, in the way that $e$
(or $\abe'$) 
corresponds to the left blob:
\[
f \mapsto \;\;\;
\raisebox{-2.cm}{
\includegraphics{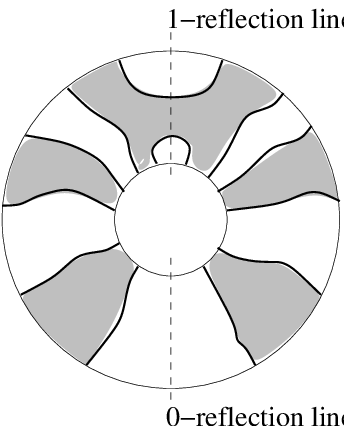}
}
\]
\subsection{General and generic results}
Assume that  $\delta_L$ is invertible in $\Ring$.
Let 
$\rho': D^{\phi}_{2m-2} 
 \longrightarrow \frac{1}{\delta_L} D^{\phi}_{2m}
 \; \subset K  D^{\phi}_{2m}$
denote the map that inserts a cup and cap astride the 
0-reflection line and then rescales by $ \frac{1}{\delta_L}  $.

Note that the map $\rho'$ is injective, with image 
the set of {\em all} (rescaled) diagrams in $D^{\phi}_{2m}$ 
with a cup and cap astride the 
0-reflection line. 
Thus
$\rho'( B^{\phi}_{2m-2})$ spans the
subalgebra 
$\frac{\abe'}{\delta_L} b_{2m}^{\phi}(\delta,...) \frac{\abe'}{\delta_L}$
of $b_{2m}^{\phi}$ 
in a similar way to the ordinary blob case.

Let $\rho$ denote the map on 
$ \rho'( D^{\phi}_{2m-2})  $  
that removes this cup and cap and normalisation (inverse to $\rho'$). 

\prl(rho hom)
The map $\rho$ extends to  an algebra isomorphism 
\eql(ebe)
\abe' \; 
b_{2m}^{\phi}(\delta,\delta_{L},\delta_{R},\kappa_L,\kappa_R,\kappa_{LR}) 
\; \abe' \;\; 
\stackrel{\rho}{\longrightarrow}  \; 
b_{2m-2}^{\phi}(
\delta,\kappa_L,\delta_{R},\delta_{L},\kappa_R,\kappa_{LR} ) . 
\eq
\end{pr}
{\em Proof:}
This follows from Propositions~\ref{chufster} and~\ref{x-hom}.
However it is useful to sketch a direct proof analogous to
Proposition~\ref{recool}. 

In order to readily distinguish $\delta_L$ and $\kappa_L$ in the
periodic realisation it is again useful to two-colour the diagrams. As
before we set the interior of the region whose closure includes the
northern interval astride the 0-reflection line to {\em white}. 
Then $\delta_L$ loops of 
$b_{2m}^{\phi}(\delta,\delta_{L},\delta_{R},\kappa_L,\kappa_R,\kappa_{LR})$ 
are white loops astride this line. 
The colour of the corresponding interval (and hence loops) astride the
1-reflection line depends on whether $m$ is odd or even. 
(Thus the image of $f.f=\delta_R f$ is a black loop if $m$ is odd.)

As in the ordinary blob case, comparing $\rho(a)\rho(b)$ to $\rho(ab)$
the underlying diagrams agree, and there is a correspondence between
the loops produced on each side, but the cup and cap removal means
that they all change colour. 
Thus, on applying $\rho$ to a pseudodiagram
(a diagram with loops which reduce to scalars), 
the roles of $\delta_L$ and $\kappa_L$ are interchanged. 
As for $\delta_R$ and $\kappa_R$, there is a colour change due to
$\rho$, which nominally interchanges them, but $\rho$ also changes $m$
between odd and even, changing back. 
Altogether, we have $\rho(a)\rho(b)=\rho(ab)$ with $\rho$ as in (\ref{ebe}). 
\Qed

We could bring the two sides of (\ref{ebe}) closer together by making
further constrained parameter choices. 
However, here we will concentrate on the generic case, i.e. 
{\em working with field $k=\C$, so that parameter space
can be endowed with the Zariski topology, } 
we assume that Zariski open subsets of points in parameter space 
all have basically the same representation theory
(in the sense that the basis for simple module $M_{\lambda}$ (say) 
is in each case the image of the {\em same} module basis over ground ring $K$
under $M^K_{\lambda} \mapsto k \otimes_K M^K_{\lambda}$ 
by specialisation).
We will verify this assumption shortly. 
Under this assumption we may consider a single meta-category 
$b_{2m}^{\phi}$-mod of left modules for $b_{2m}^{\phi}$.
Then by (\ref{glob}) 
\prl(meta-cat)
Map (\ref{ebe}) provides a full embedding $G$ of $b_{2m-2}^{\phi}$-mod in 
$b_{2m}^{\phi}$-mod.
\Qed
\end{pr}
This means in particular that we can construct prestandard modules by
a recursive procedure.
(NB, much representation theory 
can still be done with the restrictions on the field imposed here removed,
but at considerable cost in brevity.)


\begin{theo}\label{index wangy}
An index set for equivalence classes of simple $b_{2m}^{\phi}$-modules
for generic parameters is $\{ 0 \}$ for $m=0$; and 
$\Lambda^{\phi}_{m} = \{ -m,-m+1,..,0,..,m-1 \}$ for $m>0$. 
\end{theo}
{\em Proof:}
By Proposition~\ref{grebo}  
there is for each simple module (irreducible representation) 
$S_{\lambda}(2m-2)$
in  $b_{2m-2}^{\phi}$-mod a simple module in $b_{2m}^{\phi}$-mod
which is the head of prestandard module $\glob(S_{\lambda})$
(using the $\glob$ corresponding to (\ref{ebe})). 

The simple modules not constructed in this way
 are those $S$ obeying $eS=0$. Since $U_1eU_1 \sim U_1$ 
the element $U_1$ acts as zero in such representations,
and hence so do all the $U_i$s. 
Thus only 1 and (possibly) $f$ are (represented by) non-zero. 
Since $f$ is (pre) idempotent there are precisely two such simple
modules in general, both one-dimensional, one where $f$ is zero 
(which module we will give the label $\lambda=-m$) and
the other not (which module we will give the label $\lambda=m-1$). 
(NB
the `bootstrap' case at $m=1$ 
is the exception, since there, uniquely,
$fef\sim f$ and $f \cong 0$ is also forced.)
\Qed

\subsection{Combinatorics of the basis $B^{\phi}_{2m}$}

\newcommand{\BD}{B} 

Note that if a diagram in $B^{\phi}_{2m}$
has any propagating lines then it has at least
two (by the symmetry), and that there is a unique mirror image pair
that can be deformed to touch the 0-reflection line at some point
--- the pair `closest' to the 0-reflection line. 
There is thus a unique region of the diagram that touches both
elements of this pair {\em and} contains a segment of the 0-reflection line. 
This region can be black or white. We will call it the {\em inner}
region
(and call the pair of lines the inner lines).

One way of organising diagrams into subsets is by the number of
propagating lines. Within this, we may subdivide the set of those with
$2l>0$ propagating lines into those in which the 
inner region is black or white. 
Let us denote these subsets $\BD^{\phi}_{2m}[\pm l]$ respectively (note
that this works at $l=+0=-0$ since there is no inner region there). 

For example:
\eql(bwex)
 \BD^{\phi}_{2m}[-m] = \{ 1 \} ; \qquad 
 \BD^{\phi}_{2m}[-(m-1)] = \{ f \}  ; \qquad 
 \BD^{\phi}_{2m}[m-1] = \{ e \} . 
\eq
Remark: It is necessary to have such a labelling scheme for these
subsets, and this scheme will serve our purposes. However, it is not
canonical. 


Following \cite{MartinSaleur94a} let us write $\hash(d)$ for the
number of propagating lines in diagram $d$. 
Similarly we extend this to apply to any scalar multiple of $d$,
so that $\hash(dd')$ is defined for any two diagrams $d,d'$.
We write $c_i(d) \in \{ b,w \}$ for the inner region colour of $d$
(and again similarly for $dd'$). 

\begin{lem} \label{filtre}
For all $d,d'$ \\
(1)  We have $\hash(dd') \leq \hash(d)$.
\\
(2) If $\hash(dd') = \hash(d)$ then $c_i(dd')=c_i(d)$. 
\end{lem}
{\em Proof:}
(1) is straightforward.
(2) Suppose that a certain pair of lines are inner in $d$. 
These lines are still identifiable beginning at the 
northern edge  of $dd'$ (which is inherited from $d$). 
If they
remain inner in the extension of $d$ by $d'$ then obviously the colour is
unchanged. On the other hand if they are not inner in $dd'$
then they are no longer propagating and $\hash(dd') < \hash(d)$.
\Qed

For $l \geq 0$ define 
\newcommand{\BBD}{BB^{\phi}_{2m}}
\[
\BBD(l) \; = \; 
\bigcup_{0 \leq j\leq l} 
\left(  \BD^{\phi}_{2m}[j] \cup   \BD^{\phi}_{2m}[-j] \right)
\]
and $\BD^{\phi}_{2m}(0)= \BD^{\phi}_{2m}[0]$ 
and for  $l \geq 1$
\[
\BD^{\phi}_{2m}(\pm l) \; = \;  
\BD^{\phi}_{2m}[ \pm l] \cup \BBD(l-1) . 
\]


\prl(ideallypoo)
For $l \in \{ -m,-m+1,...,m-1 \}$ the set $\BD^{\phi}_{2m}(l)$
is a basis for an ideal of $b^{\phi}_{2m}$, and the subset structure
\[
\UseComputerModernTips
\xymatrix{
          &  \BD^{\phi}_{2m}(m-1) \ar @{_(->}[dl] <1ex>  
                  &  \BD^{\phi}_{2m}(m-2) \ar@{_{(}->}[l]
                                        \ar@{_(->}[ddl] <1ex>  
                     &&
\\
\BD^{\phi}_{2m}(-m) & & & 
          \cdots & &  \BD^{\phi}_{2m}(0) \ar@{_{(}->}[ul]
                                       \ar@{_{(}->}[dl] <1ex>  
\\
          &  \BD^{\phi}_{2m}(-(m-1)) \ar@{_{(}->}[ul] 
                  &  \BD^{\phi}_{2m}(-(m-2)) \ar@{_{(}->}[l] 
                                           \ar@{_{(}->}[uul] 
                     &&
}
\]
passes to a subideal structure (over any ring).
\end{pr}
{\em Proof:} This follows directly from Lemma~\ref{filtre}. \Qed


Let us name the subideals
$
I^{\phi}_{2m}(l) \; = \; k \BD^{\phi}_{2m}(l) . 
$
Noting this structure, 
associate a partial order $\lhd$ to $\Lambda^{\phi}_{m}$ by 
 $l \;\lhd\; l'$ if and only if $|l| < |l'|$
(NB, this order is not total).
We then define $I^{\phi}_{2m}[l]$ as  
the $l$-th section of this ideal structure, that is 
\[
I^{\phi}_{2m}[ l] \; 
  = \left\{ \begin{array}{ll} 
\; k \BD^{\phi}_{2m}(l) / k \BBD(-l-1) & l<0 \cr
\; k \BD^{\phi}_{2m}(l)               & l=0 \cr
\; k \BD^{\phi}_{2m}(l) / k \BBD(l-1) & l>0 . \cr
\end{array} \right.  
\]
Note that $I^{\phi}_{2m}[l]$ has basis $\BD^{\phi}_{2m}[l]$
(where the action is the algebra multiplication, but taking account of
the quotient).


Next we want to decompose these sections as left-modules,
and equip their component modules with an inner product.

\pdef{Half-diagrams:}
Note that it is always possible to cut a diagram 
$d \in B^{\phi}_{2m}$
from the eastern to
the western edge in such a way that only propagating lines are cut
(and these exactly once each). 
Note, however, that this process is not always unique (even up to isotopy),
since some diagrams 
with { no} propagating lines 
contain a noncontractible loop, which could lie
above or below the cut. 
However, considering the set of {\em half-diagrams} produced in this
way
{ ignoring} any noncontractible line, 
any top-bottom pair of half-diagrams 
with the same number of propagating lines, and inner region colour if defined,
can always be recombined to produce a full CC diagram, 
and in exactly one way, with the caveat
that the CC requirement will determine if a noncontractible loop must
be inserted. 
We will call any such loop a {\em belt}. 

\newcommand{\ket}[1]{| #1 \rangle}%
\newcommand{\bra}[1]{\langle #1 |}%

Let us denote by $\ket{\BD^{\phi}_{2m}[l]}$ the set of upper half
diagrams associated to $\BD^{\phi}_{2m}[l]$ (any $l$),
and by $\ket{d}$ (the ``ket'')
 the upper half diagram obtained from diagram $d$
(write $\bra{d}$  (the ``bra'') for the corresponding lower half diagram). 
There is an obvious isomorphism of the set  
$\ket{\BD^{\phi}_{2m}[l]}$  with the set of bottom
halves $\bra{\BD^{\phi}_{2m}[l]}$, obtained by reflecting in an east-west
line, 
that is, $\ket{a} \stackrel{\sim}{\mapsto} \bra{a^o}$. 
It follows that 
\begin{eqnarray}
\BD^{\phi}_{2m}[l] & \cong & 
    \ket{\BD^{\phi}_{2m}[l]} \times \bra{\BD^{\phi}_{2m}[l]}
\\
d  & \mapsto &  ( \; \ket{d} \; , \;\; \bra{d} \;)
\end{eqnarray}
where the map is the cut map. 
(Note that elements of $\ket{\BD^{\phi}_{2m}[+l]}$ and  
$\bra{\BD^{\phi}_{2m}[-l]}$
can be concatenated with $l>0$, but they will not produce a CC
diagram.)
We may write the inverse map as a multiplication:
\[
 ( \; \ket{d} \; , \;\; \bra{d} \;)  \mapsto  \ket{d} \bra{d} . 
\]


\newcommand{\Belt}{{\mathrm{Belt}}}%
For $\ket{a} \in \ket{\BD^{\phi}_{2m}[0]}$ write $\Belt_a$ for the
subset of elements $\bra{b} \in \bra{\BD^{\phi}_{2m}[0]}$ 
such that $\ket{a}\bra{b}$ has a belt. 
The partition of $\bra{\BD^{\phi}_{2m}[0]}$ into two parts defined by
$\Belt_a$ (as one of the parts) is independent of $a$, and written
simply as $\Belt$. 
Note that 
\begin{lem} \label{(x)-}
If  $d,d' \in \BD^{\phi}_{2m}[l]$ and 
$\hash(dd') = |l|$ 
then 
\\
(i) there is a  monomial in the parameters
 $k_{dd'}$ such that 
\eql(dd')
dd'  
= k_{dd'} \ket{d} \bra{d'} . 
\eq
\\
(ii) This $k_{dd'}$ depends on $ \bra{d}$ and $ \ket{d'}$ but 
does not depend on 
 $ \bra{d'}$ and $ \ket{d}$, except in case $l=0$ through the 
noncontractible loop caveat.
\\
(iii) No top-bottom symmetric diagram $\ket{a} \bra{a^o}$ has a belt.
If $l=0$ and $dd' =  \ket{a} \bra{b} \ket{c} \bra{a^o}$ 
write $ \bra{b} \ket{c}_a$ for $k_{dd'}$. 
Let $M_a$ be the matrix $(\bra{b} \ket{c}_a)_{b,c}$. 
Let $M_a'$ be
the matrix obtained from $M_a$ by dividing every row with $b \in
\Belt_a$ by $\kappa_{LR}$;
and let $M_a''$ be
the matrix obtained from $M_a$ by dividing every column with $c \in
\Belt_a$ by $\kappa_{LR}$. 
If $a$, $a'$ are in the same part of $\Belt$, then $M_a=M_{a'}$.
If $a$, $a'$ are in different parts of $\Belt$, then $M_a'=M_{a'}''$.
\\
(iv) Now let $d' \in B_{2m}^\phi[l]$ and $d'' \in B_{2m}^\phi[l']$ be such
that $\#(d'' d') = |l|$
($d''$ could have more than $|l|$ propagating lines).  
Then 
$$
d'' d' = \sum_{d \in |B_{2m}^\phi[l]\rangle} k'_d |d \rangle \langle d'|
,$$ 
where the $k'_d$ depend on $d''$ and $\langle d'|$, but not on
$|d' \rangle$.
\end{lem}
{\em Proof:} (i) The product
$dd'$ is some scalar times a diagram. Under the given conditions it is
clear that this diagram must have the given ket-bra form.
\\
(ii) This follows from the definition of the cut map.
\\
(iii) The first part is obvious. If  $a,a'$ are in the same part of
$\Belt$, then the calculations for each  $ \bra{b} \ket{c}_a$ and 
 $ \bra{b} \ket{c}_{a'}$ are identical. 
Otherwise there are four types of case in comparing  $ \bra{b} \ket{c}_a$ and 
 $ \bra{b} \ket{c}_{a'}$:
\\ 
1) if $\{b,c\} \cap \Belt_a = \{b,c\}$ then the matrix elements differ
precisely by a factor  $\kappa_{LR}$ on the $a$ side, which is adjusted by the
division on the $a$ side; 
\\
2) if $\{b,c\} \cap \Belt_a = \{b\}$ then both sides have at least one
factor   $\kappa_{LR}$, and are the same, and both have this factor
divided out;
\\
3) if $\{b,c\} \cap \Belt_a = \{c\}$ similarly;
\\
4) if $\{b,c\} \cap \Belt_a = \emptyset$  then the matrix elements differ
precisely by a factor  $\kappa_{LR}$ on the $a'$ side, which is adjusted by the
division on the $a'$ side. 
\\
(iv) As for (i). 
(NB, 
 there can only be one term in the sum, because the diagram
basis is a monomial basis, i.e. the product of any two diagrams
is a scalar multiple of another.)
\Qed
\\
Indeed, since 
\eql(dd'xx)
dd'  =  \ket{d}\bra{d} \; \ket{d'}\bra{d'}
= k_{dd'} \ket{d} \bra{d'}
\eq
we will sometimes write $\bra{d}  \ket{d'}$ for $k_{dd'} $ 
when no ambiguity arises.


\prl(struc1)
Let $d,d' \in \BD^{\phi}_{2m}[l]$. 
\\
(i)
There exist diagrams $a,b$, and a 
nonzero monomial in the 
parameters $k$, such that 
$adb=kd'$. 
That is, provided all the parameters are units then every diagram in  
$\BD^{\phi}_{2m}[l]$
generates
$ \BD^{\phi}_{2m}(l)$.
\\
(ii) For $a$ any diagram, if $\hash(dad')=l$ then  
$dad'=k_a \ket{d}\bra{d'}$ where $k_a$ is a nonzero monomial in the  
parameters.
\end{pr}
{\em Proof:} (i) Note that if $d^o$ is the `opposite' diagram of $d$ 
(the same diagram
drawn upside-down) then $\hash(d^o d)=\hash(d d^o)=l$. 
Consider $a = \ket{d'}\bra{d^o}$ and $b=\ket{d^o}\bra{d'}$, then
\[
adb = \ket{d'}\bra{d^o} \ket{d}\bra{d}  \ket{d^o}\bra{d'} 
= k \ket{d'}\bra{d'} =kd'
\]
For the second part of (i) it is now enough to show that we can get from some
diagram in $\BD^{\phi}_{2m}[l]$ to (some appropriate scalar multiple of)
some diagram in $\BD^{\phi}_{2m}[\pm(l-1)]$. 
This is routine for a suitable choice of diagram in each case.
For example, for $2j \leq m-2$ then $w=e U_2 U_4 \ldots U_{2j} f 
\in \BD^{\phi}_{2m}[m-2j-2]$; while for  $2j < m-2$ then 
$U_{m-1} w U_{m-1} = \kappa_R w'$ with 
$w' = e U_2 U_4 \ldots U_{2j} U_{m-1} \in \BD^{\phi}_{2m}[m-2j-3]$.
\\
(ii) Similarly we have:
\[
dad' =  \ket{d}\bra{d} \ket{a}\bra{a}  \ket{d'}\bra{d'}  
=  ( \bra{d} \ket{a} \; \bra{a}  \ket{d'} ) \; \ket{d}\bra{d'} 
\]
and  $ \bra{d} \ket{a}$ and   $\bra{a}  \ket{d'}$ are
nonzero by construction. 
\Qed

An immediate corollary to \ref{struc1}(i) is 
{\co{ \label{imco}
Subject to the same parameter restriction as in \ref{struc1}(i), no
unit multiple of any diagram is in the radical.
\Qed
}}


\newcommand{\rad}{\mbox{rad}}

For $i=0,1,2,...,m-1$ define $S_i = I_{2m}(-i)$ and 
$T_i =  I_{2m}(i) +  I_{2m}(-i)$.
Note by Proposition~\ref{ideallypoo} 
that the following is a chain of ideals in  $b_{2m}^\phi$
\eql(hch)
 S_0  
\subseteq 
 S_1  
\subseteq 
 T_1  
\subseteq
 S_2  
\subseteq 
 T_2  
\subseteq 
\cdots \subseteq 
 T_{m-1}  
\subseteq 
 S_{m} = I_{2m}(-m) = b_{2m}^\phi . 
\eq


For each $l \in \Lambda^{\phi}_m$ and  $d \in \BD_{2m}^\phi[l]$
define $b_{2m}^\phi$  modules by
$\mathcal{S}^d_{2m}(0) := b_{2m}^\phi d $
if $l=0$ and otherwise 
\eql(standby)
\mathcal{S}^d_{2m}(l) := \frac{b_{2m}^\phi d \; + \; S}{S}
\eq
where 
$S=T_{|l|-1} $.

Note that right multiplication by $d'$ gives a map 
$\gamma_{d'}: \mathcal{S}^d_{2m}(l) \mapsto \mathcal{S}^{d'}_{2m}(l)$.
Thus by Proposition~\ref{struc1}(i):
\begin{lem}
If $\delta, \delta_L, \delta_R, \kappa_L, \kappa_R, \kappa_{LR} $ 
are units
then 
the precise choice of $d$ is irrelevant, 
up to isomorphism, in $\mathcal{S}^d_{2m}(l)$.
In this case define
$ \mathcal{S}^{}_{2m}(l) = \mathcal{S}^{d}_{2m}(l)$.
\end{lem}

Note that  $d \in \BD_{2m}^\phi[l]$ can always be chosen to have no
non-contractible lines. In this case 
the module 
$ \mathcal{S}_{2m}(l)  $
has basis 
$b_{2m}^\phi d \cap \BD_{2m}^{\phi}[l] $,  
where the action is the algebra multiplication, but taking account
of the quotient. Further 
$b_{2m}^\phi d \cap \BD_{2m}^{\phi}[l] \cong \ket{\BD^{\phi}_{2m}[l]}$. 


We may extend the notation $ \bra{d} \ket{d'}$ 
to a map 
$( \bra{ \BD^{\phi}_{2m}[l] }  , \ket{ \BD^{\phi}_{2m}[l] }  ) \rightarrow K$ 
by  
$ \bra{d} \ket{d'}=0$ if $\#(dd')<l$. 
Via the involution $\ket{} \stackrel{\sim}{\rightarrow} \bra{}$ the map 
$\bra{-} \ket{-}$ extends (bi)linearly to an inner product on  
$\mathcal{S}^d_{2m}(l)$.


\newcommand{\brac}{\langle}
\newcommand{\cket}{\rangle}

Define
\eql(innuit)
\Deltag^{d''}_{2m}(l) = \mathrm{det}\left( (\langle d | | d'
\rangle)_{n\times n}\right)
\eq
where $d^o,d' \in b^{\phi}_{2m}  d ''\cap B^{\phi}_{2m}[l]$, 
for a fixed $d'' \in
B^{\phi}_{2m}$ (the basis of $\sS_{2m}(l)$)
and $n = \dim \sS_{2m}(l)$.
Note 
from Lemma~\ref{(x)-}
that $\Deltag^{d''}_{2m}(l)$ does not depend on $d''$ if 
$l \neq 0$;
and depends on $d''$ only through at most an overall factor of 
a power of $\kappa_{LR}$
if $l=0$. 
Choosing the lowest power in this overall factor, we will write
simply $\Deltag^{}_{2m}(l)$ in all cases. 


\newcommand{\lambdal}{l}

We will address the specific computation of the Gram
determinant later, but note that the 
matrix entries
are monomial in the parameters. 
We have
\prl(gen semisimp)
For each $\lambdal \in \Lambda^{\phi}_{m}$ there is a
 polynomial $P$ in the parameters such that 
the prestandard module of $b_{2m}^{\phi}$  associated to $\lambdal$ 
has Gram determinant given by evaluation of $P$ at the
appropriate specialisation. 
Every prestandard is generically simple.
\end{pr}
{\em Proof:} 
By computing $\Deltag^{}_{2m}(\lambdal)$ 
with the parameters treated as
indeterminates we obtain the polynomial $P$. 
By  Proposition~\ref{orestes} the inner product we have
defined is unique up to scalars, 
and the Gram determinant is nonvanishing 
in any given parameter choice if and only if the module is
simple there. 
`Generically' has the meaning of Zariski open here, so
it is only necessary to show that no such polynomial $P$ is
identically the zero polynomial. 
This can be done by considering asymptotic cases
of the parameters. 
(The power of $\delta$, say, is maximal on the
 diagonal for all rows of the matrix, and uniquely maximal there for
 some rows, so the determinant is not zero.)
\Qed


\begin{theo}
Let $*$ be the $\Ring$-linear involutory antiautomorphism 
on $b^{\phi}_{2m}$  defined by flipping each diagram upside-down,
i.e. by reflection in a horizontal line. 
For each $l \in  \Lambda^{\phi}_{m}$ let $M(l)$ be one of the
diagram bases of $ \mathcal{S}^{}_{2m}(l)  $;
and let $C$ be the ket-bra combination of basis elements.
Then 
$b^{\phi}_{2m}$ is cellular over $\Ring$ with cell datum 
$( \Lambda^{\phi}_{m},M,C,*)$.
\end{theo}
{\em Proof:}
Note that
$*$ is an algebra antiautomorphism
by top-bottom symmetry of the reduction rules.
Proposition~\ref{ideallypoo}, Lemma~\ref{(x)-} and 
(\ref{standby}) verify the axioms given
in \cite{GrahamLehrer96}.
\Qed


\begin{theo}
If  
$\delta, \delta_L, \delta_R, \kappa_L, \kappa_R, \kappa_{LR} $ 
are units
then
$b_{2m}^{\phi}$ is quasihereditary,
and the chain (\ref{hch}) is a heredity chain. 
\end{theo}
{\em Proof:}
It is enough to show that the chain is heredity
(one might also see \cite{KoenigXi99b}). 
We need to show for $S= S_i$ or $S= T_i$
where
$A=b_{2m}^\phi/T_{i-1}$ if $S=S_i$ and 
$A=b_{2m}^\phi/S_{i}$ if $S=T_i$:
\\ \noindent (1) that $S^2 = S$
\\ \noindent (2) that $S J S = 0$ for $J = \rad(A)$ 
\\ \noindent (3) that each section defined by the proposed hereditary chain
of ideals is projective in the quotient. 
I.e. that $S/T_{i-1}$ (resp. $S/S_i$ is projective in $A$).


Suppose $S=S_i$ for some $i$,
so $S_i$ contains all diagrams with $i$ propagating lines and white
inner region and all diagrams with fewer than $i$ propagating lines.

Take $d \in \BD_{2m}[-i]$. Then $\# (d d^{\circ}) =i$.  
Thus $S_i^2$
contains a diagram with white inner region and $i$ propagating
lines. Proposition~\ref{struc1}(i) then says that $S_i \subset S_i^2$ and so
$S_i^2 = S_i$.

We now show that $SJS =0$ where $J=\rad b^{\phi}_{2m}/T_{i-1} $. 
Consider $djd'$, with $d$, $d' \in
\BD^{\phi}_{2m}(-i)$ and  $j \in J$.
Now write $j = ( \sum m_{\alpha} j_{\alpha}) +T_{i-1}$ 
where $m_\alpha$ are in the
base ring $K$ and $j_{\alpha} \in B_{2m}^{\phi}$. 
If $\# (d j_{\alpha} d')\le i-1$ then 
$ d j_{\alpha} d' $
 is zero in the quotient
$b_{2m}^{\phi} /T_{i-1}$.

If  $\# (d j_{\alpha} d')= i$ then 
Proposition~\ref{struc1}(ii) says that 
$d j_{\alpha} d' = k_{\alpha} dd'$ for some
$k_{\alpha}$ in the ring and so
$d j d' = \left(\sum m_{\alpha} k_{\alpha}\right) dd' +T_{i-1}$
(where the sum is over those $j_{\alpha}$ such that 
$\# dj_{\alpha}d' = i$). 

If $\# (d d') <i$  then $djd'=0$ in the quotient. 
If $\# (d d') =i$  then $(dd')^r = k^r dd' \not \in T_{i-1}$ 
for some $k$ a non-zero
monomial in the parameters and for all $r$ and so $dd' + T_{i-1}
\not\in J$. 
Thus for $djd'$ to be nilpotent in the quotient we need 
$\left(\sum m_{\alpha} k_{\alpha}\right) dd'$ to be zero and thus
$d j d' =0$. 
Thus $SJS=0$.

Finally we need $S_i/T_{i-1}$ to be projective as a
$b_{2m}^{\phi}/T_{i-1}$-module.

The section $S_i/T_{i-1}$ splits up into a direct sum of
 modules
$b_{2m}^{\phi} d /T_{i-1}$ for $d \in \BD_{2m}^\phi[-i]$.
(Note that $b_{2m}^{\phi} d /T_{i-1} = b_{2m}^{\phi} d' /T_{i-1} $ 
if and only if $\bra{d}=\bra{d'}$.)
The action of the algebra on the {\em right} gives maps between these
 modules for different choices of $d$. These maps are invertible by 
Proposition~\ref{struc1}, so the summands are isomorphic.
Note that there exists at least one $d$ in each 
 $\BD_{2m}^\phi[-i]$ that is idempotent up to a unit,
thus each summand is projective. 

The argument for $S=T_i$ is very similar.
\Qed

\begin{co}
If 
$\delta, \delta_L, \delta_R, \kappa_L, \kappa_R, \kappa_{LR} $ 
are units
then the modules $ \mathcal{S}_{2m}(l) $
are the  standard modules of $b_{2m}^\phi$.
\end{co}


By Proposition~\ref{rho hom} the globalisation functor $G$ is 
$$G: b_{2m-2}^{\phi}
(\delta, \delta_L,\delta_R, \kappa_L, \kappa_R, \kappa_{LR}) 
\to b_{2m}^{\phi}
(\delta, \kappa_L,\delta_R, \delta_L, \kappa_R, \kappa_{LR}) 
$$
We may `dually' define another globalisation functor using the right
hand blob 
(in Proposition~\ref{rho hom})
rather than the left hand one. We get a functor
$$G': b_{2m-2}^{\phi}
(\delta, \delta_L,\delta_R, \kappa_L, \kappa_R, \kappa_{LR}) 
\to b_{2m}^{\phi}
(\delta, \delta_L,\kappa_R, \kappa_L, \delta_R, \kappa_{LR}).
$$

It is clear
 that $G \circ G' = G' \circ G$.
We thus get three functors from
$b_{2m-4}^{\phi}$ to 
 $b_{2m}^{\phi}$ 
(as meta-categories, i.e ignoring the swapping of parameters). 

\prl(old 824)
If $\delta, \delta_L, \delta_R, \kappa_L, \kappa_R, \kappa_{LR} $ 
are units
then
$$ G(\mathcal{S}_{2m-2} (l)) \cong \mathcal{S}_{2m} (-l)$$
$$ G'(\mathcal{S}_{2m-2} (l)) \cong \mathcal{S}_{2m} (l)$$
\end{pr}
{\em Proof:}
Since $b_{2m}^\phi$ is quasi-hereditary under the assumption on the
parameters, globalising takes standard modules to standard ones 
\cite[Proposition 4]{MartinRyom02}), 
and we
need only determine which one. Globalising does not change the number
of propagating lines. The colour of the inner region changes for
$G$ but stays the same for $G'$, hence the change in sign for $G$ but
not for $G'$.
\Qed

Note that for these algebras 
$b^{\phi}_{2m}$ 
we have shown that there is a set of
prestandard modules which are in fact standard.


\subsection{Prestandard modules by $G_e$: Low rank examples}
\newcommand{\Si}{S}
\newcommand{\Sit}[2]{S_{#2}(#1)}

Let us begin the recursion, 
implicit in the proof of Theorem~\ref{index wangy},
to construct prestandard modules using $G$:
\newline
We can apply $G$ equally to simple modules directly, or to the regular
representation (and hence to diagrams), since the regular
representation is a direct sum of projective representations, each
with an appropriate unique simple module in its head. 
Firstly, 
$b^{\phi}_0$ is spanned by the empty diagram, which is thus also a
basis for the unique simple left module, $\Sit{0}{0}$.
Applying $G$ to this will give $\Sit{0}{2}$. 
Now by virtue of Proposition~\ref{basis thang}
the image of this under $\glob$ 
is given by the image under $\rho'$, 
which is 
\[
\raisebox{-.542cm}{\includegraphics{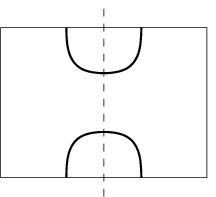}} 
\qquad \in b^{\phi}_2
\] 
In this case the multiplication map is a bijection, and the above element
generates the left prestandard module $\Sit{0}{2}$ with basis 
consisting of this and one other diagram:
\eql(bas1)
\includegraphics{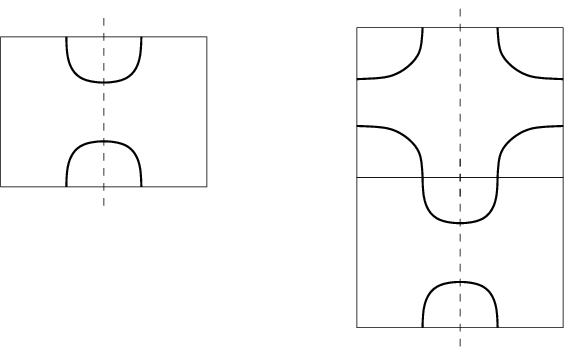}
\eq
(note that this set is spanning by the $\kappa'$ relation). 
Observe that this basis has no intrinsic dependence on the parameters.


Consider the module morphism between ideals given by $m \mapsto mf$ here.
We are considering the generic case
(so that $\kappa'$ is invertible) so 
\eql(bas2)
\includegraphics{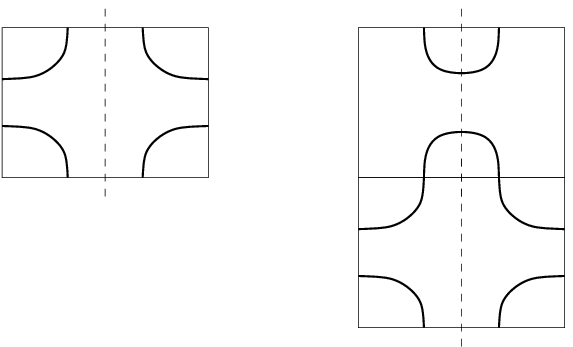}
\eq
span the isomorphic image module to that above 
(we will touch on the chiral case $RLR \not\propto R$ elsewhere). 
Note that  the
elements in (\ref{bas1}),(\ref{bas2}) span $ b^{\phi}_2 e b^{\phi}_2$.
By proposition~\ref{grebo}
the simple module missing from this construction ($\Sit{-1}{2}$) 
may be constructed as $b^{\phi}_2 / b^{\phi}_2 e b^{\phi}_2$.
Thus $\Sit{-1}{2}$ has basis
\eql(bas3)
\includegraphics{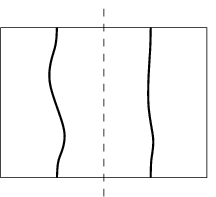}
\eq
where the action of the algebra is algebra multiplication modulo the
elements in (\ref{bas1}),(\ref{bas2}). 
Note that $| B^{\phi 0}_2| = 5$ so this is a complete decomposition of
the regular module. 


Applying $\glob$ again we determine the structure of $b^{\phi}_4$. 
The image of $b^{\phi}_2$ is as follows. 
The image of the basis elements for $\Sit{0}{2}$ in (\ref{bas1}) is the
first two elements of:
\eql(bas4)
\includegraphics{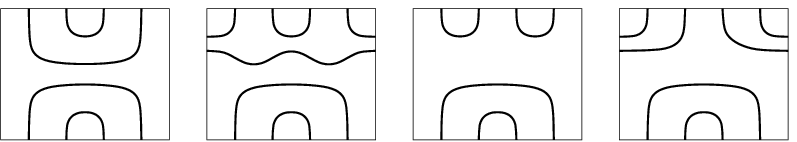}
\eq
(the other two are generated from these by the algebra action);
thus these objects are a basis for $\Sit{0}{4}$. 
Another basis for this module is obtained as the image of the elements
in (\ref{bas2}) (and two further elements generated from these):
\[
\includegraphics{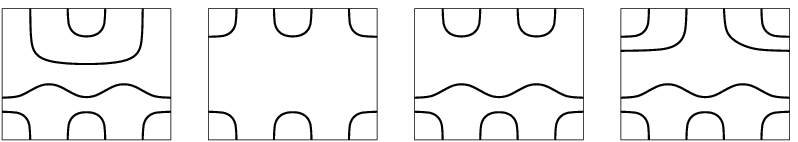}
\]
It is left as an exercise to write down two more sets of four elements
giving bases for isomorphic modules. 
The object in  (\ref{bas3}) is not strictly an element of the algebra
(because of the quotient). 
Thus we cannot apply $\rho'$ to it directly. 
However we can consider $\Sit{-1}{2}$ as 
$b^{\phi}_2 / b^{\phi}_2 e b^{\phi}_2$ and apply $\rho'$ to
$b^{\phi}_2$. 
Applying $\rho'$ to (\ref{bas3})
and to the quotienting module spanned by (\ref{bas1}),(\ref{bas2})
we have a basis for $\Sit{-1}{4}$: 
\[
\includegraphics{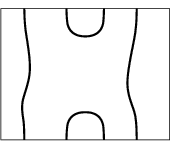}
\]

The two missing simple modules are in 
$b^{\phi}_4 / b^{\phi}_4 e b^{\phi}_4$.
Discarding all the diagrams in $ b^{\phi}_4 e b^{\phi}_4$ 
(constructed above) these are 
given by 
\[ \Sit{-2}{4} = k \; 
\raisebox{-.542cm}{\includegraphics{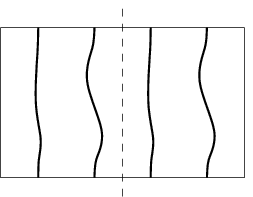}}
\]
and
\[ \Sit{+1}{4} = k \; 
\raisebox{-.542cm}{\includegraphics{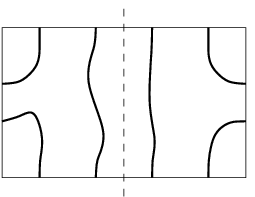}}
\]
This completes the arrangement of the basis elements for the 
left regular representation. 
We have total rank $4^2 +1+1+1$.

The $\rho'$-image of the basis elements for $\Sit{0}{4}$ in (\ref{bas4}) is the
first four elements of:
\eql(bas44)
\includegraphics{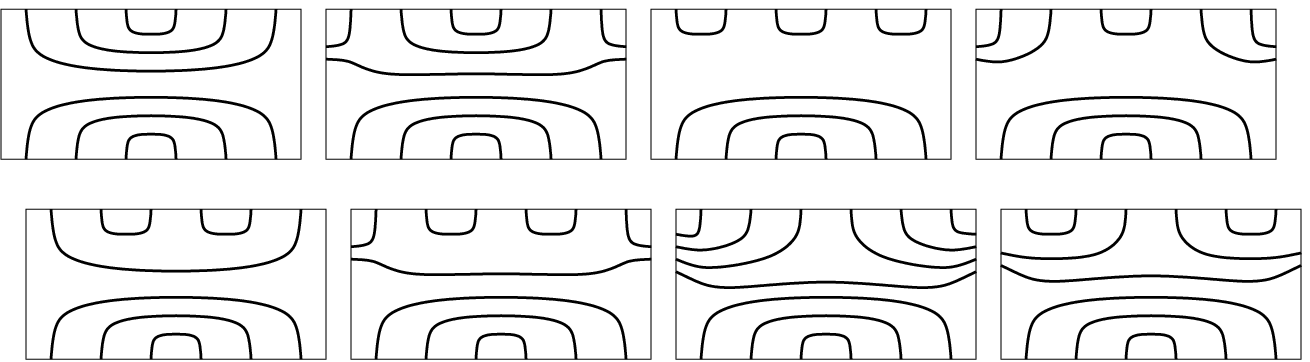}
\eq
(the other four are generated from these by the algebra action);
thus these objects are a basis for $\Sit{0}{6}$. 


The general pattern 
can now be given. 

\subsection{Top-half combinatorics}

Suppose that we view the right half of the fundamental region of a
diagram,
and concentrate for the moment on the $m$ vertices in the northern
edge in this interval
--- as it were, the upper right-hand {\em quarter} of the diagram. 
In basis elements with no propagating lines, the line from each of
these vertices descends (initially) and turns either to Left or to
Right. Considering each such line in turn, starting from the left,
say, we may construct (the upper half of) a diagram by choosing the
direction of these lines. Each direction may be chosen freely, in turn,
irrespective of the direction of previously chosen lines: we can
always choose Right since the direction of lines {\em to the right}
have yet to be chosen; we can always choose Left since either there is
a path to the left edge (the existing choices make Right-Left pairs,
i.e. cups, possibly nested, possibly together with some additional
Left directed lines); or there is a preceding Right not in a
Right-Left pair, which can then form a Right-Left pair with this new
Left. 
Some examples of these quarter diagrams 
are as follows (the shorthand on the left in each
example is L for line-turn-left; R for line-turn-right):
\eql(table key 0)
\includegraphics{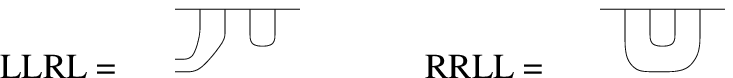}
\eq
In consequence the rank of the prestandard module $\Sit{0}{2m}$ is 
\[
| \ket{\BD^{\phi}_{2m}[0] } | \; = \; 2^{m}
\] 
This is because the composites of these elements with any element of 
$\bra{\BD^{\phi}_{2m}[0]}$ produces a basis for   $\Sit{0}{2m}$. 


Altogether the tower of bases starts as shown in figure~\ref{notbrat}, 
where $o$ denotes a propagating line.
The key for the shorthand used in this table is indicated in 
(\ref{table key 0}) above, and by 
\eql(table key)
\includegraphics{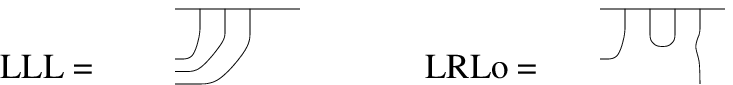}
\eq
\newcommand{\vect}[1]%
{\! \left\{ \begin{array}{c} #1 \end{array} \right\} \!}
\newcommand{\psign}{}
\newcommand{\msign}{-}
\begin{figure}
\begin{turn}{90} \vbox{
$$ \begin{array}{c|ccccccccc}
m_= &l=\psign 3 &\psign 2 &\psign 1 & {0} & \msign 1 &\msign 2 &\msign
3 
&\msign 4\\ 
\hline
0&&&&\vect{\emptyset} \\
1&&&&\vect{L\\ R}& \vect{o} \\
2&&&\vect{oR}&\vect{LL\\ LR\\ RL\\ RR}&\vect{Lo}&\vect{oo} \\
3&&\vect{ooR}&\vect{LoR}&
\vect{LLL\\ LLR\\ LRL\\ LRR\\ RLL\\ RLR\\ RRL\\ RRR}&
\vect{LLo\\ RLo\\ oRL\\ oRR}&\vect{Loo}&\vect{ooo} \\
4&\vect{oooR}&\vect{LooR}&
\vect{oRRR\\ oRLR\\ oRRL\\ RLoR\\ LLoR} &
\vect{LLLL\\ LLLR\\ LLRL\\ LLRR\\ LRLL\\ \ldots\\ RRRR} & 
\vect{LLLo\\ LRLo\\ RLLo\\ LoRL\\ LoRR} & 
\vect{LLoo\\ RLoo\\ oRLo\\ ooRL\\ ooRR} &
\vect{Looo} &
\vect{oooo}
\end{array} $$}
\end{turn}
\caption{\label{notbrat} Table of standard $\Sit{l}{2m}$ bases
up to $m=4$.
NB, this $LR$ shorthand should NOT be confused with the nonabelian
ring elements which live on strings 
--- see main text for key.}
\end{figure}
\newcommand{\ur}{\mbox{ur}_1}%
\newcommand{\uro}{\mbox{ur}_0}%

Let us write 
$\ur(d)$ for the number of lines passing out of the fundamental region
of a (half-)diagram through the 1-wall on the right,
not counting those lines which also pass out on the left
(`equatorial' lines, as it were).
Thus for example the 4-th diagram in  (\ref{bas44}) has 
$\ur(d) = 2$. 
(We will also later use $\uro(d)$ for the number of lines crossing the
0-wall of a (half-)diagram: 
the left hand edge, then, of a quarter diagram.) 

Note that the set of half-diagrams with given $m$ and $l=+x$ 
($x >0$) is that set with $x$ propagating lines and 
$\ur(d) \equiv 1$ mod~2,
while the set of half-diagrams with given $m$ and $l=-x$ 
($x >0$) is that set with $x$ propagating lines and 
$\ur(d) \equiv 0$ mod~2.

\subsection{Restriction of prestandards to blob algebra standards}\label{wa}
\newcommand{\Deltab}[2]{\Delta^{b}_{#2}(#1)}
The representation of $b_3$ induced on the basis for $\Sit{0}{6}$ in
(\ref{bas44}) is as follows (we use the isomorphism with $b_6'$
--- and by virtue of proposition~\ref{blob iso} we use $b_m$ and $b_{2m}'$
interchangeably in this section):
\[
R_0( \; \raisebox{-.1542cm}{\includegraphics{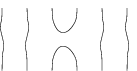}} \; )
=\mat{llllllll} 
\delta_L \\
         & \delta_L \\
        & & \delta_L \\
        & & & \delta_L \\
1 &&&& 0 \\
 &1&&&   &0 \\
&&& \kappa &&& 0 \\
&& \kappa &&&&& 0 \tam
\]
\[
R_0( \; \raisebox{-.1542cm}{\includegraphics{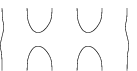}} \; )
=\mat{llllllll} 
0 &&&& \kappa_L \\
         & 0 &&&& \kappa_L \\
        & & 0 && 1\\
        & & & 0 && 1\\
 &&&& \delta \\
 &&&&   &\delta \\
&&&  &&\kappa_R& 0 \\
&&  &&&1&& 0 \tam
\]
\[
R_0( \; \raisebox{-.1542cm}{\includegraphics{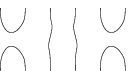}} \; )
=\mat{llllllll} 
0 && \delta_L \\
         & 0 & \kappa \\
        & & \delta && \\
        & & \delta_R& 0 && \\
 &&1&& 0 \\
 &&&&   &0&&1 \\
&&&  &&& 0 & \delta_R \\
&&  &&&&& \delta \tam
\]
(all unmarked entries zero).
Note that basis elements in positions 1,3,5 span a $b_3$-submodule
isomorphic to $\Deltab{1}{3}$
(where $\Deltab{l}{m} = \Delta_{m}(l)$
of \cite{MartinWoodcock2000}). 
The quotient by this has 
$b_3$-submodule spanned by elements 2,6,8, 
isomorphic to $\Deltab{-1}{3}$. The quotient by this has 
$b_3$-submodule spanned by element 4, 
isomorphic to $\Deltab{-3}{3}$.
The remaining quotient is isomorphic to $\Deltab{3}{3}$.


We can see this decomposition directly by looking at (\ref{bas44}).
Note for general $\Sit{l}{2m}$ that the number 
$\ur(d)$ 
cannot be increased by acting on a half-diagram
by an element of $b_m$ (i.e. of $b_{2m}'$). Thus
\prl(ur1)
The restriction $\Res(b_{2m}',b_{2m}^{\phi},) \Sit{l}{2m}$ 
has submodule structure filtered
by $\ur$. In particular 
$\Res(b_{2m}',b_{2m}^{\phi},) \Sit{0}{2m}$ 
has $m+1$ sections, since all the $\ur$ values are realised from 0 to
$m$.
Meanwhile  $\Res(b_{2m}',b_{2m}^{\phi},) \Sit{l}{2m}$
with $l=\pm x$ ($x>0$)
has  $\ur$ values from, and hence
sections indexed by: \\
$\{ m-x-1, m-x-3, \ldots, 0 \}$ if $m-x$ odd and $l<0$; \\
$\{ m-x, m-x-2, \ldots, 1 \}$ if $m-x$ odd and $l>0$; \\
$\{ m-x-1, m-x-3, \ldots, 1 \}$ if $m-x$ even and $l>0$; \\
$\{ m-x, m-x-2, \ldots, 0 \}$ if $m-x$ even and $l<0$.
\end{pr}
{\em Proof:} Only the $l \neq 0$ cases require further explanation. 
Here, since there are propagating lines, there can be no equatorial line 
(passing from one reflection wall to the other), so the
lines contributing to $\ur$ all start from the northern edge of the diagram.
Consider the lines passing out of the $m$ vertices on the northern
edge of the diagram. 
The number of propagating lines is fixed at $x$, and the total number
$\ur + \uro$ of lines passing to the reflection walls 
is of definite parity, since all other lines return to the north edge
and hence contribute to $m$ in pairs. But the parity of $\uro$ is
fixed by the sign of $l$, so the parity of $\ur$ is also fixed.
It is routine to check the extremal numbers. 
\Qed


Consider  the $r$-th $\ur$-section of   
$\Res(b_{2m}',b_{2m}^{\phi},) \Sit{l}{2m}$ ($l=\pm x$ ($x>0$)).
As already noted,
$\uro$ is of definite parity in this section.
If $\uro$ is even,
there is an injective map from the basis elements in this section 
into
the basis elements of a $b_{2m}'$ standard module 
 $\Deltab{(x+r)}{}$, 
obtained by 
deforming the ends of the $r$ lines that pass out through the 1-wall
until they pass out through the bottom of the diagram (i.e. become
propagating lines). 
It is easy to see that this extends to a module morphism. 
Since the section is a blob module the map must be
onto and hence a bijection.
There is a similar morphism for $\uro$ odd. 
\prl(ur1 2)
In the $\ur$-sections of   
$\Res(b_{2m}',b_{2m}^{\phi},) \Sit{l}{2m}$ ($l=\pm x$ ($x>0$))
an $\ur$-line acts like a propagating line. 
If $m-x$ odd and $l<0$, or $m-x$ even and $l>0$, then the inner
region is black so the first $\uro$ line  also 
acts as a propagating (and blobbed) line.
Taking into account the $x$ lines that are already propagating, 
the section with $\ur = r$ is thus isomorphic to a 
$b_{2m}'$ standard module of form $\Deltab{-(x+r+1)}{}$
if one of the `black' conditions above is satisfied; 
and of form $\Deltab{x+r}{}$ otherwise. 
(NB, the sign on the weight here does not affect the dimension of the
module.) 
A similar statement holds for $l=0$, so that 
$\Res(b_{2m}',b_{2m}^{\phi},) \Sit{0}{2m}$
is a sum of one copy of each blob standard. 
\Qed
\end{pr}


Recall that  the dimension of $\Sit{0}{2m}$ is $2^m$.

\begin{co} \label{wangy}
Consider the prestandard $\Sit{l}{2m}$, where $l = \pm x$ and $x > 0$.
Define the integer $\epsilon$ by

  $      \epsilon = 1 $  if $m-x$ odd,

$        \epsilon = 2 $  if $m-x$ even and $l > 0$,

$        \epsilon = 0 $  if $m-x$ even and $l < 0$,
\\
and let $k = (m - (x+\epsilon))/2$.  
Then the dimension of $S_l(2m)$ is given by
        $\sum_{i = 0}^k {m \choose i}$.

\end{co}
{\em Proof:}
From \cite{MartinSaleur94a}, the dimension of the $b'_{2m}$-standard module
$\Deltab{\pm c}{}$ is given by ${m \choose {(m - c)/2}}$.  The result now
follows from Proposition~\ref{ur1 2}, summing over the $r$-values specified in
Proposition~\ref{ur1}. 
\Qed
\\
Thus we have determined the complete generic representation theory.

Note that the globalisation and localisation functors act 
in a natural way on blob
categories as well as \achiralb\ categories. 
We did not need this fact here, but it is useful in computing 
non-generic representation theory.


\section{Discussion}
Having determined the generic representation theory,
and set up the homological machinery for analysing the 
exceptional (non-semisimple) cases,
in our next paper we will turn to computing the 
representation theory of the exceptional cases. 
We conclude here with a brief 
introduction to
this problem.


\newtheorem{proof}{Proof}

With the Temperley--Lieb and blob algebras, the symplectic blob
algebra (or isomorphically, $b^{\phi}_{2m}$) 
 belongs to an intriguing class of Hecke algebra quotients. 
The first two
have representation theories beautifully and efficiently described in
alcove geometrical language, where the precise geometry is determined,
in the non-semisimple cases, 
by the parameters of the algebra. 
In these first two  algebras the parameterisation appropriate to
reveal this structure is not that in which the algebras were first
described. Rather, it was discovered during efforts to put the low
rank
data on non-semisimple manifolds in parameter space in a coherent
format \cite{MartinSaleur94a}. 
The determination of the representation theory of $b^{\phi}_{2m}$ 
in the non-semisimple cases is the next important problem in the 
programme initiated in this paper. 
We therefore conclude the paper with one result on the
characterisation of non-semisimple manifolds. 
The crucial point is that because of the globalisation map,
this is derived from a 
low rank result, which then globalises to all levels in the tower. 


Throughout this section we will assume that all the parameters are
units, so $b_{2m}^{\phi}$ is quasihereditary.


Set $L_{2m}(l)$ to be the irreducible head of the standard module
$ \sS_{2m}(l)$ .
At any point in parameter space for which the 
 Gram determinant $\Deltag_{2m}(l) $
evaluates to zero, 
there is a proper submodule of $\sS_{2m}(l)$
and so $\sS_{2m}(l)$ is not simple. In this case using the fact that
$b_{2m}^\phi$ is quasihereditary we can find a non-zero map
$\sS_{2m}(j) \to \sS_{2m}(l)$ for some $j \ne l$.
Once we have found a non-zero map we can then globalise it to larger
$m$, using functor $G$ 
and 
Proposition~\ref{old 824}.

Define polynomials 
in the six parameters $\{\delta,
\delta_L, \delta_R, \kappa_L,   
\kappa_R, \kappa_{LR}\}$:
\begin{eqnarray*}
K_{0}\! \!\! &=& \!\! \kappa_{LR} \\
K_{1}\! \!\! &=& \!\! \delta_L \delta_R - \kappa_{LR} \\
K_{2}\! \!\! &=& \!\! \kappa_{LR} - \delta_L \kappa_R - \kappa_L \delta_R
               +\delta\delta_L\delta_R  \\
K_{3}\! \!\! &=& \!\! \delta^2\delta_L\delta_R-\delta\delta_L\kappa_R
               -\delta\delta_R\kappa_L-\delta_l\delta_R+\kappa_L\kappa_R \\
K_{1,3}\! \!\! &=& 
\delta^2\delta_L\delta_R-\delta\delta_L\kappa_R
-\delta\delta_R\kappa_L+\kappa_L\kappa_R-\kappa_{LR}
\end{eqnarray*}
and commuting operators $\Phi$ and $\Psi$ on the space of six-parameter
polynomials which swap the second and fourth, \resp\ third and fifth,
parameters. 
The non-trivial Gram determinants for $ b_6^\phi$ are:
$$
\Deltag_6(-1) =
\kappa_L \kappa_R K_3
$$
\begin{eqnarray*}
\Deltag_6(0) \!\!\!&=&\!\!
\kappa_{LR}^4 K_1^4 \Psi\Phi(K_1) 
\Psi(K_2) \Phi(K_2) K_{1,3}
\end{eqnarray*}


We may deduce maps:
$ S_6 (-1) \hookrightarrow S_6(0)$ for
 $K_1=0$,
$ S_6 (1) \hookrightarrow S_6(0)$ for
$\Psi\Phi(K_1) =0$,
$ S_6 (2) \hookrightarrow S_6(0)$ 
for 
$\Phi(K_2)=0$
and
$ S_6 (-2) \hookrightarrow S_6(0)$ for
$\Psi(K_2)=0$.
We may deduce that the only possible
non-zero map to $S_6(-1)$ is
$S_6(-3) \hookrightarrow S_6(-1)$ and this therefore must occur when
$K_3=0$.

We also get a non-zero map $S_6(-3) \hookrightarrow S_6(0)$ when
$K_{1,3} =0$.
Thus when both 
$K_1=0$
and $K_3=0$ 
we get two copies of $S_{6}(-3)$ in the socle of $S_6(0)$.

\begin{pr}
$b_{2m}^\phi(\delta, \delta_L,\delta_R, \kappa_L, \kappa_R,
\kappa_{LR}) $
is not semisimple when

\begin{tabular}{l}  \hspace{-.3in}(1)
$K_3=0$
and $m$ is odd and $m \ge 3$ \\
$\Phi (K_3)=0$
and $m$ is even and $m \ge 4$ \\
$\Psi (K_3)=0$
and $m$ is even and $m \ge 4$ \\
$\Phi\Psi (K_3)=0$
and $m$ is odd and $m \ge 5$. \\
\end{tabular}

\begin{tabular}{l} \hspace{-.3in}(2)
$K_{1,3}=0$
and $m$ is odd and $m \ge 3$ \\
$\Phi(K_{1,3})=0$
and $m$ is even and $m \ge 4$ \\
$\Psi(K_{1,3})=0$
and $m$ is even and $m \ge 4$ \\
$\Psi\Phi(K_{1,3})=0$
and $m$ is odd and $m \ge 5$. \\

\end{tabular}
\end{pr}
Cases (1) are proved by globalising the map $ S_6 (-3) \to S_6(-1)$ 
for
 $b_{6}^{\phi}(\delta, \delta_L,\delta_R, \kappa_L, \kappa_R,
\kappa_{LR}) $ with
$K_3=0$; 
cases (2) are proved by globalising the map $S_6(-3) \to S_6(0)$ 
for
 $b_{6}^{\phi}(\delta, \delta_L,\delta_R, \kappa_L, \kappa_R,
\kappa_{LR}) $ with
$K_{1,3} =0$.



\vspace{1.398in}

\appendix
\medskip\noindent
{\LARGE\bf Appendix} \vspace{-.1in}
\stuffb

\stuffbb

\bibliographystyle{amsplain}
\bibliography{new31}

\providecommand{\bysame}{\leavevmode\hbox to3em{\hrulefill}\thinspace}
\providecommand{\MR}{\relax\ifhmode\unskip\space\fi MR }
\providecommand{\MRhref}[2]{%
  \href{http://www.ams.org/mathscinet-getitem?mr=#1}{#2}
}
\providecommand{\href}[2]{#2}
\begin{thebibliography}{10}

\bibitem{AlvarezMartin05}
M~Alvarez and P~P Martin, \emph{Towards an algebra of 2-manifold surgery},
  preprint (2005).

\bibitem{Armstrong79}
M~A Armstrong, \emph{Basic topology}, McGraw Hill, 1979.

\bibitem{Baxter82}
R~J Baxter, \emph{Exactly solved models in statistical mechanics}, Academic
  Press, New York, 1982.

\bibitem{BeilinsonBernsteinDeligne81}
A.A. Beilinson, J.N. Bernstein, and P.~Deligne, \emph{Faisceaux pervers},
  Analyse et topologie sur les espaces singuliers, vol.~1, Ast{\'e}risque, no.
  100, 1981.

\bibitem{Bergman}
G~M Bergman, \emph{The diamond lemma for ring theory}, Adv. Math. \textbf{29}
  (1978), 178--218.

\bibitem{Bloss04}
M~Bloss, \emph{$g$-colored partition algebras as centralizer algebras of wreath
  products}, Journal of Algebra \textbf{265} (2003), 690--710.

\bibitem{Brauer37}
R~Brauer, \emph{On algebras which are connected with the semi--simple
  continuous groups}, Annals of Mathematics \textbf{38} (1937), 854--872.

\bibitem{Cherednik91}
I~Cherednik, \emph{A unification of {K}nizhnik--{Z}amolodchikov and {Dunkl}
  operators via affine {H}ecke algebras}, Invent Math \textbf{106} (1991),
  411--431.

\bibitem{Chevalley56}
C~Chevalley, \emph{Fundamental concepts of algebra}, Academic Press, 1956.

\bibitem{ClineParshallScott88}
E~Cline, B~Parshall, and L~Scott, \emph{Finite-dimensional algebras and highest
  weight categories}, J. reine angew. Math. \textbf{391} (1988), 85--99.

\bibitem{CoxGrahamMartin03}
A~G Cox, J~J Graham, and P~P Martin, \emph{The blob algebra in positive
  characteristic}, J Algebra \textbf{266} (2003), 584--635.

\bibitem{CoxMartinParkerXi03}
A~G Cox, P~P Martin, A~E Parker, and C~C Xi, \emph{Representation theory of
  towers of recollement: theory, notes and examples}, preprint to appear in
  {\em J Algebra} (math.RT/0411395), 2003.

\bibitem{CurtisReiner62}
C~W Curtis and I~Reiner, \emph{Representation theory of finite groups and
  associative algebras}, Wiley Interscience, New York, 1962.

\bibitem{DeGier02}
J~de~Gier, \emph{Loops, matchings and alternating-sign matrices}, 14th
  International Conference on Formal Power Series and Algebraic Combinatorics
  (Melbourne 2002), math.CO/0211285 (2002).

\bibitem{ErdmannGreen99}
K.~Erdmann and R.M. Green, \emph{On representations of affine
  {T}emperley--{L}ieb algebras, {II}}, Pacific J. Math. \textbf{191} (1999),
  243--273.

\bibitem{FanGreen99}
C~K Fan and R~M Green, \emph{On the affine {T}emperley--{L}ieb algebras}, J.
  LMS \textbf{60} (1999), 366--380.

\bibitem{GoodmanDelaharpeJones89}
F~M Goodman, P~de~la Harpe, and V~F~R Jones, \emph{Coxeter graphs and towers of
  algebras}, Math Sci Research Inst Publications 14, Springer--Verlag, Berlin,
  1989.

\bibitem{Graham95}
J~J Graham, \emph{Modular representations of {H}ecke algebras and related
  algebras}, Ph.D. thesis, Mathematics, University of Sydney, 1995.

\bibitem{GrahamLehrer96}
J.~J. Graham and G.~I. Lehrer, \emph{Cellular algebras}, Invent. Math.
  \textbf{123} (1996), 1--34.

\bibitem{GrahamLehrer98}
\bysame, \emph{The representation theory of affine {Temperley-Lieb} algebras},
  L'Enseignement Math\'{e}matique \textbf{44} (1998), 173--218.

\bibitem{GrahamLehrer03}
\bysame, \emph{Diagram algebras, {H}ecke algebras and decomposition numbers at
  roots of unity}, Annales Scientifiques de l'\'Ecole Normale Sup\'erieure
  \textbf{36} (2003), no.~4, 479--524.

\bibitem{Green80}
J~A Green, \emph{Polynomial representations of ${GL}_n$}, Springer-Verlag,
  Berlin, 1980.

\bibitem{Green98}
R.~M. Green, \emph{Generalized {Temperley--Lieb} algebras and decorated
  tangles}, J. Knot Theory Ramifications \textbf{7} (1998), 155--171.

\bibitem{Green98b}
R~M Green, \emph{On representations of affine {T}emperley--{L}ieb algebras},
  Canad Math Soc Conference Proc \textbf{24} (1998), 245--261.

\bibitem{Green03}
\bysame, \emph{On planar algebras arising from hypergroups}, J. Algebra
  \textbf{263} (2003), 126--150.

\bibitem{HaringOldenburg99}
R.~H\"aring-Oldenburg, \emph{The reduced {B}irman-{W}enzl algebra of {C}oxeter
  type {B}}, J Algebra \textbf{213} (1999), 437--466.

\bibitem{Humphreys90}
J~E Humphreys, \emph{Reflection groups and {C}oxeter groups}, Cambridge
  University Press, 1990.

\bibitem{JamesKerber81}
G~D James and A~Kerber, \emph{The representation theory of the symmetric
  group}, Addison-Wesley, London, 1981.

\bibitem{JonesPlanar}
V~F~R Jones, \emph{Planar algebras, {I}}, unpublished.

\bibitem{Jones94b}
\bysame, \emph{A quotient of the affine {H}ecke algebra in the {B}rauer
  algebra}, L'Enseignement Math\'ematique \textbf{40} (1994), 313--344.

\bibitem{Jucys74}
A-A~A Jucys, \emph{Symmetric polynomials and the centre of the symmetric group
  ring}, Rep Math Phys \textbf{5} (1974), 107--112.

\bibitem{KoenigXi99b}
S~Koenig and C~C Xi, \emph{When is a cellular algebra quasi-hereditary?},
  Mathematische Annalen \textbf{315} (1999), 281--293.

\bibitem{Levy91}
D~Levy, \emph{Algebraic structure of translation-invariant spin-1/2 {XXZ} and
  q-{P}otts quantum chains}, Phys Rev Lett \textbf{67} (1991), 1971--1974.

\bibitem{Lusztig89C}
G~Lusztig, \emph{Affine {H}ecke algebras and their graded version}, J Am Math
  Soc \textbf{2:3} (1989), 599--685.

\bibitem{Martin91}
P~P Martin, \emph{Potts models and related problems in statistical mechanics},
  World Scientific, Singapore, 1991.

\bibitem{Martin94}
\bysame, \emph{Temperley--{L}ieb algebras for non--planar statistical mechanics
  --- the partition algebra construction}, Journal of Knot Theory and its
  Ramifications \textbf{3} (1994), no.~1, 51--82.

\bibitem{Martin95pres}
\bysame, \emph{Pascal's triangle and word bases for blob algebra ideals},
  preprint (1995).

\bibitem{Martin96}
\bysame, \emph{The structure of the partition algebras}, J Algebra \textbf{183}
  (1996), 319--358.

\bibitem{MartinRyom02}
P~P Martin and S~Ryom-Hansen, \emph{Virtual algebraic {Lie} theory: Tilting
  modules and {Ringel} duals for blob algebras}, Proc LMS \textbf{89} (2004),
  655--675, (math.RT/0210063).

\bibitem{MartinSaleur93}
P~P Martin and H~Saleur, \emph{On an algebraic approach to higher dimensional
  statistical mechanics}, Commun. Math. Phys. \textbf{158} (1993), 155--190.

\bibitem{MartinSaleur94a}
\bysame, \emph{The blob algebra and the periodic {T}emperley--{L}ieb algebra},
  Lett. Math. Phys. \textbf{30} (1994), 189--206.

\bibitem{MartinWoodcock2000}
P~P Martin and D~Woodcock, \emph{On the structure of the blob algebra}, J
  Algebra \textbf{225} (2000), 957--988.

\bibitem{MartinWoodcock03}
\bysame, \emph{Generalized blob algebras and alcove geometry}, LMS J Comput
  Math \textbf{6} (2003), 249--296, (math.RT/0205263).

\bibitem{MartinWoodcockLevy00}
P~P Martin, D~Woodcock, and D~Levy, \emph{A diagrammatic approach to {H}ecke
  algebras of the reflection equation}, J Phys A \textbf{33} (2000),
  1265--1296.

\bibitem{Moise77}
E~E Moise, \emph{Geometric topology in dimensions 2 and 3}, Graduate Texts in
  Mathematics 47, Springer-Verlag, New York, 1977.

\bibitem{Murphy81}
G~E Murphy, \emph{A new construction of {Y}oung's seminormal representation of
  the symmetric group}, J Algebra \textbf{69} (1981), 287--291.

\bibitem{Nazarov96}
M~Nazarov, \emph{Young's orthogonal form for {B}rauer's centralizer algebra}, J
  Algebra \textbf{182} (1996), 664--693.

\bibitem{NicholsRittenbergdeGier05}
A~Nichols, V~Rittenberg, and J~de~Gier, \emph{One-boundary {Temperley--Lieb}
  algebras in the {XXZ} and loop models}, J Stat (2005, to appear),
  cond-mat/0411512.

\bibitem{OrellanaRam01}
R~Orellana and A~Ram, \emph{Affine braids, {M}arkov traces and the category
  {O}}, Newton Institute preprint NI01032-SFM, 2001.

\bibitem{PasquierSaleur90}
V~Pasquier and H~Saleur, \emph{Common structures between finite systems and
  conformal field theories through quantum groups}, Nucl Phys B \textbf{330}
  (1990), 523.

\bibitem{RuiXi04}
H~Rui and C~C Xi, \emph{The representation theory of cyclotomic
  {T}emperley--{L}ieb algebras}, Comment Math Helv \textbf{79} (2004),
  427--450.

\bibitem{Sahi99}
Siddhartha Sahi, \emph{Nonsymmetric {K}oornwinder polynomials and duality},
  Annals of Mathematics \textbf{150} (1999), 267--282.

\bibitem{TemperleyLieb71}
H~N~V Temperley and E~H Lieb, \emph{Relations between percolation and colouring
  problems and other graph theoretical problems associated with regular planar
  lattices: some exact results for the percolation problem}, Proceedings of the
  Royal Society A \textbf{322} (1971), 251--280.

\bibitem{tomDieck94}
T~tom Dieck, \emph{Symmetrische {B}r\"ucken und {K}notentheorie zu den
  {D}ynkin--diagrammen vom type {B}}, J Reine angew Math \textbf{451} (1994),
  71--88.

\bibitem{Weyl46}
H~Weyl, \emph{Classical groups}, Princeton, Princeton, 1946, ch.III--V.

\end{thebibliography}
\end{document}